\documentclass[11pt,a4paper,reqno]{article}
 \usepackage[letterpaper,margin=0.80in,bottom=0.90in,top=0.85in]{geometry}

\usepackage[latin,english]{babel}
\usepackage{comment}
\usepackage[utf8]{inputenc}
\usepackage[T1]{fontenc}
\usepackage{eufrak, color, mathrsfs,dsfont,geometry, bbm, dutchcal}
\usepackage{amsmath, amssymb, amsthm, pdfsync, hyperref}
\usepackage{mathtools, mathabx, cases, braket, upgreek}
\usepackage{babelbib}
\usepackage{math}
\mathtoolsset{showonlyrefs}

\numberwithin{equation}{section}

\theoremstyle{plain}
\newtheorem{theorem}{Theorem}[section]
\newtheorem{lemma}[theorem]{Lemma}
\newtheorem{proposition}[theorem]{Proposition}
\newtheorem{corollary}[theorem]{Corollary}
\theoremstyle{definition}
\newtheorem{definition}[theorem]{Definition}

\theoremstyle{remark}
\newtheorem{remark}[theorem]{Remark}

\renewcommand{\epsilon}{\varepsilon}
\newcommand{\Opw}{\operatorname{Op^w}}
\newcommand{\lala}{\langle\hspace{-2pt}\langle}
\newcommand{\rara}{\rangle\hspace{-2pt}\rangle}
\renewcommand{\div}{\textup{div}}
\renewcommand{\cX}{X}

\begin{document}

\title{Growth of Sobolev norms for completely resonant \\
quantum harmonic oscillators on $\R^2$}

\author{
Beatrice Langella\footnote{
	International School for Advanced Studies (SISSA), Via Bonomea 265, 34136, Trieste, Italy. \newline
	\textit{Email:} \texttt{beatrice.langella@sissa.it}, },
Alberto Maspero\footnote{
	International School for Advanced Studies (SISSA), Via Bonomea 265, 34136, Trieste, Italy. \newline
	\textit{Email:} \texttt{amaspero@sissa.it}}, \  
 Maria Teresa Rotolo\footnote{	International School for Advanced Studies (SISSA), Via Bonomea 265, 34136, Trieste, Italy. \newline 
	 \textit{Email:} \texttt{mrotolo@sissa.it}}}

\date{}
\maketitle

\begin{abstract}
We consider  time dependently perturbed  quantum harmonic oscillators in $\mathbb{R}^2$: 
$$\im \partial_t u=\frac12(-\partial_{x_1}^2-\partial_{x_2}^2 + x_1^2+x_2^2)u +V(t, x, D)u, \qquad \ x \in \mathbb{R}^2,$$
where $V(t, x, D)$ is a selfadjoint pseudodifferential operator of degree zero, $2\pi$ periodic in time.\\
We identify sufficient conditions on the principal symbol of the potential $V(t, x, D)$ that ensure existence of solutions exhibiting unbounded growth in time of their positive Sobolev norms and we show that the class of symbols satisfying such conditions is generic in the Fréchet space of classical $2\pi$- time periodic symbols of order zero. To prove our result we apply the abstract Theorem of \cite{Mas22}: the main difficulty is to find a conjugate operator $A$ for the resonant average of $V(t,x, D)$. We construct explicitly the symbol of the conjugate operator $A$, called escape function, combining techniques from microlocal analysis, dynamical systems and contact topology. 
\end{abstract}

\tableofcontents

\section{Introduction}
In this paper we construct solutions of  perturbed quantum harmonic oscillators in $\R^2$ undergoing infinite time Sobolev norm explosion. 
Precisely, we consider the equation 
\begin{equation}
    \im \partial_t u=H_0u+V(t, x, D)u, \quad x=(x_1, x_2) \in \R^2, 
\label{Quantum_Harmonic_oscillator}
\end{equation}
where $H_0:=\frac12(-\partial_{x_1}^2-\partial_{x_2}^2+x_1^2+x_2^2)$, $D:={(D_1, D_2)= }(\frac{1}{\im} \pa_{x_1}, \frac{1}{\im} \pa_{x_2})$, and $V(t, x, D)$ is a pseudodifferential operator of degree zero, selfadjoint for every $t$ and $2\pi$-periodic in time, and we prove existence of  global in time solutions fulfilling
\begin{equation}\label{s.growth}
    \lim_{t \to +\infty} \norm{u(t)}_s= + \infty \ ,  \quad s >0 \ ;
\end{equation}
here the norm $\norm{\cdot}_s$ is the one of the  Sobolev space 
\begin{equation}\label{Sob.Spaces}
    \cH^s(\R^2):=\{u \in L^2(\R^2, \C) : \|u\|_s:=\|H_0^s u\|_{L^2}<+\infty\}, \qquad  s \ \in \R \ . 
\end{equation}
Note immediately that, 
when  $V \equiv 0$, every solution of \eqref{Quantum_Harmonic_oscillator} has constant  $\cH^s$-norm for every $s$ and every time.
Hence, the existence of solutions with unbounded paths depends  on the ability of the perturbation $V(t, x, D)$ to couple together sufficiently many different Fourier modes to provoke energy cascades towards higher and  higher frequencies.
Of course, here a key role is played by the fact that the time frequency of $V$ is completely resonant with the frequencies of the harmonic oscillators.

In the recent paper \cite{Mas23}, one of the authors  has proved that, in 1-spatial dimension, generic $2\pi$-periodic pseudodifferential perturbations $V(t, x_1 , D_1)$ of the harmonic oscillator
provoke growth of Sobolev norms. 
However such 1-dimensional perturbations, when considered as perturbations on $\R^2$, are neither pseudodifferential, nor generic in the higher dimensional setting.

The goal of the present paper is to prove that, for typical bounded pseudodifferential perturbations $V(t,x,D)$ in $\R^2$, phenomena of growth of Sobolev norms occur. 
A rough statement of our result is 
\begin{theorem}\label{main.informal}
For a ``generic'' zero-order pseudodifferential operator $V(t, x, D)$, $2\pi$-periodic in time,  for any $s >0$, there exists a solution  $u(t,x) \in \cH^s$ of \eqref{Quantum_Harmonic_oscillator} such that $\|u(t)\|_s \geq C_s t^s$, for some $C_s>0$ and for all $t$ large enough.
\end{theorem}
We refer to Theorem \ref{main_theorem} below  for a precise statement. For the moment we limit ourselves to 
explain what ``generic'' means: since $V(t, x, D)$ is a pseudodifferential operator, it has a well defined symbol. Then  our result holds for symbols in an open and dense set with respect to the  Fréchet topology of symbols.
Actually we prove much more: we give an explicit sufficient condition on the symbol for its  quantization to  provoke growth of Sobolev norms.
The condition is the following: 
denote by  $v_0(t, x, \xi)$, $\xi=(\xi_1, \xi_2) \in \R^2$,   the principal symbol of $V(t,x,D)$, which we assume to be a  definitely homogeneous function of degree zero (see Definitions \ref{def:pos_hom} and \ref{def:symbol}). 
Denote by $\la v_0 \ra$  the average in time of $v_0$ along the flow $\phi^t_{h_0}$ of the classical bidimensional harmonic oscillator, i.e. 
$$
\la v_0 \ra(x,\xi) := \frac{1}{2\pi} \int_0^{2\pi} v_0(t, \phi^t_{h_0}(x, \xi)) \di t \ ;
$$
note that such symbol is time-independent. 
Then, by the  abstract result  in \cite{Mas22}, if  $\la v_0 \ra$ admits an escape function, namely a symbol $a(x,\xi)$ of order one such that, on some energy level of $\la v_0 \ra$, 
\begin{equation}\label{cond.intro}
   \{ \la v_0 \ra, a \} \geq \delta >0 \ , 
\end{equation}
the operator $V(t,x,D)$ provokes growth of Sobolev norms. 

To prove the existence of an escape function is far from trivial, so is the claim that {\em generic} symbols admit an escape function. The proof of these facts is the core of our paper, and it  relies on a delicate combination of elements from microlocal analysis, contact geometry and dynamical systems.
Here lays also the main difference with \cite{Mas23}: in the  one-dimensional setting of \cite{Mas23} both the construction of the escape function and the proof that \eqref{cond.intro} holds for generic symbols can be done in a rather direct and explicit way, which is not possible in the present 
  higher dimensional setting.

Let us briefly explain how we proceed.
The first step is to 
 extend to our framework a construction done by Colin de Verdière  \cite{Colin_de_Verdi_re_2020} for symbols on the torus  (see also \cite{CdVSR}).
We prove that the  existence of an escape function for $\la v_0\ra$   relies on   dynamical properties of the vector field 
$\tilde X_{\la v_0 \ra}$ obtained 
by restricting the Hamiltonian vector field $X_{\la v_0 \ra}$ on the energy level $\la v_0 \ra = e_0$ (which due to the homogeneity of $\la v_0 \ra$ has a cone structure) and then
projecting $X_{\la v_0 \ra}$ on the unit sphere $\S^3$ (see \eqref{Xh}). 
Then the needed dynamical property is that $\tilde X_{\la v_0 \ra}$
has Morse-Smale flow on the surface 
$Z^{\la v_0 \ra}_{e_0}:= \la v_0 \ra^{-1}(e_0) \cap \S^3$.
Actually, a more general dynamical condition is sufficient (see Proposition \ref{simple_structure_thm}), which includes Morse-Smale flows as a special case. 

Then we are left to prove that, for {\em generic} symbols $v(t,x,\xi)$, the projected vector field $\tilde X_{\la v_0 \ra}$ is Morse-Smale. 
Whereas Morse-Smale vector fields are generic among the smooth vector fields on two dimensional, closed, oriented manifolds-- a celebrated Theorem by Peixoto \cite{Peixoto1962StructuralSO}-- 
the projected vector field $\tilde X_{\la v_0 \ra}$ is rather specific, leaving the problem to ensure that a generic symbol $v$ produces a Morse-Smale vector field $\tilde X_{\la v_0 \ra}$.
To this goal, we exploit a structural property of  projected vector fields of the form $\tilde X_{h}$ when  $h$ is homogeneous of degree zero. 
Precisely, we regard $Z^{h}_{e_0}$ as a submanifold of the contact manifold $\S^3$ (with its standard contact structure inherited by the symplectic one of $\R^4$). 
The contact structure of $\S^3$ induces geometrically a natural foliation on  $Z^{h}_{e_0}$, which is called {\em characteristic foliation} (Definition \ref{char_fol_dim_2}). 
We prove that  the vector field $\tilde X_{h}$ is actually the one directing  such  characteristic foliation 
(Proposition \ref{Char_fol_Zw}). 
   Then, applying a local version of Gray's theorem (a contact version of the famous Darboux's theorem of symplectic geometry), we show that it is always possible to slightly deform $Z^{h}_{e_0}$ into a new submanifold $Z^{h_\e}_{e_0}$, with $h_\e$ an arbitrarily small perturbation of $h$, on which the vector field $\tilde X_{h_\e}$ is Morse-Smale (Proposition \ref{close.MS}). This implies the genericity result.\\

\noindent{\bf Related literature}.
The problem of exhibiting growth of Sobolev norms for linear, time dependent Schr\"odinger equations has been first addressed by Bourgain in \cite{BourgainGSN1999} for systems on  the torus (see also the recent result \cite{Chabert}). In the last years, a lot of efforts have also been done in the construction of time dependent perturbations $V$ of harmonic oscillators giving rise to solutions with unbounded trajectories. In particular, in the one-dimensional case, \cite{del} constructs the first example of bounded, pseudodifferential $V$ admitting solutions whose $\cH^s$-norm grows as $t^s$. \cite{Mas22} exhibits a wider class of such perturbations, which in \cite{Mas23} are proven to be generic among pseudodifferential operators of order zero. \cite{Maspero2018LowerBO} contains an example of $V$, still bounded and pseudodifferential, provoking $t^s$ growth of Sobolev norms for \emph{all} the solutions. In \cite{BGMR} an example of unbounded $V$ generating solutions with $\cH^s-$norm growing as $t^{2s}$ was constructed, while examples of multiplication operators provoking a logarithmic growth are exhibited in \cite{FaouRaphael}.\\
In the particular case of pseudodifferential perturbations which are quadratic polynomials both in $x$ and $D$, the growth rate at infinity of $\cH^s$-norms of the solutions in the one-dimensional setting has been completely described  in the works \cite{LLZ1, LLZ3, LLZ}, and recently extended in \cite{LLZ2} to higher dimensions.\\
We also mention \cite{BGMRV} for perturbations of the Landau Hamiltonian, \cite{HausMaspero} for  anharmonic oscillators on $\R^d$, and  \cite{Thomann2020GrowthOS} for  perturbations causing infinite growth of solutions in the Bargmann-Fock space.

Similarly, the closely related question of exhibiting a priori upper bounds on the growth of Sobolev norms for linear time-dependent systems  has been widely investigated. 
\noindent Concerning  perturbations of harmonic oscillators, in general one cannot expect better than a $t^s$ upper bound for the growth of the $\cH^s$-Sobolev norm of all solutions, which has been proven by \cite{del} in spatial dimension 1, and in the abstract work \cite{Maspero_2017} in arbitrary spatial dimension.
Note  that the solution of our Theorem \ref{main.informal} saturates such upper bound. 
Instead,  restricting to time quasi-periodic potentials and under non-resonant conditions on the time frequencies, a stronger $t^\epsilon$ upper bound, with $\e>0$ arbitrarily small,  has been proven in \cite{BGMRgrowth}, and under the further hypothesis that the perturbation has small size, reducibility results have been provided in \cite{BGMR, GrePat}, implying uniform $\cH^s-$boundedness in time of the solutions.

Coming to different models,  $t^\epsilon$ upper bounds have also been obtained for 
  systems with increasing spectral gaps (such as 
   the one dimensional anharmonic oscillator or the Laplace-Beltrami operator on Zoll manifolds) in \cite{Nenciu, DuclosLevStovicek, Maspero_2017, BGMRgrowth}, 
   and for perturbations of the Laplacian on $d$-dimensional flat tori  in
\cite{BourgainGSN1999, Delort2010, BertiMaspero, BLMgrowth}. In this last case,  it is important to exploit the specific cluster structure
of the spectrum of the Laplacian on tori.  
A generalization of such structural property has recently been identified in \cite{BLR2022, QN, QNbari} and used to extend the $t^\epsilon$ upper bound to a wider class of Schr\"odinger equations, such as perturbations of the Laplace-Beltrami operator on Lie groups and of the anharmonic oscillators on $\R^2$.

Before closing this introduction, we mention that  constructing  solutions with unbounded orbits in {\em nonlinear} Schr\"odinger-like equations is a big challenge.
We mention  the seminal works by Kuksin \cite{Kuk, Kuk2},
the breakthrough result by Colliander-Keel-Staffilani-Takaoka-Tao \cite{CKSTT},  its refinements and extensions obtained in \cite{guardia_kaloshin,
hani14,  haus_procesi15, guardia_haus_procesi16, GHMMZ, GG}, that exhibit strong forms of norm inflation.
Truly  unbounded  orbits  in nonlinear systems are known only in a handful of cases, {see \cite{gerard_grellier, hani15, GL, BanVega2022, GGG}}. We stress in particular the result of  \cite{Chabert2} dealing with time dependent perturbations of nonlinear harmonic oscillators on $\R^2$.\\
Finally we mention the recent paper \cite{MasperoMurgante} in which the method of the escape function is successfully extended to a one dimensional  fractional quasilinear Schr\"odinger equation, allowing to  construct solutions undergoing growth of Sobolev norms.

\smallskip

\footnotesize{
\noindent {\bf Acknowledgments.}
A. Maspero is supported by the European Union  ERC CONSOLIDATOR GRANT 2023 GUnDHam, Project Number: 101124921
and  by  PRIN 2022 (2022HSSYPN)   "TESEO - Turbulent Effects vs Stability in Equations from Oceanography".
B. Langella  is supported by  PRIN 2022 (2022HSSYPN)   "TESEO - Turbulent Effects vs Stability in Equations from Oceanography".
 Views and opinions expressed are however those of the authors only and do not necessarily reflect those of the European Union or the European Research Council. Neither the European Union nor the granting authority can be held responsible for them. We acknowledge the support of the INdAM groups GNAMPA and GNFM.
}
\normalsize

\subsection{Main result}
In this section we give a precise statement of our main result. 
Let us now start defining the set of functions and then the class of symbols that we will use: 

\begin{definition}[Homogeneous function]\label{def:pos_hom}
    A  smooth function $f\colon \R^4\setminus\{0\} \to \C$ is said to be {\em positively homogeneous of degree $m$} if
\begin{equation}\label{pos.hom}
f(\lambda x, \lambda \xi)= \lambda^m \, f(x, \xi), \qquad \forall \ \lambda >0 , \quad \forall \ (x, \xi)\neq 0
\end{equation} 
and {\em definitely homogeneous} if it is smooth at zero and \eqref{pos.hom} holds only for $\lambda \geq 1$ and $|(x,\xi)|\geq 1$.
\end{definition}

\begin{definition}[Symbols]
\label{def:symbol}
We define the class of symbols of order $\rho \in \R$ as
\[
\bS^{\rho}:=\left\{f \in C^{\infty}(\R^4, \C) \colon \forall \ \alpha, \beta \in \N^2_0 \ , \ \ 
\exists \ C_{\alpha, \beta}>0 \colon |\pa_x^{\alpha}\pa_{\xi}^{\beta}f(x, \xi)| \leq C_{\alpha, \beta}(1+|x|^2+|\xi|^2)^{\rho-\frac{|\alpha|+|\beta|}{2}}\right\}.
\]
We say that a symbol $f \in \bS^{\rho}$ is {\em classical} of order $\rho$ if there exist $f_{\rho}$ definitely homogeneous of degree $2\rho$ and   $\mu<\rho $
such that 
\begin{equation}\label{class_sym}
    f-f_{\rho} : =f_\mu \in \bS^{\mu}.
\end{equation}
We call $f_{\rho}$ the {\em principal symbol} of $f$, and denote by $\bS_{cl}^{\rho}$ the set of classical symbols of order $\rho$. 
\end{definition}

\begin{remark}$(i)$ The inclusion $\bold{S}_{cl}^0 \subset \bold{S}^0$ holds.\\
$(ii)$ The symbol
\begin{equation}\label{h0}
 h_0(x,\xi):=  \frac{|x|^2 + |\xi|^2}{2}    \ , \quad (x, \xi) \in \R^2 \times \R^2 \,,
\end{equation}
is  definitely homogeneous  of degree two, therefore it belongs to  $ \bS^1_{cl}$. \\
$(iii)$ If $v(x, \xi)$ is a positively homogeneous function of degree $0$, then 
\begin{equation}\label{gradh.ho}
\di v [\grad h_0] = \langle \grad v(x, \xi), \grad h_0 (x,\xi) \rangle = 0   \ .
\end{equation}
\end{remark}
The space $\bS^{\rho}$, endowed with the family of seminorms 
\begin{equation}\label{time.indep.seminorms}
    \wp_j^{\rho}(f):=\sum_{|\alpha|+|\beta| \leq j}\sup_{x, \xi \in \R^2}\frac{|\pa_x^{\alpha}\pa_{\xi}^{\beta}f(x, \xi)|}{(1+|x|^2+|\xi|^2)^{\rho- \frac{|\alpha|+|\beta|}{2}}}, \qquad \forall j \in \N,
    \end{equation}
is a Fréchet space with distance defined as 
\begin{equation}\label{time.indep.distance}
    \td^{\rho}(f, g):= \sum_{j\geq 0}\frac{1}{2^j}\frac{\wp_j^{\rho}(f-g)}{1+\wp_j^{\rho}(f-g)}.
\end{equation}
Analogously, we define also  the Fréchet space $C^0(\T, \bS^{\rho})$ of continuous maps $f \colon t \to f(t, \cdot)\in \bold{S}^{\rho}$, with finite seminorms: 
\begin{equation}
\wp_{\T, j}^{\rho}(f):=\sum_{\substack{|\alpha|+|\beta| \leq j} }\sup_{\substack{x, \xi \in \R^2\\t \in \T}}\frac{\Big|\pa_x^{\alpha}\pa_{\xi}^{\beta}f(t, x, \xi)\Big|}{(1+|x|^2+|\xi|^2)^{\rho-\frac{|\alpha|+|\beta|}{2}}}, \qquad \forall j \in \N.
\label{seminorm}
\end{equation}
The distance between two functions $f, g \in C^0(\T, \bold{S}^{\rho})$ is hence given by:
\begin{equation}\label{distance}
    \td^{\rho}_{\T}(f, g):= \sum_{j\geq 0}\frac{1}{2^j}\frac{\wp_{\T,j}^{\rho}(f-g)}{1+\wp_{\T,j}^{\rho}(f-g)}.
\end{equation}
More generally, given a  Fréchet space $\cF$ equipped with a sequence of seminorms $\{\wp_j\}_{j \in \N}$, we define on $\cF$ the natural distance
\begin{equation}\label{d.frechet}
\td_{\cF}(f,g) := \sum_{j\geq 0}\frac{1}{2^j}\frac{\wp_{j}(f-g)}{1+\wp_{j}(f-g)}\,.
\end{equation}

We now define pseudodifferential operators of order $\rho$:

\begin{definition}[Pseudodifferential operators]\label{def:pseudodiff}
We say that the operator $F$ belongs to $\cS^{\rho}$, the class of pseudodifferential operators of order $\rho$, if there exists a symbol $f \in \bold{S}^{\rho}$ such that $F=\Opw(f)$, where
\begin{equation}\label{Weyl_quantization}
\Opw(f)[\psi](x):=\frac{1}{2\pi}\iint_{ y, \xi \in \R^2}e^{\im (x-y)\xi} \, f\Big( \frac{x+y}{2}, \xi\Big)\, \psi(y) \, \di y \di \xi
\end{equation}
is the Weyl quantization of the symbol $f$. 
\end{definition}

Finally we define the resonant average of a symbol $v$ and the average of the operator $V$: we denote by $\phi^t_{h_0}(x_1, x_2, \xi_1,\xi_2)$ the classical flow of the harmonic oscillator in $\R^4$ with initial datum $(x_1, x_2,  \xi_1, \xi_2)$, where $h_0$ is defined as in \eqref{h0}:
\begin{equation}\label{h0.flow}
    \begin{split}
    \phi^t_{h_0}(x, \xi)& = (x_1\cos t+\xi_1 \sin t \,, \, 
    x_2\cos t+\xi_2 \sin t \, , \,  -x_1 \sin t+\xi_1\cos t \,  ,  \, -x_2 \sin t+\xi_2\cos t) \ .
    \end{split}
\end{equation}
One  easily checks the following result:
\begin{lemma}\label{composition_flow}
    Let $a(t, x, \xi) \in C^0(\T, \bold{S}^{0})$, then $a(t, \phi_{h_0}^t(x, \xi)) \in C^0(\T, \bold{S}^{0})$.
\end{lemma}

\begin{definition}[{Resonant average}]\label{resonant_average}
Given a symbol  $v \in C^0(\T, \bold{S}^0)$, we define the resonant average $\la v \ra$ with respect to the flow \eqref{h0.flow} of the harmonic oscillator as:
\begin{equation}\label{average}
\la v \ra (x, \xi):=\frac{1}{2\pi}\int_0^{2\pi} v(t, \phi^t_{h_0}(x, \xi)) \, \di t, 
\end{equation}
which is a time independent symbol of order zero. 
In addition if $v \in C^0(\T, \bS^0_{cl})$, then $\la v \ra \in \bS^0_{cl}$.

\noindent We define the average of the pseudodifferential operator $V \in C^0(\T, \cS^0)$, $V(t, x, D)=\Opw(v(t, x, \xi))$, as: 
\begin{equation}
 \la V \ra:=\Opw(\la v \ra) \ ,   
 \label{average_operator}
\end{equation}
which is a time independent operator, $\la V \ra \in \cS^0$. 
\end{definition}

\noindent{\bf Notation:} We denote 
by $\cX_h$ the Hamiltonian vector field of the function $h(x, \xi)$, i.e. 
$$
\cX_h:= \begin{pmatrix} 0 & 1 \\ -1 & 0 \end{pmatrix} \begin{pmatrix} \grad_x h \\ \grad_\xi h \end{pmatrix} = J \grad_{x,\xi} h \ 
$$
and the  Poisson bracket
\begin{equation}\label{poisson.bracket}
    \{ h , f \} := \grad_\xi h \cdot \grad_x f - \grad_x h \cdot \grad_\xi f  \  = \di f [ \cX_h] \ .
\end{equation}    

We are finally ready to state our main theorem: 
\begin{theorem}\label{main_theorem}
Denote by 
    \begin{equation}\label{setV}
    \begin{aligned}
        \cV := \Big\lbrace & v \in C^{0}(\T, \bS^0_{cl}) \colon \exists \  I\subset \Ran(\la v_0 \ra),  \exists \ a \in \bS^1_{cl} \mbox{ and } \exists \ \delta >0  \\
        & \mbox{ such that }  \{ \la v_0 \ra, a \} \geq \delta \quad   \forall \, (x, \xi) \in \la v_0 \ra ^{-1}(I) \cap \{|(x,\xi)| \geq 1 \}  \Big \rbrace  \,,
        \end{aligned}
    \end{equation}
    {where $v_0$ denotes the principal symbol of $v$}.
    Then 
    \begin{enumerate}
       \item[(i)] {\bf Genericity:} The set $\cV$   is generic  in $C^0(\T, \bS^0_{cl})$; precisely it is open and dense with respect to the metric $\td^{0}_{\T}$ in \eqref{distance}.
       
        \item[(ii)] {\bf Instability:} for every $v \in \cV$, for every $s>0$, there exist a function $u_0 \in \cH^s(\R^2)$ and a time $T>0$ such that the solution $u(t, x)$ of equation \eqref{Quantum_Harmonic_oscillator} with $V(t, x, D)=\Opw(v)$ and initial datum $u_0$
        satisfies 
\begin{equation}\label{lower.bound}
   \| u(t)\|_s \geq C_{s, u_0}\la t \ra^s, \qquad \forall \ t>T, 
\end{equation} 
for some positive constant $C_{s, u_0}$.
    \end{enumerate}
\end{theorem} 
Note that the lower bound in \eqref{lower.bound} is sharp, indeed in \cite{Maspero_2017} the authors prove that, for all $s>0$ and for all $u_0 \in \cH^s$,  there exists a constant $\tilde{C}_s >0$ such that 
$
\| u(t)\|_s \leq \tilde{C}_s \la t \ra^s \norm{u_0}_s$, showing that our unstable solutions grow with the optimal speed.
\vspace{1em}

\noindent{\bf Plan of the paper.} We conclude this section by giving an outline of the proof. The
   central element, as already discussed,  is the construction of the escape function $a$ for generic symbols $v \in C^0(\T, \bS^0_{cl})$, for which we proceed as follows.

In Section 3 we consider an arbitrary function $h$, positively homogeneous of degree zero, a regular value $e_0$ for $h$ and the projection $\tilde X_h$ of the  Hamiltonian vector field $X_h$ on the manifold $Z^h_{e_0}:= h^{-1}(e_0) \cap \S^3$. Following Colin de Verdière \cite{Colin_de_Verdi_re_2020}, we prove that if the flow of $\tilde \cX_h$ has a weakly hyperbolic simple structure (see Definitions \ref{def:wh}, \ref{def_simple_structure}) then there exists an escape function $a$ for $h$ at $e_0$, i.e. 
    $
    \{ h, a \} \geq \delta>0$ for any $(x, \xi) $ on $ h^{-1}(e_0)$ (Proposition \ref{simple_structure_thm}).
    
    In Section 4 we prove genericity of the set $\cV$. The main difficulty is to prove its density. This is achieved in several steps. 
    First we prove that  $\tilde \cX_h$ directs the characteristic foliation of $Z_{\mathbcal{e}_0}^h$ (Proposition \ref{Char_fol_Zw}). 
    Then, via a local Gray's theorem, we slightly perturb $h$ to a new positively homogeneous function  $h_{\epsilon}$ whose projected vector field $\tilde{\cX}_{h_{\epsilon}}$ is Morse-Smale (Proposition \ref{close.MS}), and actually has even  weakly hyperbolic simple structure (Proposition \ref{MS.WH}).
    This implies the genericity result in $  C^0(\T, \bS^0_{cl})$ (Proposition \ref{prop.density}). 
    
 In Section 5 we prove instability using the abstract result of \cite{Mas22}, reported in Theorem \ref{thm:ab0}. 
We verify its assumptions   by  proving  that  any operator $V$ with symbol $v \in \cV$ fulfills   the  Mourre estimate \eqref{mourre.est}.

\section{Preliminaries and abstract results}
We start by giving some preliminary definitions and results that we will need in the rest of the work.

\subsection{Symbolic calculus and the abstract result of \cite{Mas22} }\label{sec:pseudo}
In this section we state some results of symbolic calculus and next we report the abstract result of \cite{Mas22} (which in particular also applies to our pseudodifferential setting), that we use to prove point $(ii)$ of our main Theorem \ref{main_theorem} (see Section 5).

\paragraph{Symbolic calculus.}
We refer for example to \cite{shubin2001pseudodifferential} for the results of symbolic calculus:
\begin{proposition}\label{symbolic.calculus}
    Let $a \in \bold{S}^m$, $b \in \bold{S}^n$, $m, n \in \R$. Then:
    \begin{enumerate}
        \item For any $s \in \R$, there exist $C$, $M>0$ such that: $
          \|\Opw(a)u\|_{s-m} \leq C\wp_{M}^{m}(a)\|u\|_s;$  
        \item\label{aggiunto_Weyl}
        $\Opw(a)^*=\Opw(\overline{a}),$ (here $^*$ denotes the adjoint operator);
        \item There exists $c \in \bold{S}^{m+n}$ such that 
        $\Opw(a)\circ \Opw(b)=\Opw(c).$ Furthermore $c-ab \in \bold{S}^{m+n-1}$, and \begin{equation}\label{comp.symbols}
            \forall k \in \N \text{ , } \exists N, C>0 : \wp_k^{m+n-1}(c-ab) \leq C \wp_N^m(a) \wp_N^n(b);
        \end{equation}
        \item 
        There exists $d \in \bold{S}^{m+n-1}$ such that 
        $ \im [\Opw(a),\Opw(b)] = \Opw(d)$. 
        Moreover $d -\{a, b\}\in \bold{S}^{m+n-3}$, and \begin{equation}\label{comm.symbols}
        \forall k \in \N \text{ , } \exists N, C>0 : \wp_k^{m+n-3}(d-\{a, b\}) \leq C \wp_N^m(a) \wp_N^n(b);
        \end{equation}
        \item
        If $m<0$, then $\Opw(a)$ is compact as an operator from $L^2(\R^2)$ to itself. 
    \end{enumerate}
\end{proposition}

\begin{theorem}[Exact Egorov theorem]\label{egorov}
    Let $a \in \bold{S}^{\rho}$,  $\rho \in \R$.
    Then 
    $e^{\im tH_0}\Opw(a)e^{-\im tH_0} = \Opw(a \circ \phi^t_{h_0})$ for any  $ t \in \R$, 
where  $\phi^t_{h_0}$ is the Hamiltonian flow generated by $h_0$, see \eqref{h0.flow}.
\end{theorem}

Next we state a strong Gårding inequality \cite{shubin2001pseudodifferential}, referring to \cite{Mas23} for a proof of this statement. 

\begin{theorem}[Strong Gårding inequality]\label{strong.garding}
    Let $a \in \bS^0$ be a symbol of order zero, and suppose that there exists $R>0$ such that
    \begin{equation}\label{pos.sym}
        a(x, \xi) \geq 0, \qquad \forall \ |(x, \xi)| \geq R.
    \end{equation}
    Then there exists a positive constant $C$ such that 
    $$
        \la \Opw(a) u, u \ra \geq -C \norm{u}_{-\frac12}^2, \qquad \forall \ u \in L^2(\R^2).
    $$
\end{theorem}
We conclude stating the following result about composition of symbols with compactly supported functions. It is an application of pseudodifferential functional calculus, for which we refer to \cite{Davies.functional.calculus} and \cite[Appendix B]{Mas22}.
\begin{lemma}\label{lemma_compactness_cutoff}
Consider a smooth compactly supported real function $g \in C_c^{\infty}(\R, \R)$.
\begin{itemize}
    \item[(i)] Let $a \in \bold{S}^0$ be real valued. Then the operator $g(\Opw(a))$, defined via functional calculus, is a pseudodifferential  operator of order zero: $g(\Opw(a)) \in \cS^0$. 
Moreover 
$g(\Opw(a))-\Opw(g(a)) \in \cS^{-1}$.
    \item[(ii)] Let $A$, $B$ be bounded selfadjoint operators in  an Hilbert space. If $A-B$ is compact, so is  $g(A)-g(B)$. 
\end{itemize}
\end{lemma}

\paragraph{The abstract result of \cite{Mas22}.}
\label{sec:alg} We now report the abstract result proven in \cite{Mas22}, {in a small variant suited for our setting}.
Given an Hilbert space $\cH$, endowed  with the scalar product $\la \cdot, \cdot \ra$,   let  $K_0$ be a   selfadjoint,   positive operator,
with compact resolvent. 
Define a scale of Hilbert spaces by $\cH^r:=D(K_0^r)$
(the domain of the operator $K_0^r$) if $r\geq 0$, and
$\cH^{r}=(\cH^{-r})^\prime$ (the dual space) if $ r<0$, and denote by $\cH^{-\infty} = \bigcup_{r\in\R}\cH^r$ and $\cH^{+\infty} =
\bigcap_{r\in\R}\cH^r$. 
We endow $\cH^r$ with the norm
$\norm{\psi}_r:= \norm{K_0^r \psi}_{0}$, where $\norm{\cdot}_0$ is
the norm of $\cH^0 \equiv \cH$.       

On $\cH^r$ consider the  abstract Schr\"odinger equation 
\begin{equation}\label{eq.ab0}
    \im \pa_t  \psi = K_0 \psi + V(t) \psi ,
    \end{equation}
with  $V(t)$ a time dependent selfadjoint operator.
      
Next assume to have a graded algebra $\cA$ of operators which mimic  pseudodifferential
operators.  For $m\in\R$ let $\cA_m$ be a linear subspace of $
\bigcap_{s\in\R}\cL(\cH^s,\cH^{s-m})$  (where $\cL(\cH^s,\cH^{s-m})$ denotes the space of linear and bounded operators from $\cH^s$ to $\cH^{s-m}$) and define
$\cA:=\bigcup_{m\in\R}\cA_m$.  We notice that the space
$\bigcap_{s\in\R}\cL(\cH^s,\cH^{s-m})$, is a Fr\'echet space equipped
with the seminorms: $\Vert A\Vert_{m,s} := \Vert
A\Vert_{\cL(\cH^s,\cH^{s-m})}$.
 We say that $A$ is of {\em order $m$} if
$A \in \cA_m$.

We shall also use the following notations.    
For $\Omega\subseteq \R^d$ and $\cF$ a Fr\'echet space, 
we denote by $C_b^0(\Omega, \cF)$ the space of
$C^m$ maps $f: \Omega\ni x\mapsto f(x)\in\cF$ such that,  for every
seminorm $\norm{\cdot}_j$ of $\cF$, one has
       \begin{equation}
       \label{star}
       \sup_{x\in\Omega}\Vert\partial_x^\alpha f(x)\Vert_{j} < +\infty
       \ , \quad \forall \alpha\in \N^d\  :\ \left|\alpha\right|\leq m  \ .
       \end{equation}
If \eqref{star} is true
$\forall m$, we say $f \in C^\infty_b(\Omega, \cF)$. 
Similarly we denote by $C^\infty(\T, \cF)$ the space of smooth maps from the torus $\T = \R/(2\pi\Z)$ to the Fr\'echet space $\cF$.  
Given  two operators $A, B \in \cL(\cH)$, we write $A \geq B$  with the meaning   $\la A \vf, \vf \ra \geq \la B \vf, \vf \ra$ $\, \forall \vf \in \cH$.\\

We now state a list of assumptions on $\cA, K_0, V(t)$ required to apply the abstract theorem in \cite{Mas22}.

\paragraph{Assumption I: Pseudodifferential algebra}\label{ass.1}
     \begin{itemize}
      \item[(i)] For each $m\in \R$, $K_0^m\in\cA_m$; in particular
        $K_0$ is an operator of order one.
\item[(ii)]  For each $m\in\R$, $\cA_m$ is a Fr\'echet space for a family of  filtering semi-norms $\{ \wp^m_j \}_{j\geq 0}$  such that the embedding $\cA_m\hookrightarrow  \bigcap_{s \in \R} \cL(\cH^s,\cH^{s-m})$ is continuous\footnote{A family of seminorms $\{ \wp^m_j \}_{j\geq 0}$  is called filtering if for any $j_1, j_2 \geq 0$  there exist  $k \geq 0$ and $c_1, c_2 >0$  such that the two inequalities $\wp^m_{j_1}(A) \leq c_1 \wp^m_{k}(A)$ and $\wp^m_{j_2}(A) \leq c_2 \wp^m_{k}(A)$ hold for any $A \in \cA_m$.}. 
 If $m^\prime \leq m$ then $\cA_{m^\prime}\subseteq\cA_m$ with a continuous embedding.

    \item[(iii)]    
    For all $m,n\in \R$, if $A\in \cA_m$   and $B\in\cA_n$ then $A B\in\cA_{m+n}$ 
     and the map $(A,B)\mapsto AB$ is continuous from 
     $\cA_{m}\times\cA_{n}$ into $\cA_{m+n}$. 
     \item[(iv)] 
     If $A\in \cA_m$   and $B\in\cA_n$, then the commutator $[A,B]\in\cA_{m+n-1}$  and the map $(A,B)\mapsto [A,B]$
      is continuous from $\cA_{m}\times\cA_{n}$ into $\cA_{m+n-1}$.
      
     \item[(v)] $\cA$ is closed under  perturbation by smoothing operators in the following sense:
     let $A$ be a linear map: $\cH^{+\infty}\rightarrow\cH^{-\infty}$. If   there  exists $m\in\R$ such that  for  every $N>0$  we  have  a   decomposition        $A= A^{(N)}+ S^{(N)}$,  with $A^{(N)}  \in\cA_m$  and $S^{(N)}$  is  $N$-smoothing\footnote{Namely, if  $\forall s \in \R$, it  can be extended  to an operator in   $\cL(\cH^{s}, \cH^{s+N})$.  When this is true  for  every $N\geq 0$, we call it a    smoothing  operator},  then
      $A\in\cA_m$. 
       \item[(vi)] If $A \in \cA_m$ then also the adjoint operator $A^* \in \cA_m$.
  The  duality here is defined by the scalar product 
      $\langle \cdot, \cdot\rangle$ of $\cH=\cH^0$. The adjoint $A^*$ is defined by 
      $\langle  u, Av\rangle  = \langle A^*u, v\rangle$  for $u,v\in\cH^\infty$  and extended by continuity.  
     \end{itemize}

     \noindent
       {\bf Assumption II: Properties of $K_0$}
       \begin{itemize}\label{ass.2}
       \item[(i)] The operator $K_0$ has purely discrete spectrum fulfilling,        for some  $\lambda \geq 0$, 
       $  \sigma(K_0) \subseteq \N +\lambda$.
 \item[(ii)] For any $m \in \R$ and $A\in\cA_m$,   the map defined on $\R$ by 
         $\tau\mapsto A(\tau):={\rm e}^{\im \tau K_0}\, A \, {\rm e}^{-\im \tau
         K_0}$ belongs to $ C^\infty_b(\R, \cA_m)$ and one has:
   $  \forall  j \quad \exists N \ s.t.\  \ \sup_{\tau \in \R}  \wp^{m}_j(A(\tau)) \leq C \, \wp^m_N(A) $ 
  for some positive constant $C(m,j)$.
\end{itemize}

    \noindent
       {\bf Assumption III: Properties of the potential $V(t)$}\label{ass.3} \\ 
      The operator   $V \in C^{{0}}(\T, \cA_{0})$, $V(t)$  selfadjoint $\forall t$, and  its resonant average 
      \begin{equation}\label{Vaverage}
     \lala V \rara := \frac{1}{2\pi} \int_0^{2\pi}  e^{\im s K_0} \, V(s) \, e^{- \im s K_0} \di s  
\end{equation}
      fulfills:
      \begin{itemize}
      \item[(i)]  There exists  an interval  $I_0\subset \R$  such that  $\abs{\sigma(\lala V \rara)\cap I_0} >0 $; here $\abs{\cdot}$ denotes the Lebesgue measure.
      \item[(ii)]  {\em Mourre estimate} over  $I_0$:  there exist a  selfadjoint  operator $A \in \cA_{1}$ and a function   $g_{I_0} \in C_c^\infty(\R, \R_{\geq 0})$ with $g_{I_0}\equiv 1$ on $I_0$ such that 
\begin{equation}\label{mourre.est}
g_{I_0}(\lala V \rara) \, \im [\lala V \rara, A] \, g_{I_0}(\lala V \rara ) \geq \theta \, g_{I_0}(\lala V \rara)^2  + K
\end{equation}
for some $\theta >0$ and $K$ a selfadjoint compact operator.
      \end{itemize}
{The following is a small variant of the} abstract result in \cite{Mas22}:
\begin{theorem}
\label{thm:ab0}
Assume that $\cA$ is a graded algebra as in {\rm Assumption {\rm I}}, and that $K_0 $ and $V(t) \in C^{{0}}(\T, \cA_0)$ satisfy {\rm Assumptions} {\rm II} and {\rm III}.  Then 
for any  $r >0$  there exist a solution $\psi(t)$ of \eqref{eq.ab0} in $\cH^r$ and constants $C , T  >0$ such that 
$
\norm{\psi(t)}_{r} \geq C \la t \ra^{r} , \quad \forall t \geq  T \ .
$ 
\end{theorem} 
\begin{proof}
{The proof goes through exactly as the one of Theorem 2.7 of \cite{Mas22}, which is proven in Section 3 therein. The only difference is that here we only suppose $V \in C^0(\T; \cA^0)$, instead of requiring the stronger hypothesis $V \in C^\infty(\T; \cA^0)$, assumed in \cite{Mas22}. In particular, all the results of Section 3 of \cite{Mas22} remain unchanged, with the only exception of Lemma 3.3, which is the only point where this hypothesis might enter. Given a time periodic family $W \in C^\infty(\T; \cA^m)$ for some $m$, such lemma ensures the existence of a solution $X$ of the homological equation
$$
\partial_t X + \im [K_0, X(t)]= W(t) - \int_0^{2\pi} e^{-\im s K_0} W(t + s) e^{\im s K_0} \ ds\,.
$$
However, for our purpose it is sufficient to note that $X$ gains regularity in time w.r.t. $W$, namely $X \in C^1(\T; \cA^m)$ if $W \in C^0(\T; \cA^m)$.}
\end{proof}

\subsection{Morse-Smale vector fields}
\label{dynamical_systems}

In this section we recall the definition of  Morse-Smale vector fields  and state a genericity result by Peixoto, that we shall use in Section 4. 

Given $M$ a smooth manifold and $X$ a smooth vector field on $M$, we denote by $\phi^t_X(x)$ the flow of $X$ that passes through $x$ at time $t=0$. 
We refer to  \cite{palis1998geometric} for most of the definitions related to dynamical systems.

\begin{definition}[Attractor and repellor]\label{def:attractor}
Let $M$ be a smooth closed manifold and $X$ a smooth vector field on $M$. 
A compact, invariant subset  $K^+ \subset M$  [respectively $K^-$] is called an {\em attractor} [respectively a {\em repellor}] if there exists an open neighbourhood $U^+$ of $K^+$ [respectively $U^-$ of $K^-$] in $M$ such that 
\begin{equation}
    K^+=\bigcap_{t \geq 0} \phi^t_X(U^+) \quad \big[\mbox{respectively} \ \ \ 
  K^-=\bigcap_{t \leq 0} \phi^t_X(U^-)
    \big] \ .
    \label{attractor}
\end{equation}
If $K^+$ is an attractor, the {\em basin} of $K^+$ is the set of points $x$ such that $\phi^t_X(x) \to K^+$ as $t \to +\infty$. 
One defines, mutatis mutandis, the basin of $K^-$.
\end{definition}

We recall the definition of  Morse-Smale vector field on a surface, referring to \cite[Definition 4.6.8]{geiges2008introduction}: 
\begin{definition}[Morse-Smale vector field]\label{Morse-Smale}
A vector field $X$ on a closed, orientable surface $M$ is said to be of {\em Morse–Smale}
type if it has the following properties:
\begin{enumerate}\label{characterization MS}
    \item There are finitely many singularities and closed orbits, all of which
are non-degenerate \footnote{{Recall} that a fixed point $p$ is non-degenerate if any eigenvalue of the linearized vector field at $p$ has non-vanishing real part. A closed orbit $\gamma$ is called non-degenerate if given $p \in \gamma$, $p$ is a non-degenerate fixed point of the Poincaré return map.};
    \item There are no trajectories connecting  saddle points;
    \item The $\alpha$ and $\omega$–limit set of each flow line is either a singular point or
a closed orbit.
\end{enumerate} 
\end{definition}
An important property of Morse-Smale vector fields is that, on any given closed and orientable surface, they are dense in the set of smooth vector fields, as it was proven by Peixoto  \cite[Theorem 2]{Peixoto1962StructuralSO}.

\begin{theorem}[Peixoto]\label{generic.MS}
    Let $S$ be a smooth closed orientable surface. The set of smooth Morse-Smale vector fields is open and dense in the set of smooth vector fields on $S$, with respect to the $C^{\infty}$–topology.
\end{theorem}

We also recall that, given compact smooth manifolds $M$ and $N$, the $C^k(M; N)$ norm of a function $f:M \to N$ is defined once finite atlases $(M_i, \phi_i)$, $(N_i, \psi_i)$ and a partition of the unity $\{\chi_i\}_i$ subordinated to $\{M_i\}_i$ are fixed on $M$ and $N$. Then $\norm{f}_{C^k(M; N)}:=\max_{i,j} \sup_{\substack{x \in M_i\\ 0 \leq |\alpha| \leq k}} |\pa^{\alpha}_x (\psi_j \circ \left(\left. {(  \chi_i f)} \circ \phi_i^{-1}\right)\right|_{M_j})|$ and $\|f\|_{C^\infty(M; N)}$ is defined according to \eqref{d.frechet}, with $\wp_k(f):= \|f\|_{C^k(M; N)}$.
When $M = N$, we will simply write $\| \cdot\|_{C^k(M)}$.\\
Furthermore, given a vector field $X$ on a compact smooth manifold $M$ of class $C^k$, we define its norm as 
\begin{equation}\label{distance.vf}
   \norm{X}_{C^k}:= \max_{i} \norm{X_i}_{C^k(M; \R)}
\end{equation}
where $X_i$ are the components of $X$. 
Then the $C^{\infty}$ norm is defined according to \eqref{d.frechet} with $\wp_k(X):=\norm{X}_{C^k}$.



\subsection{Contact geometry}
In this section we  recall some basic concepts of contact geometry, referring to \cite{geiges2008introduction, Giroux1991}. 

\begin{definition}[Contact structure]\label{contact.man}
Let $M$ be a odd dimensional smooth manifold. We say that a smooth field of hyperplanes  $\kappa \subset TM$ is a {\em contact structure} on $M$ if there exists a 1-form $\alpha$ such that $\kappa = \operatorname{ker}(\alpha)$ and $\di \alpha|_{\ker (\alpha)}$ is non-degenerate.
We say that $\alpha$ is a {\em contact form} on $M$ and we say that the pair ($M$, $\kappa$) is a {\em contact manifold} .  
\end{definition}

\begin{remark}\label{rem:contact}
The contact form is not uniquely determined by the contact structure: indeed if $\alpha$ is a contact form, then for every $f: M \to \R$ never vanishing, $f \alpha$ is also a contact form for $(M, \kappa)$. 
\end{remark}

\begin{remark}\label{char.contact}
    Let $\alpha$ be a 1-form on $M$, $\dim M = 2n+1$,
    then 
    $
    \di \alpha|_{\ker (\alpha)} \mbox{ is non-degenerate if and only if } \alpha \wedge (\di \alpha)^n \neq 0$, (where by $(\di \alpha)^n$ we denote the wedge product of $\di \alpha$ repeated $n$ times) see e.g. \cite[
    Section I.1.A]{Giroux1991}.
\end{remark}
The following lemma is well known (see \cite[Lemma 1.4.5]{geiges2008introduction}), but we sketch its proof for completeness.    
\begin{lemma}\label{Lemma_contact_transversality}
    Let $(V, \Omega)$ be an exact symplectic manifold of dimension $2n$ and let $X$ be a  {\em Liouville vector field} on $V$, namely a vector field fulfilling $\cL_{X} \Omega = \Omega$. 
    Every hypersurface $E\subset V$ of dimension $2n-1$ and transverse to $X$ (i.e. $X(p) \notin T_pE, \ \forall \ p \in E$) is a contact hypersurface with contact form $\alpha:=i_X\Omega \vert_E$ and contact structure $\kappa_p=\ker(\alpha_p)$.
\end{lemma}

\begin{proof}
We  show the equivalent condition  $\alpha \wedge (\di \alpha)^{(n-1)} \neq 0$.
Using Cartan's formula and that $X$ is a Liouville vector field we have
    $$
\di \alpha=\di(i_X\Omega)=\cL_X\Omega=\Omega \quad \implies 
    \alpha\wedge (\di \alpha)^{n-1}=i_X\Omega \wedge \Omega^{n-1}=\frac{i_X(\Omega^n)}{n}. 
    $$
We now claim that $i_X(\Omega^n) \neq 0$ as $2n-1$ form on $E$, thus proving the lemma. 
Indeed, since $E$ is transversal to $X$,  it is possible to pick, for every point $p \in E\subset V$, a basis for $T_pV$ of the form 
$
\cB=\{e_1,\ldots,e_{2n-1}, X(p)\},
$ 
where $\{e_j\}_{j=1}^{2n-1}$ is a basis of $T_pE$.
Evaluating the volume form $\Omega^n$ at this basis we have that 
$  0 \neq \Omega^n(X(p), e_1,\ldots, e_{2n-1})=i_X(\Omega^n)(e_1, \ldots ,e_{2n-1})$,
    proving that $i_X (\Omega^n) \neq 0$.
\end{proof}
As a consequence we have:
\begin{corollary}\label{contact.sphere}
 The sphere $\S^3\equiv\{(x, \xi) \in \R^4 : |(x, \xi)|=1\}$ is a contact manifold, with contact form  locally given by 
  $\alpha:=\left.\frac12 (-\xi_1 \di x_1 -\xi_2 \di x_2 +x_1 \di \xi_1 + x_2 \di \xi_2)\right|_{\S^3}.$
\end{corollary}
\begin{proof}
We regard $\S^3$ as a submanifold of the exact symplectic manifold $(\R^4, \Omega)$, with $\Omega$ the standard symplectic form. A  Liouville vector field is given by 
$\frac12 \nabla h_0,$
with $h_0$ as in \eqref{h0}, which is transverse to $\S^3$. 
Lemma \ref{Lemma_contact_transversality} yields the result, with 
$ \alpha:= i_{\frac{\nabla h_0}{2}}\Omega \vert_{\S^3}$ having the claimed form.
\end{proof}

The following result is a local reformulation of the result in \cite[Theorem 2.2.2]{geiges2008introduction}.
\begin{proposition}[{Local Gray Theorem}]\label{Gray} 
    Let $M$ be a smooth closed contact manifold and $S \subset M$ a codimension one smooth compact submanifold. Let $\delta>0$. Suppose that there exist $\cU(S) \subset M$ tubular neighbourhood of $S$, and $\{\alpha_t\}_{t \in [0, 1]}$ smooth family of contact forms defined on $\cU(S)$ such that\footnote{{Here to define $\norm{\alpha}_{C^{\infty}(\cU(S))}$ we pick a finite open covering $\{U_i\}$ of $\overline{\cU(S)}$ on which local coordinates $y=(y_1, y_2, y_3)$ are defined, and we write locally
$\alpha=\sum_{i=1}^3 \alpha_i(y) \di y_i$ in such coordinates. Then, for every $k \in \N$, the $C^k$ norm of $\alpha$ is 
$
\norm{\alpha}_{C^k}:=\max_{i} \sum_{j=1}^3 \sup_{\substack{  y \in U_i\\
0 \leq |\beta| \leq k}} |\pa_y^\beta \alpha_j(y)|,
$
}
and the $C^\infty$ norm is defined via the Fréchet distance \eqref{d.frechet}.} 
\begin{equation}
\label{norm.alpha.dot}
    \norm{\dot{\alpha}_t}_{C^{\infty}(\cU(S))} < \delta'\,.
    \end{equation}
    Then there exists $\delta'_0>0$ such that, if $\delta'< \delta'_0$, there exist a neighborhood $\cU'(S) \subset \cU(S)$ of $S$, and two smooth families of maps $\{\psi_t\}_{t \in [0,1]}$, $\psi_t: \cU'(S) \rightarrow \cU(S)$, $\{\lambda_t\}_{t \in [0,1]},$ $\lambda_t : \cU(S) \rightarrow \R^+$, such that,  for any $t$, $\psi_t$ is a diffeomorphism on $\psi_t(\cU'(S)),$
    \begin{gather}
    \label{psi.id}
        \|\psi_t - \uno\|_{C^\infty(\cU'(S))} < \delta \,,\\
        \label{iso.psi}
        \psi^*_t \alpha_t = \lambda_t \alpha_0
    \end{gather}
    as 1-forms on $\cU'(S)$\,.
\end{proposition}

\begin{proof}
    The idea is to look for $\psi_t$ as the flow  at time $t$ of a  suitable vector field $Y_t \in \kappa_t:=\ker(\alpha_t)$, defined in $\cU(S)$. More precisely, suppose that 
    \begin{equation}\label{yt.small}
\|Y_t\|_{C^\infty(\cU(S))} < \delta''
    \end{equation} for some positive $\delta''$. Take local coordinates $(z, \zeta)$ around $S$ such that $S$ can be identified with the surface $\{ z = 0\}$ and, given $r>0$, let $\cU_r(S)$ be the neighborhood of $S$ corresponding to $|z| < r$. If $r$ and $\delta''$ are small enough, then $\cU_{r + \delta''} (S) \subset \cU(S)$.
    Then the Cauchy problem $\dot \psi_t = Y_t(\psi_t),$ $\psi_0 = \uno$, is well defined on $\cU'(S):= \cU_r(S)$, $\psi_t(z) \in \cU(S)$ for any $z \in \cU'(S)$ and $t \in [0,1].$ Furthermore, by \eqref{yt.small} one deduces that $\forall k \in \N$ there exists $D_k>0$ such that $\|\psi_t - \uno\|_{C^k(\cU'(S))} \leq \delta'' t D_k$ thus, recalling
    \eqref{d.frechet} with $\wp_k(\cdot ) := \|\cdot \|_{C^k(\cU'(S))}$, and eventually shrinking $\delta''$,  also $\|\psi_t - \uno\|_{C^{\infty}(\cU'(S))}\leq \delta t.$ 
    Thus, to prove the result it is sufficient to find $Y_t$ such that \eqref{iso.psi}, \eqref{yt.small} hold.
    
    \noindent Now, derive \eqref{iso.psi} with respect to time. Since  $\frac{\di}{\di t} \psi^*_t(\alpha_t) = \psi^*_t\left(\dot{\alpha_t} + \cL_{Y_t} (\alpha_t)\right)$ (see for example \cite[Lemma 2.2.1]{geiges2008introduction}), using Cartan's formula, and the fact that $Y_t \in \kappa_t$, we obtain
\begin{equation}\label{gray1}
        \frac{\di}{\di t} \psi^*_t(\alpha_t) = \psi^*_t( \dot{\alpha}_t +i_{Y_t} \di \alpha_t) = \dot{\lambda_t} \alpha_0\,.
    \end{equation}
    Since
    $$\dot{\lambda_t} \alpha_0= \frac{\dot{\lambda_t}}{\lambda_t} \psi_t^*(\alpha_t)=\psi_t^*(\mu_t \alpha_t), \quad \mbox{ with } \mu_t:= \pa_t(\log(\lambda_t))\circ \psi_t^{-1},$$
    equation \eqref{gray1} reads: 
    \begin{equation}\label{gray2}
       \psi^*_t( \dot{\alpha}_t +i_{Y_t} \di \alpha_t)=\psi_t^*(\mu_t \alpha_t) \quad \iff \quad \dot{\alpha}_t +i_{Y_t} \di \alpha_t=\mu_t \alpha_t\,. 
    \end{equation}
    Let now $R_{\alpha_t}$  be the Reeb vector field associated to  $\alpha_t$, i.e. the vector field on $\cU(S)$ uniquely defined by $\alpha_t(R_{\alpha_t})=1$, $\di\alpha_t(R_{\alpha_t}, \cdot)=0$. 
    Evaluating \eqref{gray2} in $R_{\alpha_t}$ gives 
    $
    \mu_t=\dot{\alpha_t}(R_{\alpha_t}) \,,
    $
    which defines $\mu_t$.
    Hence, by \eqref{gray2}, we obtain \begin{equation}\label{gray3}
     \di \alpha_t(Y_t, \cdot)\vert_{\kappa_t}= \left(\mu_t \alpha_t -\dot{\alpha_t}\right)\vert_{\kappa_t} = -\dot\alpha_t\vert_{\kappa_t}\,,   
    \end{equation}
    which uniquely determines $Y_t$ due to the non degeneracy of $\di \alpha_t\vert_{\kappa_t}.$ Finally, the smallness condition \eqref{yt.small} with $\delta''=C \delta'$ for some positive $C$, depending only on the coefficients of $\di \alpha_t$,
    follows from \eqref{norm.alpha.dot}.
\end{proof}

\begin{lemma}
Let $\{\psi_t\}_{t \in [0,1]}$ as in Proposition \ref{Gray}. For any $t \in [0,1]$ there exists $\delta_0>0$ such that, if $\delta < \delta_0$, then $S \subset \psi_t(\cU'(S))$.
\end{lemma}
\begin{proof}
   Take a tubular neighborhood $\cU''(S)$ of $S$ such that $\overline{\cU''(S)} \subset \cU'(S)$ and for any $\zeta \in S$ consider the map $F_{t,\zeta}: \overline{\cU''(S)} \rightarrow \cU(S)$ defined by $F_{t,\zeta}(\eta) := \zeta - (\psi_t(\eta) - \eta)$. Then by \eqref{psi.id}, if $\delta< \delta_0$ is small enough, $F_{t,\zeta}$ maps $\overline{\cU''(S)}$ into itself and it is a contraction. Moreover, compactness of $S$ ensures that $\delta_0$ can be taken independently of $\zeta$. Then the map $F_{t,\zeta}$ admits a  fixed point $\eta_0$, namely a point such that $\psi_t(\eta_0) = \zeta$.
\end{proof}
\paragraph{Characteristic foliation.}
First recall that
a 1-dimensional oriented  foliation $\mathscr{F}$ on a manifold $S$ is the equivalence class of vector fields on $S$ positively proportional to a given vector field  $X$ on $S$, which we say to {\em direct} the foliation. 
    Note that $\mathscr{F}$ is a singular foliation since  $X$ might vanish.

Consider now a hypersurface $S$ embedded in a contact manifold $(M, \kappa)$. We define its {\em characteristic foliation} \cite{Giroux1991, geiges2008introduction}:
\begin{definition}\label{char_fol_dim_2}
Let $S$ be a {oriented} hypersurface in a contact {3-dimensional} manifold $(M, \kappa)$ and fix $\mu$ a volume form on $S$ {(compatible with the given orientation)}.
The characteristic foliation $\mathscr{F}$ is the foliation directed by the vector field $X$ such that 
\begin{equation}\label{def.charX}
i_X \mu = \alpha\vert_S,
\end{equation}
where $\alpha$ is a contact form for $(M, \kappa)$. 
\end{definition}
Note  that different choices of volume form $\mu$ and of  contact
form $\alpha$  differ by a positive function (Remark \ref{rem:contact}), {so the vector field $X$  in \eqref{def.charX} is defined only up to  a positive function, however the induced foliation is the same.}

\section{Construction of the escape function}\label{sec:escape}

In this section we consider a smooth function $h\colon \R^4\backslash\{0\} \to \R$ positively homogeneous of degree zero (see Definition \ref{def:pos_hom}) and prove that,  under suitable dynamical conditions on its flow, 
it is possible to construct an escape function for it (see Definition \ref{def:escape}). 
Our construction extends to our setting the one done by Colin de Verdière in \cite{Colin_de_Verdi_re_2020} for pseudodifferential operators on the two-dimensional torus.

 \begin{remark}\label{homo.flow}
If $h$ is a positively homogeneous function of degree 0, its 
Hamiltonian vector field $\cX_h$  is homogeneous of degree $-1$; consequently the flow 
 $\Phi^t_{\cX_h}(x, \xi)$ 
    satisfies 
    \begin{equation}
        \Phi^t_{\cX_h}(\lambda x, \lambda \xi)=\lambda \Phi^{t / \lambda^2}_{\cX_h}(x, \xi), \qquad \forall  \lambda>0, \  t \in \R. 
    \end{equation}
\end{remark}

\subsection{Setting and main result}\label{subsec:escape}
We first give the definition of escape function.
\begin{definition}[Escape function]\label{def:escape}
    Let $h: \R^4 \backslash\{0\} \to \R$ be a  positively homogeneous function of degree zero (see Definition \ref{def:pos_hom}), and let  $\mathbcal{e}_0 \in \R$ be a regular value for $h$ (i.e. $\di h(\mathbcal{e}_0)\neq 0$).
    Denoting by 
    \begin{equation}\label{Sigma0}
    \Sigma_{\mathbcal{e}_0}^h :=h^{-1}(\mathbcal{e}_0)
    \end{equation}
    the energy shell associated to $\mathbcal{e}_0$, we say that a function $a: \R^4 \setminus \{0\} \to \R$,
    positively  homogeneous of degree two, is an {\em escape function} for $h$ at $\mathbcal{e}_0$ if the Poisson brackets $\{h, a\}  $ (see \eqref{poisson.bracket}) satisfy: 
    \begin{equation}\label{Poiss.esc}
    \{h, a\} \geq \delta >0  \quad \mbox{ on } \Sigma_{\mathbcal{e}_0}^h  \ . 
    \end{equation}
\end{definition}

\noindent Let $\boldsymbol{\pi}$ be  the projection of $\R^4\setminus\{0\}$ on $\S^3$, i.e., 
\begin{equation} 
\boldsymbol{\pi}\colon \R^4\setminus\{0\} \to \S^3, \qquad 
\boldsymbol{\pi}(x, \xi)=\Big(\frac{x}{\rho}, \frac{\xi}{\rho}\Big)=:\zeta, \qquad \rho : = \sqrt{|x|^2 + |\xi|^2} \,
\label{projection}
\end{equation}
and denote by $Z_{\mathbcal{e}_0}^h$ {the set}
\begin{equation}\label{Zwo}
Z_{\mathbcal{e}_0}^h:=\boldsymbol{\pi}(\Sigma_{\mathbcal{e}_0}^h)= \Sigma_{\mathbcal{e}_0}^h \cap \S^3\,.
\end{equation}
Note that $Z_{\mathbcal{e}_0}^h$ is a smooth manifold of dimension 2 since, due to the homogeneity property of $h$, $\mathbcal{e}_0$ is still a regular value of $h \vert_{\S^3}$ (which follows from \eqref{gradh.ho} with $v \leadsto h$), {and that, due to the homogeneity of $h$, 
the energy shell $\Sigma_{\mathbcal{e}_0}^h$ is a cone.}

\begin{remark}\label{rmk.finite}
    Note that $Z^h_{\mathbcal{e}_0}$ is compact, thus it only has a finite number of connected components\footnote{{To see this, recall that the union of all the connected components of $Z^{h}_{\mathbcal{e}_0}$ is an open covering of $Z^h_{\mathbcal{e}_0}$, see for instance Proposition 1.8(d) of \cite{lee2003introduction}; then the finiteness of the covering follows by compactness.}}.
   Moreover, if $Z^h_{\mathbcal{e}_0}$ has $n$ connected components, then $\Sigma^h_{\mathbcal{e}_0}$ is the disjoint union of $n$ open connected cones. Therefore in the present section we are going to assume without loss of generality that $Z^h_{\mathbcal{e}_0}$ and $\Sigma^h_{\mathbcal{e}_0}$ are connected, since if they are not, the construction of the escape function presented in this section can be repeated separately on each one of the finitely many open cones forming $\Sigma^h_{\mathbcal{e}_0}$.
\end{remark}

Due to its cone structure, we can parametrize 
$\Sigma_{\mathbcal{e}_0}^h $
with ``polar'' coordinates of the form $(\rho, \zeta) \in \R^+ \times Z_{\mathbcal{e}_0}^h$. We now express the  Hamiltonian vector field $\cX_h$ in such coordinates.
 
\begin{lemma}\label{lemma_Xh}
Let $h$ be a positively homogeneous function of degree $0$, and let $\cX_h$ be its Hamiltonian vector field.
    There exists a  vector field $\tilde \cX_h$ on $\Sigma_{\mathbcal{e}_0}^h$, homogeneous of order $0$, such that $\tilde \cX_h|_{Z_{\mathbcal{e}_0}^h} \in TZ_{\mathbcal{e}_0}^h$ and
    \begin{equation}\label{Xh}
        \cX_h = \frac{  \{ h, h_0\}}{ \rho^2} \grad h_0  + \frac{1}{\rho} \tilde \cX_h \ , 
    \end{equation}
where $h_0$ is defined in \eqref{h0}.
\end{lemma}
\begin{proof}
First note that $\cX_h$ is homogeneous of degree $-1$.
Hence the vector field 
\begin{equation}\label{tilde_Xh}
\tilde \cX_h := \rho \Big(\cX_h - \frac{\{h, h_0\}}{\rho^2} \grad h_0\Big)
\end{equation}
is homogeneous of order zero. We show that $\tilde \cX_h|_{Z_{\mathbcal{e}_0}^h} \in TZ_{\mathbcal{e}_0}^h$. In view of \eqref{Zwo}, it is sufficient to show that, at each point of $Z_{\mathbcal{e}_0}^h$, $\tilde \cX_h$ is tangent to both $\Sigma_{\mathbcal{e}_0}^h$ and $\S^3$. 
It is tangent to $\Sigma_{\mathbcal{e}_0}^h$ since 
$$
\di h[\tilde{\cX}_h] \stackrel{\eqref{tilde_Xh}}{=} \rho \di h[\cX_h]- \frac{\{h, h_0\}}{\rho} \di h[\nabla h_0]\stackrel{\eqref{gradh.ho}}{=}0.
$$
Next we show that 
the vector field $\tilde \cX_h$ is tangent to 
 $\S^3\equiv \{ (x, \xi) \colon h_0(x, \xi) = \frac12 \}$.
Indeed, as $\rho^2 = 2 h_0$, we also have
$$
\di h_0[\tilde \cX_h] \stackrel{\eqref{tilde_Xh}}{=} 
\rho \Big( \di h_0[\cX_h] - \frac{\{ h, h_0\}}{\rho^2} \underbrace{\di h_0 [\grad h_0]}_{ = 2 h_0} \Big)
=
\rho
\big( \di h_0[\cX_h] -  \{ h, h_0\} \big) {=} 0
$$
proving \eqref{Xh}.
\end{proof}

\begin{remark}\label{dirhoX}
    Similarly it follows that   
 $   \di \rho[\tilde \cX_h] = 0 $ and $
    \di \rho[\cX_h] =\frac{ \{h, h_0\}}{\rho} $.
\end{remark}



It is convenient to decompose also the Liouville measure on the energy level $\Sigma_{\mathbcal{e}_0}^h$ in a ``radial part'' and a component supported on  $Z_{\mathbcal{e}_0}^h$.
So let us denote by  $\nu_0$  the Liouville measure 
induced on the energy level $\Sigma_{\mathbcal{e}_0}^h$ (see \eqref{Sigma0}) by the standard Liouville measure $\nu$ on $\R^4$; it is given explicitly by 
\begin{equation}\label{muSigma0}
    \nu_0 := i_V \nu_{\vert_{\Sigma_{\mathbcal{e}_0}^h}} \ , \ \ \  \mbox{ where } \quad V:= \frac{\grad h}{| \grad h |} \ 
\end{equation}
and $i$ is the contraction operator on differential forms. 
We have
\begin{lemma}\label{lemma.contrazioni}
The Liouville measure $\nu_0$ in \eqref{muSigma0} can be written as 
\begin{equation}\label{muS.dec}
    \nu_0 = \rho^3 \di \rho \wedge \tilde \mu ,
\end{equation}
where $\tilde \mu$ is the two dimensional volume form on $Z_{\mathbcal{e}_0}^h$ defined by
\begin{equation}\label{tildemu}
\tilde \mu := -i_V \nu_{\S^3} \ ,
\end{equation}
with $V$ in \eqref{muSigma0} and $\nu_{\S^3}$ the spherical measure on $\S^3$.
 One has
 \begin{equation}\label{tildemu.prop}
 i_{\grad h_0} \tilde \mu  = 0 \ .
 \end{equation}
\end{lemma}
\begin{proof}
In spherical coordinates $\nu = \rho^3 \di \rho \wedge \nu_{\S^3}$.
Then using \eqref{muSigma0} and 
$
\di \rho[V]  = \frac{1}{|\grad h|}\di \rho[\grad h] = 0$
(which follows from  \eqref{gradh.ho}), we obtain  \eqref{muS.dec} with $\tilde{\mu}$ in \eqref{tildemu}. \\
We now show \eqref{tildemu.prop}. From the definition of $\tilde \mu$ and since $i_X\, i_Y = -i_Y \, i_X$, we have
$i_{\grad h_0}\tilde \mu = i_V \, i_{\grad h_0} \,  \nu_{\S^3} = 0$, concluding the proof.
\end{proof}

We denote by $\mathscr{F}^h$ the 1-dimensional,  oriented, singular  foliation generated by $\tilde \cX_h$ on $Z_{\mathbcal{e}_0}^h$, namely
\begin{equation}\label{def:cF}
    \mathscr{F}^h:= \{ X \in T Z_{\mathbcal{e}_0}^h\ |\ \exists \lambda\in C^\infty(Z_{\mathbcal{e}_0}^h; \R^+) \ \ \text{s.t.} \ \ X = \lambda \tilde\cX_h \} \,. 
 \end{equation} 
It turns out that the existence of an  escape function (according to Definition 
\ref{def:escape}) is guaranteed by certain 
 dynamical conditions on the foliation $\mathscr{F}^h$ that we now describe. 
 Following \cite{Colin_de_Verdi_re_2020}, we now recall the definitions of simple structure and of weak hyperbolicity (recall also Definitions \ref{def:attractor} of attractor, repellor, and of their basins of attraction). 
\begin{definition}[Simple structure]\label{def_simple_structure}
An oriented one dimensional singular foliation $\mathscr{F}^h$ of the compact manifold $Z_{\mathbcal{e}_0}^h$ admits a simple structure $(K^+, K^-)$ if $Z_{\mathbcal{e}_0}^h=K^+ \cup K^- \cup G$ as a disjoint union where: 
\begin{itemize}
    \item[-] $K^+$ is a compact  attractor of the oriented foliation $\mathscr{F}^h$;
    \item[-] $K^-$ is a compact repellor of the oriented foliation $\mathscr{F}^h$;
    \item[-] all leaves of points in $G$ converge to $K^+$ for positive times and to $K^-$ for negative times. In particular the basin of attraction of $K^+$ is $G \cup K^+$ and the basin of attraction of $K^-$ for the reversed orientation of $\mathscr{F}^h$ is $G \cup K^-$. 
\end{itemize} 
\end{definition}

We introduce a further property of an invariant set, namely weak hyperbolicity.
In order to do this, we recall that, given a volume form $\mu$ on a manifold $M$ and a vector field $X$ on {$TM$}, the divergence 
$\div_\mu (X)$ is the function defined by 
\begin{equation}\label{def.div}
    \div_\mu(X) \mu = \di (i_X \mu) \ . 
\end{equation}
We also recall the following standard formulas: if $f$ is a smooth function {on $M$}, 
 \begin{equation}\label{proprietà_divergenza}
 \textup{div}_{\mu}(fX)=\di f(X)+ f \textup{div}_{\mu}(X)  \ , 
 \qquad 
       \textup{div}_{f \mu}(X)=\frac{\di f(X)}{f}+\textup{div}_{\mu}(X) \ . 
       \end{equation}

\begin{definition}[Weakly hyperbolic sets]\label{def:wh}
Let  $(K^+, K^-)$ be a simple structure of a one dimensional singular foliation $\mathscr{F}^h$ on $Z_{\mathbcal{e}_0}^h$. 
The invariant sets $K^\pm $ are called weakly hyperbolic if there exist  neighbourhoods $U^\pm\supset K^\pm$ and a smooth density $\mu$, absolutely continuous with respect to $\tilde{\mu}$  (i.e. a density of the form $\tilde f \tilde{\mu}$, with $\tilde f$ smooth and strictly positive function defined on $Z_{\mathbcal{e}_0}^h$), such that
 $\div_{\mu}(\tilde \cX_h)<0$ on $U^+$ respectively
 $\div_{\mu}(\tilde \cX_h)>0$ on  $U^-$.
 Here $\tilde \cX_h$ is the vector field in \eqref{Xh} and $\tilde \mu$ the measure on $Z_{\mathbcal{e}_0}^h$ in \eqref{tildemu}.
\label{weak_hyp}
\end{definition}

\noindent  We are now able to state the main result of this section. 

\begin{proposition}
    Let $h:\R^4\setminus\{0\} \to \R$ be a positively homogeneous function of degree 0 (according to  Definition \ref{def:pos_hom}).
        If $\mathbcal{e}_0 \in \R$ is a regular value for $h$ such that the foliation $\mathscr{F}^h$ induced by $\tilde{\cX}_h$  on $Z_{\mathbcal{e}_0}^h$ (see \eqref{Zwo} and \eqref{Xh}) has simple structure $(K^+, K^-)$ with $K^+$ and $K^-$ weakly hyperbolic, then there exists an escape function  for $h$ at $\mathbcal{e}_0$ (see Definition \ref{def:escape}).
\label{simple_structure_thm}
\end{proposition}

Let us 
give now an example of a definitely homogeneous function having a weakly hyperbolic simple structure; see Appendix  \ref{app:example} for the proof.
\begin{lemma}\label{lem.ms.ex}
The definitely homogeneous function $h_\star: \R^4\setminus\{0 \} \to \R$ defined by 
\begin{equation}
    h_\star(x, \xi) := \frac12 \frac{x_1^2}{h_0(x, \xi)}
    \end{equation}
    induces on the surface $Z_{\frac12}^{h_\star}$ a foliation $\mathscr{F}^{h_\star}$ with weakly hyperbolic simple structure.
\end{lemma}

The rest of the  section is devoted to the proof of Proposition \ref{simple_structure_thm}.
The construction is rather involved and  follows the following  steps: 
\begin{enumerate}
    \item we construct local escape functions for $h$ in small neighbourhoods of $\Gamma^+:= \boldsymbol{\pi}^{-1}(K^+)$ and $\Gamma^-:= \boldsymbol{\pi}^{-1}(K^-)$;
    \item we extend the local escape functions to the basins of attraction $B(\Gamma^\pm)$ of $\Gamma^\pm$ (see \eqref{def.basins}), and we show that their union covers the entire $\Sigma_{\mathbcal{e}_0}^h$; 
    \item we glue together the local escape functions defined on $B(\Gamma^\pm)$ obtaining a global escape function on $\Sigma_{\mathbcal{e}_0}^h$, which we finally extend to the whole $\R^4 \setminus \{0\}$.
\end{enumerate}

Since the function $h$ is fixed, from now on we shall simply denote 
$\Sigma_{\mathbcal{e}_0}:= \Sigma_{\mathbcal{e}_0}^h$, $Z_{\mathbcal{e}_0}:=Z_{\mathbcal{e}_0}^h$ and $\mathscr{F}:=\mathscr{F}^h$.

\subsection{Construction of a local escape function }

The goal of this section is to construct  a local escape function for $h$ in a neighbourhood of  
$$\Gamma^+:=\boldsymbol{\pi}^{-1}(K^+) \qquad \mbox{ and } \qquad \Gamma^-:=\boldsymbol{\pi}^{-1}(K^-).$$
In all this section, we shall denote by $U^{\pm}$ the open sets defined in Definition \ref{def:wh} of weak hyperbolicity. First we write the Hamilton equations of $h$ in the "radial" variables $(\rho, \zeta) \in \R^+ \times Z_{\mathbcal{e}_0}$: 
\begin{lemma}\label{ham_field_in_rho.zeta}
The variables $(\rho, \zeta) \in \R^+ \times Z_{\mathbcal{e}_0}$ defined in \eqref{projection} satisfy
\begin{equation}\label{eq_rho_zeta}
    \begin{pmatrix} \dot{x} \\ \dot{\xi} \end{pmatrix} =\cX_h(x, \xi)
    \ \ \Leftrightarrow \ \ 
    \begin{cases}
    \dot{\rho}=\dfrac{\{h, h_0\}}{\rho}\\
    \dot{\zeta}=\dfrac{1}{\rho^2}\tilde{\cX}_h(\zeta).
\end{cases}
    \end{equation}
\end{lemma}

\begin{proof}
    We compute the time derivatives of $\rho$ and $\zeta$:
   \begin{equation}\label{dotrho}
        \dot{\rho}= \frac{1}{\rho}\langle \nabla h_0, (\dot{x}, \dot{\xi}) \rangle =\frac{1}{\rho}\di h_0[\cX_h] =\frac{1}{\rho}\{h, h_0\}.
   \end{equation}
   Similarly 
  \begin{equation}
        \dot{\zeta}= \frac{\rho (\dot{x}, \dot{\xi})-\dot{\rho}(x, \xi)}{\rho^2}\stackrel{\eqref{dotrho}}{=}\frac{1}{\rho^2}\Big(\rho \cX_h-\frac{\{h, h_0\}}{\rho}\grad h_0 \Big)\stackrel{\eqref{tilde_Xh}}{=}\frac{1}{\rho^2}\tilde{\cX}_h.
    \end{equation} 
\end{proof}

\noindent Next we also prove the following identity. 
\begin{lemma}
    One has
\begin{equation}\label{ident}
    2\{ h, h_0 \} +  \, \rho\, \div_{\tilde \mu} (\tilde \cX_h) = 0  \quad \mbox{ on } \Sigma_{\mathbcal{e}_0} \ .
\end{equation}
\end{lemma}
\begin{proof}
We start from the identity 
$ \div_{\nu_0 } \cX_h = 0$ on $\Sigma_{\mathbcal{e}_0}$, 
which follows from the fact that the Hamiltonian flow of $h$ preserves the Liouville measure. 
We have
\begin{align}
\div_{\nu_0  } (\cX_h) \, \nu_0   & \stackrel{\eqref{def.div}}{=} 
\di ( i_{\cX_h} \nu_0  )  =
 \di \big( 
\rho^2 \{ h, h_0 \}  \tilde \mu - \rho^2 \di \rho \wedge i_{\tilde \cX_h} \tilde \mu
\big)\,,
\end{align}
where in the last equality we used \eqref{muS.dec}, Remark \ref{dirhoX} and 
$$
i_{\cX_h} \tilde \mu \stackrel{\eqref{Xh}}{=} \frac{\{ h, h_0\}}{\rho^2} i_{\grad h_0} \tilde \mu  + \frac{1}{\rho} i_{\tilde \cX_h} \tilde \mu \stackrel{\eqref{tildemu.prop}}{=} \frac{1}{\rho} i_{\tilde \cX_h} \tilde \mu \ . 
$$
Finally, taking the external differential 
{and observing that, since $\{h, h_0\}$ is an homogeneous function of degree zero and $\tilde{\mu}$ is a volume form on $Z_{e_0}^h$, then, by \eqref{gradh.ho}, $\di \{h, h_0\}\wedge \tilde \mu = 0$, we obtain}
$$
\div_{\nu_0  } (\cX_h) \, \nu_0   = 
2 \rho \, \{h, h_0\} \, \di \rho \wedge \tilde \mu  +
\rho^2 \di \rho \wedge  \di (i_{\tilde \cX_h} \tilde \mu) \stackrel{\eqref{muS.dec}\eqref{def.div}}{=} 
\big( 2  \, \{h, h_0\}\rho^{-2}  +
\rho^{-1} \, \div_{\tilde \mu}(\tilde \cX_h) \big) \nu_0  \ .
$$
As $\div_{\nu_0}( \cX_h) = 0$, we get \eqref{ident}.
Remark that $\div_{\tilde \mu}(\tilde \cX_h)$ is a function homogeneous of order $-1$.
\end{proof}
\begin{remark}
    Identity \eqref{ident} on $Z_{\mathbcal{e}_0}$ reads 
     $2\{ h, h_0 \} +   \div_{\tilde \mu} (\tilde \cX_h) = 0$.
\end{remark}

We begin now the construction of the global escape function for $h$. 
The first step is to construct two functions $k^{\pm}$, which are local escape functions in neighbourhoods of  $\Gamma^\pm.$
\begin{proposition}\label{localK}
Assume that the foliation $\mathscr{F}$ has a simple structure $(K^+, K^-)$ with $K^\pm$ weakly hyperbolic and let $\tilde f: Z_{\mathbcal{e}_0} \to \R$ be the strictly positive, smooth function given by Definition \ref{def:wh}.
Define the functions 
\begin{equation}\label{def:kpm}
  k^\pm: \Sigma_{\mathbcal{e}_0} \to \R, \qquad k^\pm:=\pm\frac{h_0}{\tilde f}, 
\end{equation}
then there exists $\delta >0$ such that
\begin{equation}\label{hh0}
    \{ h, k^\pm\} \geq \delta > 0 \ \ \mbox{ on } \ \ \boldsymbol{\pi}^{-1}(U^\pm).   \,
\end{equation} 
\end{proposition}

\begin{proof}
 We prove the statement with ``$+$''.
From weak hyperbolicity of $K^+$ one has 
\begin{equation}\label{divK+}
\div_{\tilde f \tilde{\mu}}(\tilde{\cX}_h)<0 \mbox{ in } \ U^+   \ .
\end{equation}
We compute: 
\begin{equation}\label{h,k+}
\{h, k^+\}=\di k^+[\cX_h]=\di \Big(\frac{h_0}{\tilde f}\Big)[\cX_h]=\frac{1}{\tilde f}\di h_0[\cX_h]-\frac{h_0}{\tilde f^2}\di \tilde f[\cX_h].
\end{equation}
Evaluating separately the term $\di \tilde f[\cX_h]$ we have: 
\begin{equation}\label{dF.Xh}
\di \tilde f[\cX_h] \stackrel{\eqref{Xh}}{=} \frac{\{h, h_0\}}{\rho^2}\di \tilde f[\nabla h_0]+\frac{1}{\rho} \di \tilde f[\tilde{\cX}_h]=\frac{1}{\rho}\di \tilde f[\tilde{\cX}_h],
\end{equation}
where the last inequality holds by \eqref{gradh.ho}, recalling that $\tilde f$ is homogeneous of degree $0$. 
Inserting this last identity in \eqref{h,k+} and using $h_0 = \frac12 \rho^2$ gives
\begin{equation*}
    \{h, k^+\}\stackrel{\eqref{dF.Xh}}{=}\frac{1}{\tilde f}\{h, h_0\}-\frac{h_0}{\rho}\frac{\di \tilde f[\tilde{\cX}_h]}{\tilde f^2}=\frac{1}{\tilde f}\Big[\{h, h_0\}-\frac{\rho}{2}\frac{\di \tilde f[\tilde{\cX}_h]}{\tilde f}\Big].
\end{equation*}
By identity \eqref{ident}, $\{ h, h_0 \}  = -\frac{\rho}{2} \div_{\tilde \mu}(\tilde \cX_h)$ so we obtain 
\begin{equation*}
    \{h, k^+\}=-\frac{\rho}{2\tilde f}\Big[\frac{\di \tilde f[\tilde{\cX}_h]}{\tilde f}+\div_{\tilde{\mu}}\tilde{\cX}_h\Big]
    \stackrel{\eqref{proprietà_divergenza}}{=}-\frac{\rho}{2 \tilde f}\div_{\tilde f \tilde{\mu}}(\tilde{\cX}_h).
\end{equation*}
Then, by \eqref{divK+} and since the function 
$\frac{\rho}{2 \tilde f}\div_{\tilde f \tilde{\mu}}(\tilde{\cX}_h)$ is homogeneous of degree 0, 
there exists  a strictly positive $\delta $ such that 
$$
\{ h, k^+\} = -\frac{\rho}{2\tilde f}\div_{\tilde f \tilde{\mu}}(\tilde{\cX}_h) \geq \delta >0   \quad \mbox{ on } \ \ \boldsymbol{\pi}^{-1}(U^+)  \ ,
$$
proving \eqref{hh0}. 
The statement with ``$-$'' is analogous.
\end{proof} 
\subsection{Extension to the basins of attraction of $\Gamma^{\pm}$}

The next goal is to extend the local escape functions of Proposition \ref{localK} to the whole {\em basins of attraction} 
 $B(\Gamma^\pm)$ of $\Gamma^\pm = \boldsymbol{\pi}^{-1}(K^\pm)$ defined as 
 \begin{equation}\label{def.basins}
     B(\Gamma^\pm) := \{z=(x, \xi) \in \Sigma_{\mathbcal{e}_0} \colon  \ \textup{dist}\left(\Phi^t_{\cX_h}(z), \Gamma^\pm \right) \to 0 \mbox{ as } t \to \pm \infty \} \ .
 \end{equation}
 We first characterize $B(\Gamma^\pm)$ and prove several properties of these sets. 

\begin{proposition}\label{global_flow}
One has $B(\Gamma^+) \cup B(\Gamma^-) = \Sigma_{\mathbcal{e}_0}$. Precisely:
\begin{itemize}
    \item[(i)] There exist open sets $U^\pm \subset Z_{\mathbcal{e}_0}$, $U^{\pm} \supset K^{\pm},$ such that for any $z_0 \in \boldsymbol{\pi}^{-1}(U^\pm)$ the flow $\Phi^t_{\cX_h}(z_0)$ is globally defined for $\pm t \geq 0$ and 
 $\Phi^t_{\cX_h}(\boldsymbol{\pi}^{-1}(U^\pm)) \subseteq \boldsymbol{\pi}^{-1}(U^\pm)$ for $\pm t \geq 0$.
    \item[(ii)] If $z_0 \in \Gamma^+$, the flow $ \Phi^t_{\cX_h}(z_0) $ exists globally positively in time and $ \Phi^t_{\cX_h}(z_0) \subseteq \Gamma^+$ for any $t \geq 0$, whereas there exists a time $T^+(z_0)<0$ such that 
    $$\Phi^t_{\cX_h}(z_0)  \xrightarrow[t \searrow T^+(z_0)] {} 0;$$
    \item[(iii)] 
    If $z_0 \in \Gamma^-$, 
    the flow $ \Phi^t_{\cX_h}(z_0) $ exists globally negatively in time and $ \Phi^t_{\cX_h}(z_0) \subseteq \Gamma^-$ for any $t \leq 0$, whereas 
    there exists a time $T^-(z_0)>0$ such that 
    $$\Phi^t_{\cX_h}(z_0)  \xrightarrow[t \nearrow T^-(z_0)] {} 0.$$
    \item[(iv)] If $z_0 \in \Sigma_{\mathbcal{e}_0}\backslash (\Gamma^+ \cup \Gamma^-)$, then the flow $\Phi^t_{\cX_h}(z_0)$ exists globally both positively and negatively in time, and one has 
    $$ \Phi^t_{\cX_h}(z_0)  \xrightarrow[t \to +\infty]{}\Gamma^+, \qquad  \Phi^t_{\cX_h}(z_0)  \xrightarrow[t \to -\infty]{}\Gamma^-;$$
    in particular $z_0 \in B(\Gamma^+) \cap B(\Gamma^-)$.
\end{itemize}
\end{proposition}
The long proof of this Proposition is postponed to Appendix \ref{app:global_flow}.

We shall also use the following lemma.

\begin{lemma}\label{lem:fm}
There exists a  positively homogeneous function $\fm$ of degree zero, defined on $\Sigma_{\mathbcal{e}_0}$, such that 
\begin{equation}\label{fm}
\mbox{$\fm$} = \{ h, k^\pm \} \ \ \mbox{on} \ \ \boldsymbol{\pi}^{-1}(U^\pm)   \ , \qquad  \mbox{ and } \qquad \fm(z) \geq \frac{\delta}{2} \quad \ \forall z \in \Sigma_{\mathbcal{e}_0} \ , 
\end{equation}
where $\delta$ is the positive value fixed by \eqref{hh0} and $U^{\pm}$ are as in Proposition \ref{global_flow}. 
\end{lemma}

\begin{proof}
 Take open sets $V^\pm$ with $V^\pm \supset U^\pm$, $V^\pm \cap U^\mp = \emptyset$ and complete $\{V^+, V^-\}$ to a 
covering of $Z_{\mathbcal{e}_0}$ via finitely many open sets $\{V_i\}_i$. Up to restricting $V^\pm$ one has $\{h, k^\pm\} \geq \frac{\delta}{2}$ on $V^\pm$, with $\delta$ in \eqref{hh0}. 
Consider a smooth, non-negative, partition of the identity
$\varphi^+, \varphi^-$, $\{\varphi_i\}_i$ with the properties
$\textup{supp}(\varphi^\pm) \subset V^\pm$ and $\varphi^\pm \equiv 1$ on 
$U^\pm$, 
$\textup{supp}(\varphi_i) \subset V_i$.
Define the smooth function
$$
\tilde \fm:Z_{\mathbcal{e}_0} \to \R, \qquad \tilde \fm := \{h, k^+ \} \varphi^+ + \{h, k^- \} \varphi^- + \delta \sum_i \varphi_i \,.
$$
One has 
 $\tilde \fm \equiv  \{h, k^\pm\} \geq \frac{\delta}{2}$ on $V^\pm$ and thus  $\tilde \fm \geq \frac{\delta}{2}$ on $Z_{\mathbcal{e}_0}$. 
 Then  extend $\tilde \fm$   by homogeneity to a function $\fm$  on $\Sigma_{\mathbcal{e}_0}$,  namely put  $\fm(z):= \tilde \fm(\boldsymbol{\pi}(z))$. 
\end{proof}
We are finally ready to  extend the  local escape functions $k^\pm$ in \eqref{localK} to escape functions on the whole  basins of attraction
 $B(\Gamma^\pm)$.
\begin{lemma}\label{ell.pm}
For $z \in  B(\Gamma^\pm)$, define 
\begin{equation}\label{def:ellpm}
\ell^\pm(z):=\lim_{t \to \pm \infty}\Big( k^\pm(\Phi_{X_h}^t(z))-\int_0^t  \fm(\Phi_{X_h}^s(z))\di s \Big) \ .
\end{equation}
Then
\begin{itemize}
    \item[(i)] the functions $\ell^\pm$ are  well defined on $B(\Gamma^\pm)$, smooth and  positively homogeneous of degree two;
    \item[(ii)] $\ell^\pm$ are  escape functions on $B(\Gamma^\pm)$, namely there exists $\delta >0$ such that 
    \begin{equation}\label{ellpm.esc}
    \{ h, \ell^\pm \} \geq \delta >0  \ \quad \mbox{ in } \ B(\Gamma^\pm) \ . 
    \end{equation}
\end{itemize}
\end{lemma}

\begin{proof}
We prove the statement with ``+'', the other being analogous.\\
$(i)$ \underline{Well defined:} 
Let $z \in B(\Gamma^+)$ and let $t_0(z)>0$ be such that $\Phi_{X_h}^{t_0(z)}(z) \in \boldsymbol{\pi}^{-1}(U^+)$, with $U^+$ as in Proposition \ref{global_flow}. By definition of $B(\Gamma^+)$, $t_0(z)$ is finite.  
We show that for any  $t>t_0(z)$, the argument of the limit in \eqref{def:ellpm} is constant.  
Since,  by  Proposition \ref{global_flow} $(i)$, the set $\boldsymbol{\pi}^{-1}(U^+)$ is positively invariant, by \eqref{fm} we have
$\fm(\Phi_{X_h}^t(z))= \{ h, k^+\}(\Phi_{X_h}^t(z))$  for any  $ t\geq t_0(z)$.
Therefore for $t\geq t_0(z)$ one has
\begin{align}
    \notag
        & k^+(\Phi_{X_h}^t(z))-\int_{0}^t \fm(\Phi_{X_h}^s(z))\di s
        \notag
        = k^+(\Phi_{X_h}^t(z))-\int_0^{t_0(z)} \fm(\Phi_{X_h}^s(z))\di s-\int_{t_0(z)}^t \{h, k^+\}(\Phi_{X_h}^s(z))\di s\\
        \notag
        = & k^+(\Phi_{X_h}^t(z))-\int_0^{t_0(z)} \fm(\Phi_{X_h}^s(z))\di s- k^+(\Phi_{X_h}^t(z))+ k^+(\Phi_{X_h}^{t_0(z)}(z))\\
        \label{ell2}
        = & -\int_0^{t_0(z)} \fm(\Phi_{X_h}^s(z))\di s + k^+(\Phi_{X_h}^{t_0(z)}(z)),
    \end{align}
    and thus the argument of the limit is definitely constant, proving that $\ell^+$ is well defined. 

\noindent \underline{Smoothness:} 
fix $z_0 \in B(\Gamma^+)$ and let $V_{z_0} \subset B(\Gamma^+)$ be a neighbourhood of $z_0$. We observe that, again by Proposition \ref{global_flow}, there exists $t_1(V_{z_0})>0,$ eventually larger than $t_0(z_0),$ such that $ \Phi^{t}_{\cX_h}(z) \in \boldsymbol{\pi}^{-1}(U^+)$ for all $z \in V_{z_0}$ and $t \geq t_1(V_{z_0})$. 
Then,  by \eqref{ell2}, on the set $V_{z_0}$ the function $\ell^+$ reduces to
\begin{equation}\label{ell.piu.v}
\ell^+(z)=  -\int_0^{t_1(V_{z_0})} \fm(\Phi_{X_h}^s(z))\di s + k^+(\Phi_{X_h}^{t_1(V_{z_0})}(z)) , \quad \forall z \in V_{z_0}  \ ,
\end{equation}
which is smooth.

\noindent \underline{Homogeneity:}
It follows using the homogeneities of $k^+$ and $\fm$ and 
Remark \ref{homo.flow}.

\noindent $(ii)$
Given $z_0 \in B(\Gamma^+)$ and $V_{z_0} \subset B(\Gamma^+)$ a small neighbourhood of $z_0$, let $t_1(V_{z_0})$ as above. We have
\begin{equation*}
    \begin{split}
        \{h, \ell^+\}(z) 
        &\stackrel{\eqref{ell.piu.v}}{=} \{h, k^+ \circ \Phi_{X_h}^{t_1(V_{z_0})}\}(z)-
        \int_0^{t_1(V_{z_0})} \{ h, \fm \circ \Phi_{X_h}^s(z)) \} \di s\\
        &= \{h, k^+  \}(\Phi_{X_h}^{t_1(V_{z_0})}(z))
        -\frac{\di}{\di t}\Big|_{t=0}\int_0^{t_1(V_{z_0})} \fm(\Phi_{X_h}^{t+s}(z)) \di s  \\
        &= \{h, k^+  \}(\Phi_{X_h}^{t_1(V_{z_0})}(z)) -\fm(\Phi_{X_h}^{t_1(V_{z_0})}(z))+ \fm(z)\\
        &{=} \fm(z) {\geq} \frac{\delta}{2}  \ ,
    \end{split}
\end{equation*}
where in the last two passages we used \eqref{fm}. {Finally rename $\delta \leadsto 2\delta$.}
\end{proof}

\subsection{Construction of the  global escape function}
In this section we  glue together the local escape functions $\ell^+$ and $\ell^-$ of Lemma \ref{ell.pm}, defined in $B(\Gamma^+)$ and $B(\Gamma^-)$ respectively,  to obtain a global escape function on $\Sigma_{\mathbcal{e}_0}$.
We start by proving the following result.

\begin{proposition}\label{prop:glue}
    Given  $\varepsilon >0$ sufficiently small, there exists a smooth function 
    $\eta_\varepsilon \colon \Sigma_{\mathbcal{e}_0} \to \R$ such that
    \begin{itemize}
        \item[(i)] $\eta_{\epsilon}$ is  homogeneous of degree $0$;
        \item[(ii)] $\eta_{\epsilon}(\Gamma^+) \equiv 1$, $\eta_{\epsilon}(\Gamma^-) \equiv 0$;
        \item[(iii)] $|\di \eta_{\epsilon}(\zeta)[\tilde \cX_h(\zeta)] | \leq \varepsilon$ for any $\zeta \in Z_{\mathbcal{e}_0}$.
    \end{itemize}
\end{proposition}
 To prove this proposition we shall use the following result,  whose proof is postponed  to  Appendix \ref{app:V-}.
\begin{lemma}\label{prop:V-}
There exist a smooth function $\tilde{\cL}:Z_{\mathbcal{e}_0} \rightarrow \R$ and an interval $(\alpha,\ + \infty) \subset \R$ such that the open set $V^-:= \tilde{\cL}^{-1}((\alpha, \ + \infty))$ has the following properties: 
\begin{itemize}
    \item[$(i)$] $V^- \supset K^-$ is invariant by the flow of $\tilde{\cX}_h$ for negative times, and $Z_{\mathbcal{e}_0}\setminus V^- \supset K^+$ is invariant by the flow of $\tilde \cX_h$ for positive times; 
    \item[$(ii)$] $\partial V^-$ is   transverse to the flow of $\tilde{X}_h$, in particular  $\tilde{\cX}_h[\tilde{\cL}]<0$ on $\partial V^-$. 
    \end{itemize}
\end{lemma}
We use this result to define the gluing function. Take  $\zeta \in Z_{\mathbcal{e}_0} \setminus (K^+ \cup K^-)$. Since 
$\Phi_{\tilde{X}_h}^{t}(\zeta) \to  K^\pm $  when  $t \to \pm \infty$,
 by Lemma \ref{prop:V-} there will be a unique time $t_{\zeta} \in \R$ such that $\Phi_{\tilde{X}_h}^{-t_\zeta}(\zeta)$ intersects $\partial V^-$; we then put 
\begin{equation}
\tilde t(\zeta):=
\begin{cases}
t_{\zeta} , &  \zeta \in  Z_{\mathbcal{e}_0} \setminus (K^+ \cup K^-)  \\
+\infty, &  \zeta  \in K^+ \\
-\infty, & \zeta \in K^-.\\
\end{cases}
\label{tildet(z)}
\end{equation}
Next, given $\varepsilon >0$ sufficiently small,  take a smooth non-decreasing function  $ \phi_\varepsilon\colon \R \to \R$ such that 
\begin{equation}\label{def:phivare}
{\phi}_\varepsilon(\tau):=
\begin{cases}
0, & \tau\leq -\varepsilon  \\
\varepsilon \tau , & \varepsilon <\tau \leq \varepsilon^{-1}- \varepsilon  \\
1,  & \tau> \varepsilon^{-1} + \varepsilon\\
\end{cases} \ 
\qquad
\mbox{ and }
\qquad 
\sup_{\tau \in \R} |\phi'_\varepsilon(\tau)  | \leq 10 \epsilon \ .
\end{equation}
The gluing function of Proposition \ref{prop:glue} is
\begin{equation}\label{def:glue}    \eta_\varepsilon(z):=\phi_{\frac{\varepsilon}{10}}(\tilde t(\boldsymbol{\pi}(z))) \ , \quad z \in \Sigma_{\mathbcal{e}_0} \ ,
\end{equation}
where $\tilde t$ is defined in \eqref{tildet(z)} and $\phi_{\epsilon}$ in \eqref{def:phivare}.

\begin{proof}[Proof of Proposition \ref{prop:glue}]
It is clear from its definition that $\eta_\varepsilon$ is homogeneous of degree $0$.
Moreover from \eqref{tildet(z)}, \eqref{def:phivare} we have 
 $
 \eta_\varepsilon(z) 
 := \phi_{\frac{\varepsilon}{10}}(\tilde t(\boldsymbol{\pi}(z))) = 1$  for any $z \in \Gamma^+$;  
similarly $\eta_\varepsilon(\Gamma^-) \equiv 0$. This proves items $(i)$ and $(ii)$.

We prove now  item $(iii)$. Let $\zeta \in Z_{\mathbcal{e}_0}$ and note that, 
by \eqref{tildet(z)}, for any $\tau \in \R$ one has 
\begin{equation}\label{time.tilde}
 \tilde t( \Phi^{\tau}_{\tilde{X}_h}(\zeta) ) = \tilde t(\zeta) + \tau \ ;   
\end{equation}
hence we compute 
 \[
\di \eta_\varepsilon(\zeta)[\tilde{X}_h(\zeta)]
=\frac{\di}{\di \tau}\Big|_{\tau=0}\eta_\varepsilon\Big(\Phi^{\tau}_{\tilde{X}_h}(\zeta)\Big)
\stackrel{\eqref{time.tilde}}{=}\frac{\di}{\di \tau}\Big|_{\tau=0}\phi_{\frac{\varepsilon}{10}}(\tilde t(\zeta)+\tau)
= \phi'_{\frac{\varepsilon}{10}}(\tilde t(\zeta)) \ ,
\]
and $(iii)$ follows from \eqref{def:phivare}.

We are left to prove that $ \eta_\varepsilon$ is smooth. Remark that, due to the  homogeneity of $ \eta_\varepsilon$, we can assume $z=:\zeta \in Z_{\mathbcal{e}_0}.$
We are going to prove that $\eta_\varepsilon$ is constant
on sufficiently small  neighbourhoods of $K^{\pm}$,  whereas outside such neighbourhoods we will prove smoothness of $\tilde{t}$ (and thus of $\eta_\varepsilon$) by implicit function theorem.
Take a point $\underline{\zeta} \in Z_{\mathbcal{e}_0} \setminus (K^+ \cup K^-)$ and consider the corresponding point
$
\underline{y}=\Phi_{\tilde{X}_h}^{-\tilde t(\underline{\zeta})}(\underline \zeta) \in \partial V^-,$
with $V^-$ the set of Lemma \ref{prop:V-}.
Recall that $\partial V^-= \tilde{\cL}^{-1}(\alpha)$. 
Then, let $U_{\underline \zeta}$ be a small neighbourhood of $\underline \zeta$ with 
$U_{\underline \zeta} \cap (K^+ \cup K^-) = \emptyset$:
 the function 
\begin{equation}
    F: U_{\underline \zeta} \times \R \to \R  \ , \quad 
    F(\zeta, t):=\tilde{\cL}(\Phi_{\tilde{X}_h}^{-t}(\zeta)) - \alpha
\end{equation}
 is smooth, $F(\underline{\zeta}, \tilde t(\underline{\zeta}))=0$ by the very definition of $\tilde t$ and, being $\partial V^-$ transverse to $\tilde{X}_h$ by Lemma \ref{prop:V-} $(iii)$, 
\begin{equation}
    \partial_t {F}(\underline \zeta, \tilde t(\underline{\zeta}))=
    - \di \tilde{\cL}(\underline y)[\tilde{X}_h(\underline y)]  >  0 \ .
\end{equation}
Then, by implicit function theorem, in 
$U_{\underline{\zeta}}$ the function $\zeta \mapsto \tilde t(\zeta)$ is smooth. 
We conclude that $\eta_\varepsilon$ is smooth  in every point of $Z_{\mathbcal{e}_0}\setminus (K^+ \cup K^-)$. 

The final claim is that there exist ${W}^\pm$ neighbourhoods of $K^\pm$  such that $\eta_\varepsilon\equiv 0$ in $W^-$ and $\eta_\varepsilon \equiv 1$ in $W^+$. We prove this claim with "+". 
To this aim it suffices to define 
\[
{W}^+:=\Phi^{\varepsilon^{-1}+1}_{\tilde{X}_h}(Z_{\mathbcal{e}_0}\setminus V^-) \ .
\]
First remark that $W^+$ is a neighbourhood of $K^+$ since  $K^+ \subset Z_{\mathbcal{e}_0}\setminus V^-$, and since $K^- \not\subset Z_{\mathbcal{e}_0}\setminus V^-$, one has $ Z_{\mathbcal{e}_0}\setminus V^-$ is contained in the basin of attraction of $K^+$; then  using  \eqref{attractor}
$$
K^+=\bigcap_{t>0}\Phi^t_{\tilde{X}_h}(Z_{\mathbcal{e}_0}\setminus V^- ) \subset \Phi^{\varepsilon^{-1} + 1}_{\tilde{X}_h}(Z_{\mathbcal{e}_0}\setminus V^-) = W^+.
$$
Next, for any   $\zeta \in W^+$, there is  $ \upsilon  \in Z_{\mathbcal{e}_0}\setminus V^-$ such that 
$
\zeta =\Phi^{\varepsilon^{-1}+1}_{\tilde{X}_h}(\upsilon) $ and with $\tilde t(\upsilon) \geq 0$, 
hence 
  $
\tilde t(\zeta)\stackrel{\eqref{time.tilde}}{=} \epsilon^{-1}+1+\tilde t(\upsilon)> \epsilon^{-1} + \epsilon \ ; 
$
 therefore $\eta_\varepsilon(\zeta)=\phi_{\frac{\varepsilon}{10}}(\tilde t(\zeta))=1$ by the properties of $\phi_\varepsilon$ in \eqref{def:phivare}. \\
The  same argument  can be repeated for $K^-$ taking
$W^-:=\Phi^{-1}_{\tilde{X}_h}(V^-)$, 
obtaining that  $\eta_\varepsilon \equiv 0$ in $W^-$. 
\end{proof}

We are now in position to construct a global escape function.
First we construct it on $\Sigma_{\mathbcal{e}_0}$ and then extend it to the whole $\R^4 \setminus \{0\}$. 
Using all previous notation for the local escapes and the gluing functions, we prove the following result:

\begin{proposition}\label{escape.sigma0}
There exists a function $a\colon \R^4\setminus \{0\} \to \R$ positively homogeneous of degree 2, such that $a$  is an escape function for $h$ at $\mathbcal{e}_0$ in the sense of Definition \ref{def:escape}.
\end{proposition}

\begin{proof}
   \underline{Step 1:} First we construct a smooth function $\tilde{a}: \Sigma_{\mathbcal{e}_0} \to \R$, positively homogeneous of degree 2, that satisfies \eqref{Poiss.esc} for some $\delta>0$. 
To this goal, it suffices to define, for $\epsilon$ sufficiently small,
 \begin{equation}\label{def:globalescape}
    \tilde{a}(z):=\eta_\varepsilon(z)\ell^+(z)+(1-\eta_\varepsilon(z))\ell^-(z) \ ,
\end{equation}
with $\eta_\epsilon$ defined in \eqref{def:glue} and $\ell^\pm$ in \eqref{def:ellpm}. 
Remark that $\tilde{a}$ is well defined: indeed, from Lemma \ref{prop:V-},
$\eta_\epsilon \equiv 0$ in  $\Gamma^- $where  $\ell^+$  is not defined, whereas  $ 1-\eta_\varepsilon \equiv 0$ in
$\Gamma^+$  where  $\ell^- $ is not defined.

Next note that $\tilde{a}$ is smooth, since all terms in its expression are smooth and  $\tilde{a}$ is positively homogeneous of degree 2, since both $\ell^+$ and $\ell^-$ are (see Lemma \ref{ell.pm}) and $\eta_{\epsilon}$ is positively homogeneous of degree zero (see Proposition \ref{prop:glue}). 
We now show that $\tilde{a}$ satisfies \eqref{Poiss.esc}. 
One has: 
\begin{equation}\label{PBha}
    \{h, \tilde{a}\}=\{h, \eta_\varepsilon\}\ell^+ +
    \{h, \ell^+\}\eta_\varepsilon
    +\{h, \ell^-\}(1-\eta_\varepsilon)+\{h, (1-\eta_\varepsilon)\}\ell^- \stackrel{\eqref{ellpm.esc}}{\geq}
    \{h, \eta_\varepsilon\}(\ell^+ - \ell^-)+\delta \ .
\end{equation}
Next we remark that, by \eqref{Xh} we have
\begin{equation}
    \{h, \eta_\varepsilon\}=\di \eta_\varepsilon[X_h]
    =\frac{\{h, h_0\}}{\rho^2} \di \eta_\varepsilon[\grad h_0] + \frac{1}{\rho} \di \eta_\varepsilon[\tilde{X}_h]= \frac{1}{\rho}\di \eta_\varepsilon[\tilde{X}_h],
\end{equation}
where the last equality is obtained using  \eqref{gradh.ho} and the homogeneity of $\eta_\varepsilon$. 
Since $\di \eta_\varepsilon[\tilde{X}_h]$ is homogeneous of degree $-1$, using Proposition \ref{prop:glue} $(iii)$ and writing $z  \equiv  \rho \boldsymbol{\pi}(z)\in \Sigma_{\mathbcal{e}_0}$, we get 
\[
    \abs{\left\{h, \eta_\e \right\}(z)} = \left|\frac{1}{\rho}\di \eta_\varepsilon[\tilde{X}_h](z)\right| = 
\left| \frac{1}{\rho^2}\di \eta_\varepsilon[\tilde{X}_h](\boldsymbol{\pi}(z))\right| \leq \frac{\epsilon}{\rho^2} \, . 
\]
In conclusion, using also that $\ell^\pm$ are positively homogeneous functions of degree 2,  we have the following bound:
\[
 |\{h, \eta_\varepsilon\}(z)(\ell^+(z)-\ell^-(z))|
 \leq  \frac{\varepsilon}{\rho^2} \sup_{z \in \Sigma_{\mathbcal{e}_0}}|\ell^+(z)-\ell^-(z)|
 \leq   \epsilon  \sup_{\zeta \in Z_{\mathbcal{e}_0}}|\ell^+(\zeta)-\ell^-(\zeta)| \leq C \varepsilon \ .
\]
Hence, shrinking $\e$,  we obtain from \eqref{PBha} the lower  bound 
\begin{equation}\label{atilde}
    \{h, \tilde{a}\}>\delta - C\varepsilon>\frac{\delta}{2}\quad  \mbox{ on } \Sigma_{\mathbcal{e}_0} \ . 
\end{equation}

\noindent \underline{Step 2:}
We extend $\tilde{a}$ to the whole $\R^4\setminus \{0\}$.
To this aim  consider first the restriction 
$
    \tilde{a}|_{Z_{\mathbcal{e}_0}}: Z_{\mathbcal{e}_0} \to \R.$ 
    It is a smooth function defined on the  closed submanifold $Z_{\mathbcal{e}_0}$ of $\S^3$.
   By classical  extension Lemma (see e.g.  \cite[Lemma 2.26]{lee2003introduction}), given an open neighbourhood  $\cU$  of $Z_{\mathbcal{e}_0}$ in $\S^3$, there exists a smooth function $\hat{a}$,
    $$\hat{a}: \S^3 \to \R, \mbox{ such that } \hat{a}=\tilde{a} \mbox{ on } Z_{\mathbcal{e}_0} \mbox{ and } \mbox{ supp}(\hat{a}) \subset \cU.$$
Then  we define $a$ as the extension of $\hat{a}$ by homogeneity on the whole $\R^4\backslash \{0\}$. Due to the cone structure of $\Sigma_{\mathbcal{e}_0}$ and the fact that $\tilde{a}$ and $a$ are both homogeneous of degree two, $\tilde{a}=a$ on the whole $\Sigma_{\mathbcal{e}_0}$, since they coincide on $Z_{\mathbcal{e}_0}$. 
    In particular $\{h, a\}=\{h, \tilde{a}\} \geq \frac{\delta}{2}$ on $\Sigma_{\mathbcal{e}_0}$.
\end{proof}
\begin{remark}  If $\Sigma_{\mathbcal{e}_0}$ is not connected, one proves {Step 1} of the proof of Proposition \ref{escape.sigma0} separately for each one of its $n$ connected components $\Sigma_j$, obtaining  $n$ functions $\tilde {a}_j: \Sigma_j \rightarrow \R$, homogeneous of degree 2, such that $\{h, \tilde a_j\} > \frac{\delta_j}{2}$ for all $j$ and for some positive $\delta_1, \dots, \delta_n$. Then one defines $\tilde a$ on the whole $\Sigma_{\mathbcal{e}_0}$ as $\tilde a = \tilde a_j$ on $\Sigma_j$; \eqref{atilde} holds with $\delta = \min\{\delta_j\}>0$, and Step 2 follows unchanged. \end{remark}
We conclude this section proving that if $a$ is an escape function for $h$ at $\mathbcal{e}_0$ in the sense of Definition \ref{def:escape}, then $a$ is also an escape function for $h$ at $\omega$, for $\omega$ sufficiently close to $\mathbcal{e}_0$. Precisely we have: 

\begin{corollary}\label{interval.escape}
    Let $\mathbcal{e}_0$ be a regular value of $h$. 
    If $a$ is an escape function for $h$ at $\mathbcal{e}_0$ (see Definition \ref{def:escape}), there exists a neighbourhood $I$ of $\mathbcal{e}_0$, contained in the image of $h$, such that $a$ is still an escape function for $h$ at all values $\omega \in I$. 
\end{corollary}

\begin{proof}
Let $\delta >0$ be such that \eqref{Poiss.esc} is satisfied on $\Sigma_{\mathbcal{e}_0}$; we start by showing that there exists a sufficiently small  interval $I \ni \mathbcal{e}_0$ such that 
\begin{equation}\label{escape.nbh}
\{ h, a \} \geq \frac{\delta}{2}, \qquad \forall \ (x, \xi) \in \bigcup_{\omega \in I} \Sigma_{\omega}. 
\end{equation}
Since the function $\{h, a\}$ is homogeneous of degree zero, it is sufficient to prove \eqref{escape.nbh} for any $ \zeta \in \bigcup_{\omega \in I} Z_{\omega}$,  
with $Z_{\omega}:= h^{-1}(\omega) \cap \S^3$. 
From continuity of $\{h, a\}|_{\S^3}$ and \eqref{Poiss.esc}, there exists an open neighbourhood $\cU_{Z_{\mathbcal{e}_0}} \subset \S^3$ of $Z_{\mathbcal{e}_0}$ such that 
$$ \{h, a\} \geq \frac{\delta}{2} \qquad \forall \ \zeta \in \cU_{Z_{\mathbcal{e}_0}}.$$
We show that, choosing $\epsilon>0$ small enough, the interval $I:=( \mathbcal{e}_0 -\epsilon, \mathbcal{e}_0 +\epsilon)$ fulfills $h^{-1}(I) \subset \cU_{Z_{\mathbcal{e}_0}}$, thus proving \eqref{escape.nbh}.
Suppose by contradiction that $\exists \ N \in \N$, such that, defining 
\begin{equation}\label{def.In}
I_n:=\Big(\mathbcal{e}_0-\frac{1}{n},\ \mathbcal{e}_0 + \frac{1}{n}\Big),
\end{equation}
one has 
\begin{equation}\label{contradiction}
   \{ h^{-1}(I_n) \cap \S^3 \} \not\subset \cU_{Z_{\mathbcal{e}_0}}, \qquad \forall \ n>N.
\end{equation}
Take now a sequence $\{x_n\}_{n \in \N}$ such that $x_n \in \{h^{-1}(I_n) \cap \S^3 \} \ \backslash\ \cU_{Z_{\mathbcal{e}_0}} $.
By compactness, up to a subsequence it converges to a point $ x_{\infty} \in \S^3 \setminus \cU_{Z_{\mathbcal{e}_0}}$.
 From continuity of $h$ we get 
$
h(x_n) \to h(x_{\infty}),$ 
but $I_n \ni h(x_n) \to \mathbcal{e}_0$ from the definition \eqref{def.In} of $I_n$, whilst $x_{\infty} \notin Z_{\mathbcal{e}_0}$, hence $h(x_{\infty}) \neq \mathbcal{e}_0$, obtaining a contradiction.  

Finally remark that, since $\mathbcal{e}_0$ is a regular value of the smooth function $h$, up to shrinking $\epsilon$ we can assume that $I$ is contained in the image of $h$ and only consists of regular values, concluding the proof.
\end{proof}

\section{Genericity of the set $\cV$}

The goal of this section is to show that the set $\cV$ defined in \eqref{setV} is generic in $C^0(\T, \bS^0_{cl})$ with respect to the distance $\td^{0}_{\T}$ defined in \eqref{distance}, proving
Theorem \ref{main_theorem} $(i)$.

\subsection{The set $\cV$ is open}
The main result of this section is the following: 

\begin{proposition}\label{V_open}
    The set $\cV$ defined in \eqref{setV} is open with respect to the distance $\td^{0}_{\T}$ defined in \eqref{distance}.
\end{proposition}
We separately state the following preliminary lemma.

\begin{lemma}\label{small.average}
Let $v \in C^0(\T, \bS_{cl}^0)$ and let $\la v \ra$ be its average along the flow of $h_0$ (see Definition \ref{resonant_average}), then
$$
\sum_{|\alpha|+|\beta| \leq 1}\sup_{|(x, \xi)| \geq 1} {\Big|\pa_x^{\alpha}\pa_{\xi}^{\beta}\la v \ra(x, \xi)\Big|} \leq 4 \sum_{|\alpha|+|\beta| \leq 1}\sup_{\substack{|(x, \xi)| \geq 1 \\ t \in \T}} {\Big|\pa_x^{\alpha}\pa_{\xi}^{\beta}v(t, x, \xi)\Big|}\,.
$$


\end{lemma}

\begin{proof}
We only show the bound for the case $\alpha =(1, 0)$ and $\beta = (0,0)$, all the others being analogous. Using the explicit form  \eqref{h0.flow} of the classical flow  $(x(t), \xi(t)):= \phi^t_{h_0}(x, \xi)$ of  the harmonic oscillator, 
and the definition \eqref{average}, we have
\begin{equation}
\begin{split}
     \sup_{|(x, \xi)| \geq 1} |\pa_{x_1} \la v \ra (x, \xi)| 
    & \leq   \frac{1}{2\pi} \int_0^{2\pi} \sup_{\substack{|(x, \xi)| \geq 1\\ t \in \T}}(|\pa_{x_1}v(t, x(t), \xi(t))| +|\pa_{\xi_1} v(t, x(t), \xi(t))|) \, \di t\\
    & \leq  \sup_{\substack{|(x, \xi)| \geq 1\\ t \in \T}}(|\pa_{x_1}v(t, x, \xi)| +|\pa_{\xi_1} v(t, x, \xi)|)  \ . 
    \end{split}
\end{equation}
Analogously one bounds the other terms with $|\alpha|+|\beta|\leq 1$, obtaining the claim. 
\end{proof}

\begin{proof}[Proof of Proposition \ref{V_open}]
Take $v \in \cV$. We show the existence of  $\epsilon>0$ sufficiently small such that 
\begin{equation}\label{vicino.w}
    \td^{0}_{\T}(v, v+w) < \epsilon \quad \Rightarrow \quad  v+w \in \cV \ 
\end{equation}
for any $w \in C^0(\T,  \bS^0_{cl})$. 
Decompose $v=v_0+v_{-\mu},$ and $ w=w_0+w_{-\mu'}$ following Definition \ref{def:symbol}, where $v_0$ and $w_0$ are the principal symbols of $v$ and $w$ respectively, and $v_{-\mu} \in C^0(\T, \bS^{-\mu})$ and $w_{-\mu'} \in C^0(\T, \bS^{-\mu'})$, with $\mu, \mu'>0$.\\
By the very definition of $\cV$ in \eqref{setV}, since $v \in \cV$ there exist $a \in \bS^1_{cl}$, $I \subset \Ran(\la v_0 \ra)$ and $\delta>0$  such that
\begin{equation}
    \{ \la v_0 \ra, a \} \geq \delta, \qquad \forall \ (x, \xi) \in \la v_0 \ra^{-1}(I) \cap \{|(x, \xi)| \geq 1\}.
\end{equation}
We are going to show that there exist $J \subset \Ran(\la v_0 +w_0 \ra)$ and $\delta'>0$ such that 
\begin{equation}\label{claim}
    \{ \la v_0 + w_0 \ra, a \} \geq \delta', \qquad \forall \ (x, \xi) \in \la v_0 + w_0 \ra^{-1}(J) \cap \{|(x, \xi)| \geq 1\}\,.
\end{equation}
To this aim, we claim  that
\begin{equation}\label{R.seminorm}
    \sum_{|\alpha|+|\beta| \leq 1}\sup_{\substack{|(x, \xi)| \geq 1 \\t \in \T}}\Big|\pa_x^{\alpha}\pa_{\xi}^{\beta} w_0(t, x, \xi)\Big|\leq 10 \epsilon\,.
\end{equation}
Assuming this estimate, we continue with the proof of \eqref{claim}.
We decompose $\{\la v_0 + w_0 \ra, a \} 
    = \{ \la v_0\ra, a\}+\{\la w_0 \ra, a\}\,$ and we observe that, by Lemma \ref{small.average},
    \begin{equation}\label{}
 \sup_{|(x, \xi)| \geq 1}|\{ \la w_0 \ra, a\}| \leq 2 \wp^{1}_1(a)  \sum_{|\alpha|+|\beta| \leq 1}\sup_{|(x, \xi)| \geq 1}\Big|\pa_x^{\alpha}\pa_{\xi}^{\beta} \la w_0 \ra (x, \xi)\Big| \stackrel{\eqref{R.seminorm}}{\leq} 80 \wp^{1}_1(a) \varepsilon\,.
    \end{equation}
   Thus, choosing $\e \leq \frac{\delta}{160 \wp^{1}_1(a)}\,,$ estimate \eqref{claim} follows with $\delta'= \frac{\delta}{2}$ in $\la v_0 \ra^{-1}(I) \cap \{|(x, \xi)| \geq 1\}$.
   In order to conclude the proof of \eqref{claim} we show that there exists $J \subset I$ such that, for all $w$ satisfying the bound in \eqref{vicino.w},
\begin{equation}\label{on.J}
\la v_0+w_0 \ra ^{-1}(J) \cap \{|(x, \xi)| \geq 1\} \subset \la v_0 \ra ^{-1}(I) \cap \{|(x, \xi)| \geq 1\}\,.
\end{equation}
From equation $\eqref{R.seminorm}$ and Lemma \ref{small.average} we have that 
\begin{equation}\label{sup.sm}
\sup_{|(x, \xi)| \geq 1} |\la w_0 \ra (x, \xi) | \leq 4 \sup_{|(x, \xi)| \geq 1 \atop t \in \T} | w_0(t, x, \xi) | \leq 40 \epsilon\,.
\end{equation}
Hence, denoting by $I=: [a, b]$, we choose 
$J : =[ a+ 40 \epsilon, b - 40 \epsilon]$. Up to shrinking $\epsilon$ we have $ \emptyset \neq  J \subset I$, and by \eqref{sup.sm} $J$ satisfies
\eqref{on.J}. Thus \eqref{claim} is proven.\\
It remains to prove \eqref{R.seminorm}.
We first notice that, from the definition \eqref{distance} of the distance, if $\td^{0}_{\T}(v, v+w) < \epsilon$, up to shrinking $\e$ we have

\begin{equation}\label{seminorm0}
    \wp_{\T,1}^{0}(w)=\sum_{|\alpha|+|\beta| \leq 1}\sup_{\substack{(x, \xi) \in \R^4 \\t \in \T}} \frac{\Big|\pa_x^{\alpha}\pa_{\xi}^{\beta} w(t, x, \xi)\Big|}{(1+|x|^2+|\xi|^2)^{-\frac{|\alpha|+|\beta|}{2}}} \leq 4 \epsilon\,. 
\end{equation}

\noindent Recall the decomposition $w=w_0 + w_{-\mu'};$ since $w_{-\mu'} \in C^0(\T, \bS^{-\mu'})$ (see Definition \ref{def:symbol}), we have 
\begin{equation}\label{decay.w}
\frac{\Big|\pa_x^{\alpha}\pa_{\xi}^{\beta} w_{-\mu'}(t, x, \xi)\Big|}{{(1+|x|^2+|\xi|^2)^{-\frac{|\alpha|+|\beta|}{2}}}} \leq \frac{C_{\alpha, \beta}(w)}{(1+|x|^2+|\xi|^2)^{\mu'}}, \qquad \forall \alpha, \beta \in \N^2.  
\end{equation}

\noindent Fix now $\lambda_{\epsilon} = \lambda_{\e}(w)\geq 1$ sufficiently large such that
\begin{equation}\label{bound.tail}
     \sum_{|\alpha| + |\beta| \leq 1} \frac{C_{\alpha, \beta}(w)}{(1+\lambda_{\epsilon}^2)^{\mu'}}< \epsilon.
\end{equation}
We have
\begin{equation}\label{5epsilon}
\begin{aligned}
 \sum_{|\alpha|+|\beta| \leq 1} \sup_{\substack{|(x, \xi)|\geq \lambda_{\epsilon} \\ t \in \T}} \frac{\left|\partial^{\alpha}_{x} \partial^{\beta}_{\xi} w_0(t, x, \xi)\right|}{(1 + |x|^2 + |\xi|^2)^{-\frac {|\alpha| + |\beta|}{2}}} & 
 \leq \sum_{|\alpha|+|\beta| \leq 1} \sup_{\substack{|(x, \xi)|\geq \lambda_{\epsilon} \\ t \in \T}} \frac{\left|\partial^\alpha_x \partial^\beta_\xi w_{-\mu'}(t, x, \xi)\right|}{(1 + |x|^2 + |\xi|^2)^{-\frac {|\alpha| + |\beta|}{2}}} + \wp^{0}_{\T,1}(w)\\
&\stackrel{\eqref{seminorm0}, \eqref{bound.tail}}{\leq} \e + 4 \e = 5 \e\,. 
\end{aligned}
\end{equation}

\noindent Finally remark that 
\begin{equation}\label{quasi.homo}
\frac{\Big|\pa_x^{\alpha}\pa_{\xi}^{\beta} w_0(t, x, \xi) \Big|}{(1+ |x|^2+|\xi|^2)^{-\frac{|\alpha|+|\beta|}{2}}} \leq  2 \frac{\Big|\pa_x^{\alpha}\pa_{\xi}^{\beta} w_0(t, x, \xi) \Big|}{(|x|^2+|\xi|^2)^{-\frac{|\alpha|+|\beta|}{2}}} \leq 2 \frac{\Big|\pa_x^{\alpha}\pa_{\xi}^{\beta} w_0(t, x, \xi) \Big|}{(1+ |x|^2+|\xi|^2)^{-\frac{|\alpha|+|\beta|}{2}}} ,    \qquad \forall \ |(x, \xi)|\geq 1, 
\end{equation}

\noindent and that, since $w_0$ is definitely homogeneous of degree zero (see Definition \ref{def:pos_hom}), then
\begin{equation}\label{homo.lam}
\sup_{|(x, \xi)| \geq \lambda_\e \atop t \in \T} \frac{\Big|\pa_x^{\alpha}\pa_{\xi}^{\beta} w_0(t, x, \xi) \Big|}{(|x|^2+|\xi|^2)^{-\frac{|\alpha|+|\beta|}{2}}} = \sup_{|(x, \xi)| \geq 1 \atop t \in \T} \frac{\Big|\pa_x^{\alpha}\pa_{\xi}^{\beta} w_0(t, x, \xi) \Big|}{(|x|^2+|\xi|^2)^{-\frac{|\alpha|+|\beta|}{2}}} \,,
\end{equation}
hence
\begin{equation}
\begin{split}
   &\sum_{|\alpha|+|\beta| \leq 1} 
   \sup_{\substack{|(x, \xi)|\geq 1 \\t \in \T}} \frac{\Big|\pa_x^{\alpha}\pa_{\xi}^{\beta} w_0(t, x, \xi) \Big|}{(1+ |x|^2+|\xi|^2)^{-\frac{|\alpha|+|\beta|}{2}}}   \stackrel{\eqref{quasi.homo}}{\leq} 2 \sum_{|\alpha|+|\beta| \leq 1} 
   \sup_{\substack{|(x, \xi)|\geq 1 \\t \in \T}} \frac{\Big|\pa_x^{\alpha}\pa_{\xi}^{\beta} w_0(t, x, \xi) \Big|}{(|x|^2+|\xi|^2)^{-\frac{|\alpha|+|\beta|}{2}}} \\
  &\stackrel{\eqref{homo.lam}}{\leq}  2 \sum_{|\alpha|+|\beta| \leq 1} 
   \sup_{\substack{|(x, \xi)|\geq \lambda_{\epsilon} \\t \in \T}} \frac{\Big|\pa_x^{\alpha}\pa_{\xi}^{\beta} w_0(t, x, \xi) \Big|}{(1+|x|^2+|\xi|^2)^{-\frac{|\alpha|+|\beta|}{2}}}\stackrel{\eqref{5epsilon}}{\leq} 10 \epsilon, 
\end{split}
\end{equation}
proving \eqref{R.seminorm}. 
\end{proof}

\subsection{The set $\cV$ is dense}
The starting point to prove density of $\cV$ is to use Proposition \ref{simple_structure_thm} to characterize a subset $\cA \subset \cV$, that we will prove to be dense in $C^0(\T, \bS^0_{cl})$.\\
We recall the definitions of classical symbol \ref{class_sym}, resonant average of a symbol \ref{resonant_average} and that of $Z_{\mathbcal{e}_0}$ \eqref{projection}.\\
Throughout this section we will profusely use the results of Section \ref{sec:escape}, where the functions are positively homogeneous and not definitely homogeneous as in the definition of classical symbols. 
In particular we need to transform the definitely homogeneous symbol $\la v_0 \ra$ (averaged with respect to time), into a  positively homogeneous function; we do this  defining 
\begin{equation}\label{ext.sym}
    [v_0] := \boldsymbol{\pi}^* \la v_0 \ra = \la v_0 \ra \circ \boldsymbol{\pi}\, , \quad \boldsymbol{\pi} \mbox{ in } \eqref{projection} \ .
\end{equation}

\begin{remark}\label{rmk.out.ball}
Note that $[v_0] = \la v_0 \ra$ for $ |(x, \xi)| \geq 1 $.
\end{remark}
We shall use repeatedly the notations and definitions  of Section \ref{subsec:escape}.
The main result of this section is the following:
\begin{proposition}\label{WH_fol_generic}
With the notation in \eqref{ext.sym}, denote by 
    \begin{equation}\label{setA}
    \begin{aligned}
        \cA := & \{v \in C^{0}(\T, \bS^0_{cl}) \colon \ \exists \ v_0 \mbox{ 
           principal symbol of } v  \mbox{ such that   
 the vector field } \tilde \cX_{[ v_0] } \\
       & \mbox{ induces a  foliation  on } Z^{[v_0]}_{\mathbcal{e}_0} \mbox{ with weakly hyperbolic simple structure,}\\
       & \mbox{ for some regular value } \mathbcal{e}_0 \in \Ran([v_0]) \} \ .
         \end{aligned}
    \end{equation}
 Then  $\cA\subseteq \cV$, with $\cV$ in \eqref{setV}, and $\cA$ is dense in $C^0(\T, \bS^0_{cl})$ with respect to the distance $\td^{0}_{\T}$ (see \eqref{distance}).
\end{proposition}
In the rest of this section we prove Proposition \ref{WH_fol_generic}, that, together with Proposition \ref{V_open},  proves Theorem \ref{main_theorem} point $(i)$. To prove Proposition \ref{WH_fol_generic} we shall follow the following scheme:
\begin{enumerate}
    \item given any positively homogeneous function $h$ of degree zero and a regular value $\mathbcal{e}_0$ for $h$, the foliation $\mathscr{F}^h$ induced by $\tilde \cX_h$  on $Z_{\mathbcal{e}_0}^h$ is the characteristic foliation of $Z_{\mathbcal{e}_0}^h$ regarded as a submanifold of a contact manifold;
    \item given $h$ as in point 1, there exists $h_{\epsilon}$, positively homogeneous of degree zero and arbitrarily close to $h$ (when restricted to $\S^3$),  such that the foliation induced by $\tilde{\cX}_{h_{\epsilon}}$ is Morse-Smale; 
    \item given $h_{\epsilon}$ as in point 2,  the induced characteristic foliation $\mathscr{F}^{h_\e}$ has actually weakly hyperbolic simple structure;
     \item given $ v \in C^0(\T, \bS^0_{cl})$, let $h=[v_0]$, and $h_\e$ as in point 3, there exists a definitely homogeneous symbol $v_{\epsilon} \in C^0(\T, \bS^0_{cl})$ arbitrarily close to $ v_0$ and such that $[ v_{\epsilon} ]= h_{\epsilon}$. This concludes the proof of density of the set $\cA$.
     \item Prove that $\cA \subset \cV$.
\end{enumerate} 
We now prove each step of the scheme.

\subsubsection{The induced foliation is a  characteristic foliation}

To prove the first step we are going to identify $Z^h_{\mathbcal{e}_0}$ as a submanifold of $\S^3$, which we regard as a contact manifold with contact  structure induced by the symplectic structure of $(\R^4, \Omega)$, where $\Omega$ is the standard symplectic form, $\Omega:= \di x_1 \wedge \di \xi_1 +\di x_2 \wedge \di \xi_2$.

\begin{proposition}\label{Char_fol_Zw}
Let $h:\R^4 \setminus \{0\} \to \R$ be a positively homogeneous function of degree zero, fix a regular value $\mathbcal{e}_0$ for $h$ and let $\tilde{\cX}_h$ be the vector field on $Z^h_{\mathbcal{e}_0}$ (see \eqref{Zwo}) defined in \eqref{Xh}. 
Let $\tilde{\mu}$ be the volume form  on $Z^h_{\mathbcal{e}_0}$ defined in \eqref{tildemu}.
Then $\tilde{\cX}_h$ directs the characteristic foliation of $Z^h_{\mathbcal{e}_0}$ induced by the contact manifold $(\S^3, \kappa)$ with $\kappa=\ker(\alpha)$ and $\alpha$ defined by Corollary \ref{contact.sphere}. In particular, 
\begin{equation}\label{tildeX.char}
    i_{\tilde{\cX}_h} \tilde{\mu}= \left( 2 |\nabla h|\   \alpha\right){\vert_{Z^h_{\mathbcal{e}_0}}}\,.
\end{equation} 
\end{proposition}

\begin{proof}
We compute $i_{\tilde{\cX}_h} \tilde{\mu}$.
Using the definition \eqref{tilde_Xh} of $\tilde \cX_h$ and 
identity \eqref{tildemu.prop}, we have  $i_{\tilde{\cX}_h} \tilde{\mu}= i_{\rho \cX_h}\tilde{\mu}$.
Next, from the definition \eqref{tildemu}, 
$$
\tilde{\mu}= -i_{\frac{\nabla h}{|\nabla h|}} i_{\nabla h_0}\nu =  i_{\frac{\nabla h}{|\nabla h|}} i_{\nabla h_0} \frac{ \Omega \wedge \Omega }{2}
$$
with $\Omega$ the standard symplectic form in $\R^4$.
Recalling the expression of $\alpha$ in Corollary \ref{contact.sphere}, the equality \eqref{tildeX.char} reads
$$
 i_{\rho\cX_h} i_{\frac{\nabla h}{|\nabla h|}} i_{\frac{\nabla h_0}{2}}
\Omega \wedge \Omega = (2 |\nabla h| \, i_{\frac{\nabla h_0}{2} } \Omega )\vert_{Z_{\mathbcal{e}_0}}.
$$
We now verify it. Recalling that $\Omega \Big(\frac{\nabla h}{|\nabla h|}, \cX_h \Big) = |\nabla h|$, a direct computation yields 
\begin{equation}
        i_{\cX_h} i_{\frac{\nabla h}{|\nabla h|}} i_{\frac{\nabla h_0}{2}}
(\Omega \wedge \Omega) = 
 2|\nabla h| i_{\frac{\nabla h_0}{{2}}} \Omega - \Omega \Big( \frac{\nabla h}{|\nabla h|}, {\nabla h_0} \Big) i_{\cX_h}\Omega +  \Omega\Big( \cX_h, {\nabla h_0} \Big) i_{\frac{\nabla h}{|\nabla h|}}\Omega\,.
\end{equation}
Now note that the coefficient 
$\Omega\Big( \cX_h, {\nabla h_0}\Big) = 0$ due to \eqref{gradh.ho}, whereas the one form
$
i_{\cX_h}\Omega\vert_{Z^h_{\mathbcal{e}_0}}  = \di h\vert_{Z^h_{\mathbcal{e}_0}}$ vanishes on $TZ^h_{\mathbcal{e}_0}$ since $Z^h_{\mathbcal{e}_0}$ is inside the energy manifold $h^{-1}(\mathbcal{e}_0)$.
\end{proof}

\subsubsection{Density of functions $h$ such that $\tilde{\cX}_h$ has Morse-Smale foliation}

\noindent The main result of this section is the following: 

\begin{proposition}\label{close.MS}
For any $h: \R^4 \setminus \{0\} \to \R$ positively homogeneous of degree zero, for any regular value $\mathbcal{e}_0$ of $h$ and for any $\epsilon >0$, there exists $h_{\epsilon}$, positively homogeneous of degree zero,  such that $\tilde{\cX}_{h_{\epsilon}}$ directs a Morse-Smale characteristic foliation on $Z_{\mathbcal{e}_0}^{h_\e}$, and
$
\td_{C^{\infty}(\S^3)}(h, h_{\epsilon})< \epsilon.$
\end{proposition}

\noindent We need a preliminary result.

\begin{lemma}\label{diffeo.MS}
    Let $M$, $N$ be smooth compact manifolds of dimension two and let $f:M \to N$ be a smooth diffeomorphism. If $X$ is a Morse-Smale vector field on $N$ (see Definition \ref{Morse-Smale}), then also $f^*(X)$ is a Morse-Smale vector field on $M$.
\end{lemma}

\begin{proof}
Recalling the definition of the vector field $f^*(X)$ as 
    $
    f^*(X)(p)=\di f^{-1}(f(p)) [X(f(p))],
    $
   one starts with observing that the number of critical elements of $f^*(X)$ coincides with the number of critical elements of  $X$, hence it is finite.
  Furthermore, since
    $f^{-1} \circ \phi^t_X \circ f =\phi^t_{f^*(X)}$, one checks that 
 if every orbit of $X$ has one and only one critical element as $\alpha$ limit and one and only one as $\omega$ limit, then the same holds for every orbit of $f^*(X)$.
    Finally, since the Jacobian matrix of $X$ is similar to that of $f^*(X)$, the hyperbolicity of the critical points is preserved. If there are no saddle connections for the flow of $X$, the same holds for the flow of $f^*(X)$. 
\end{proof}

We say that a smooth family of maps $\psi_t \colon M \to M$, $t \in [0,1]$,  has size $\e>0$ if $\td_{C^\infty}(\psi_t, \uno) \leq \e$ for all $t \in [0,1]$.
In particular it is possible to write
\begin{equation}\label{small.isotopy}
    \psi_1 = \uno + \e \phi_1 ,  \ \ \ \ \td_{C^\infty}( \phi_1 ,0 ) \leq  1 \ .
\end{equation}
\begin{proof}[Proof of Proposition \ref{close.MS}]
    Fix $h: \R^4 \setminus\{0\} \to \R$ positively homogeneous function of degree zero, $\mathbcal{e}_0$ regular value for $h$ and $\epsilon>0$. Let $\tilde{\cX}_h$ be the vector field on $Z_{\mathbcal{e}_0}^h$ defined in \eqref{Xh}. 
    Using the genericity result of Theorem \ref{generic.MS} by Peixoto, for any $\epsilon_0>0$ there exists a Morse-Smale vector field $X_1$ on $Z_{\mathbcal{e}_0}^h$, such that, recalling the definition of distance between vector fields in \eqref{distance.vf}, 
    \begin{equation}\label{X1.MS}
         \td_{C^{\infty}(Z_{\mathbcal{e}_0}^h)}(\tilde{\cX}_h, X_1) < \epsilon_0.
    \end{equation}
    {Note that $X_1$ does not necessarily direct the} {characteristic} foliation on $Z^h_{\mathbcal{e}_0}$ (hence we cannot guarantee that it is a vector field induced by a positively homogeneous function). However, we shall prove that it is possible to slightly deform $Z^h_{\mathbcal{e}_0}$ to a new submanifold $Z_{\mathbcal{e}_0}^{h_{\epsilon}}$ whose characteristic foliation (which is directed by $\tilde X_{h_\e})$ is  actually also directed by the pullback of $X_1$, and so it is Morse-Smale.
    
    {Precisely,} we are going to show that for any $\epsilon>0$, provided $\epsilon_0$ is small enough, there exist a positively homogeneous function of order zero $h_{\epsilon}: \R^4 \setminus \{0\} \to \R$, a diffeomorphism $\psi: Z_{\mathbcal{e}_0}^{h_{\epsilon}} \rightarrow Z^h_{\mathbcal{e}_0}$ and a smooth {strictly} positive function $f$ on $Z_{\mathbcal{e}_0}^{h_{\epsilon}}$ such that
    \begin{equation}\label{to.prove}
    \psi^*(X_1)= f \tilde{\cX}_{h_{\epsilon}} \ \mbox{ as vector fields on }Z_{\mathbcal{e}_0}^{h_{\epsilon}}
    \quad \mbox{ and } \quad \td_{C^{\infty}(\S^3)}(h, h_{\epsilon}) <\epsilon.
    \end{equation}
    This concludes the proof since, from Lemma \ref{diffeo.MS}, the flow of $\psi^*(X_1)$ on $Z_{\mathbcal{e}_0}^{h_{\epsilon}}$ is still Morse-Smale, and $ f \tilde{\cX}_{h_{\epsilon}}$ directs the same foliation of $\tilde \cX_{h_\e}$. \\
    We now prove \eqref{to.prove}.
    {The key idea is to regard $Z^h_{\mathbcal{e}_0}$ as a  submanifold of the contact manifold  $\S^3$, and  use  the local Gray theorem \ref{Gray} to deform  it to a new submanifold whose characteristic foliation is Morse-Smale.}
    {So define  $\alpha_0:= 2 |\nabla h| \alpha$ and  $\alpha$ the Liouville 1-form as in Corollary \ref{contact.sphere}.
    Note that $\alpha$ and $\alpha_0$ define the same contact structure in a neighbourhood $\cU(Z^h_{\mathbcal{e}_0})\subset \S^3$ since they differ for function  never-vanishing on $\cU(Z^h_{\mathbcal{e}_0})$ (provided the neighbourhood is chosen sufficiently small, since $|\grad h|_{\vert Z^h_{\mathbcal{e}_0}}$ is non-vanishing).}  
    
    So define on $Z^h_{\mathbcal{e}_0}$ the 1-form given by $\beta:= i_{X_1 - \tilde{\cX}_h} \tilde \mu.$ By \eqref{X1.MS} on a sufficiently small neighborhood $\cU(Z^h_{\mathbcal{e}_0})$ it is possible to extend $\beta$ to a 1-form, which we still denote by $\beta$, such that $\|\beta\|_{C^\infty(\cU(Z^h_{\mathbcal{e}_0}))} \leq 2 \epsilon_0$. 
    

Thus for any $t \in [0,1]$ the 1-form
    \begin{equation}
        \alpha_t := \alpha_0 + t \beta
    \end{equation}
    is well defined on $\cU(Z^h_{\mathbcal{e}_0})$, and it is still a contact form, since
    $$
    \alpha_t \wedge \di \alpha_t = \alpha_0 \wedge \di \alpha_0 + t \alpha_0 \wedge \di \beta + t \beta \wedge \di \alpha_0 + t^2 \beta \wedge \di \beta
    $$
    is non-degenerate on $\cU(Z^h_{\mathbcal{e}_0})$, due to the fact that $\alpha_0$ is positively proportional 
    on $\cU(Z^h_{\mathbcal{e}_0})$ to the contact form $\alpha$ and to \eqref{X1.MS}.
    Moreover on $Z^{h}_{\mathbcal{e}_0}$, in view of \eqref{tildeX.char} and the very definition of $\beta$, 
    \begin{equation}\label{X1.rel}
        i_{X_1} \tilde \mu  = \alpha_1 \vert_{Z^h_{\mathbcal{e}_0}}\,.
    \end{equation}
    Since $\alpha_t$ is a family of contact forms in $\cU(Z^h_{\mathbcal{e}_0})$ with $\|\dot \alpha_t\|_{C^\infty(\cU(Z^h_{\mathbcal{e}_0}))}= \|\beta\|_{C^\infty(\cU(Z^h_{\mathbcal{e}_0}))}\leq 2 \varepsilon_0$, we can now use Gray's Theorem \ref{Gray}: 
    {for any $\varepsilon_1>0$ (which will be fixed later), if $\varepsilon_0$ is small enough,}
    there exists a family of smooth diffeomorphisms $\{ \psi_t\}_{t \in [0,1]}$ on $\cU(Z^h_{\mathbcal{e}_0})$ and a family of {strictly} positive  smooth functions $\lambda_t$ on $\cU(Z_{\mathbcal{e}_0}^h)$ such that
        $\psi^*_t(\alpha_t)= \lambda_t \alpha_0$ {and \eqref{small.isotopy} holds, with $\varepsilon$ replaced by $\varepsilon_1$.} 
  Define now the positively homogeneous function
    \begin{equation}\label{h.epsilon}
    h_\e := h \circ \psi_1 \circ \boldsymbol{\pi}\ ;
    \end{equation}
    clearly  $Z^{h_\e}_{\mathbcal{e}_0} = \psi_1^{-1}(Z^h_{\mathbcal{e}_0})$. 
    Now we show that the diffeomorphism in equation \eqref{to.prove} can be taken as  $\psi=\psi_1$. 
    We first evaluate $\psi^*_1(\alpha_1)$ on $Z^{h_\epsilon}_{\mathbcal{e}_0}$:
    $$
    \psi^*_1(\alpha_1) \stackrel{\eqref{X1.rel}}{=} \psi^*_1(i_{X_{1}}\tilde{\mu})=i_{\psi^*_1({X_{1}})}\psi^*_1(\tilde{\mu})\,.
    $$
    Combining this identity with  \eqref{iso.psi} and using the definition of $\alpha_0$, we have:
    \begin{equation}\label{psi.id1}
        i_{\psi_1^*(X_{1})}\psi^*_1(\tilde{\mu})= \lambda_1 \alpha_0 = 2 |\nabla h| \lambda_1  \alpha \quad \mbox{on} \quad Z^{h_\epsilon}_{\mathbcal{e}_0}\,.
    \end{equation} 
We define on $Z_{\mathbcal{e}_0}^{h_{\epsilon}}$ the volume form  $\tilde \mu_\epsilon := - i_{\frac{\nabla h_\epsilon}{|\nabla h_\epsilon|}} \nu_{\S^3},$ with $\nu_{\S^3}$ as in Lemma \ref{lemma.contrazioni}; by Proposition \ref{Char_fol_Zw} with $h$ replaced by $h_\epsilon$ and $\tilde \mu$ replaced by $\tilde \mu_\epsilon$, we have \begin{equation}\label{xhe.char}
    i_{\tilde \cX_{h_\epsilon}} \tilde \mu_\epsilon= 2 | \nabla h_\epsilon| \alpha \quad \mbox{on} \quad Z^{h_\epsilon}_{\mathbcal{e}_0}\,.
\end{equation}
Moreover, we observe that both $\tilde \mu_\epsilon$ and $\psi_1^*(\tilde \mu)$ are volume forms on $Z^{h_\epsilon}_{\mathbcal{e}_0},$ thus there exists a smooth function $g:Z_{\mathbcal{e}_0}^{h_{\epsilon}} \to \R^+$ such that $\psi^*_1(\tilde{\mu})= g \tilde{\mu}_{\epsilon}$. Then comparing \eqref{xhe.char} and \eqref{psi.id1} one gets that on $Z^{h_\epsilon}_{\mathbcal{e}_0}$
$$
i_{\psi_1^*(X_1)} \tilde \mu_{\epsilon} = g ^{-1} i_{\psi_1^*(X_1)} \psi_1^* (\tilde \mu) = 2 g^{-1} |\nabla h | \lambda_1 \alpha = g^{-1} \lambda_1  |\nabla h | |\nabla h_\e|^{-1} i_{\tilde \cX_ {h_{\epsilon}}} \tilde \mu_\epsilon\,,
$$
thus the two vector fields $\psi_1^*(X_1)$ and $\tilde \cX_{h_\epsilon}$ are proportional. 

To prove \eqref{to.prove}, it only remains to show that $\di_{C^\infty}(h, h_\epsilon) < \epsilon.$
 First choose $N \in \N$ large enough 
such that $\sum_{n=N+1}^{+\infty} \frac{1}{2^n} \leq \frac{\epsilon}{2}$, so that:
\begin{equation}\label{dist.cinf}
\td_{C^\infty}(h,h_{\epsilon})\leq \sum_{n=0}^{N}\frac{1}{2^n} \frac{\norm{h-h_{\epsilon}}_{C^n(\S^3)}}{1+\norm{h-h_{\epsilon}}_{C^n(\S^3)}} +\frac{\epsilon}{2}.
\end{equation}
Note that it is sufficient to prove that $\forall n = 0, \dots, N$ there exists $C_n>0$ such that
\begin{equation}\label{seminorms.small}
\norm{h-h_{\epsilon}}_{C^n(\S^3)} \leq C_n {\epsilon_1}\,,
\end{equation}
since then \eqref{dist.cinf} and \eqref{seminorms.small} give $\td_{C^\infty(\S^3)}(h, h_\e) \leq \sum_{n=0}^N C_n {\e_1} + \frac{\e}{2} < \e\,,$ up to choosing ${\e_1} < {\e}(2 \sum_{n=0}^N C_n)^{-1}$.
To deduce \eqref{seminorms.small}, with a direct computation one inductively checks that, $\forall n \in \N$ and $\alpha, \beta \in \N^2$ with $|\alpha| + |\beta| = n$, one has
\begin{equation}\label{inductive}
\partial_x^\alpha \partial_\xi^\beta (h-h_\e)(x, \xi) = f^{(1)}_{{\e_1}, \alpha, \beta}(x,\xi) - f^{(1)}_{{\e_1}, \alpha, \beta} \left((x,\xi) + {\e_1} \phi_1(x, \xi)\right) + {\e_1} f^{(2)}_{{\e_1}, \alpha, \beta}( x, \xi)\,,
\end{equation}
with $f^{(j)}_{{\e_1}, \alpha, \beta}$ smooth functions, whose $C^k(\S^3)$ norms are finite and have uniform upper bound in ${\e_1}$ for any $k$, and $\phi_1$ as in \eqref{small.isotopy}. Then \eqref{inductive} gives
$$
\sup_{\S^3} |\partial_x^\alpha \partial_\xi^\beta (h-h_\e)(x, \xi)| \leq {\e_1} \left(\sup_{\S^3} |\nabla f^{(1)}| \sup_{\S^3} |\phi_1| +  \sup_{\S^3} |f^{(2)}_{{\e_1}, \alpha, \beta}| \right) \leq {\e_1} C_{\alpha, \beta}
$$
for some positive $C_{\alpha, \beta}$, and taking the sup over all $\alpha, \beta$ with $|\alpha| + |\beta| \leq n$ \eqref{seminorms.small} follows, concluding the proof.
\end{proof}

\subsubsection{If $\mathscr{F}$ is Morse-Smale, it has weakly hyperbolic simple structure}
In this section we show that if the foliation $\mathscr{F}^h$ induced by $\tilde{\cX}_h$ on $Z_{\mathbcal{e}_0}^h$  is Morse-Smale, then it has weakly hyperbolic simple structure (see Definition \ref{weak_hyp}). 
To state precisely the result we need the following preliminary definition: 

\begin{definition}\label{def.ws}
    Let $\zeta_0$ be a hyperbolic critical point of the foliation $\mathscr{F}^h$. Then $\zeta_0$ is called weakly stable if there exists a strictly positive function $f: Z_{\mathbcal{e}_0}^h \to \R$ such that
    $\div_{\tilde{\mu}}(f\tilde{\cX}_h)(\zeta_0) <0\,.$
    Similarly, $\zeta_0$ is called weakly unstable if  there exists a positive function $f:Z_{\mathbcal{e}_0}^h \to \R$ such that 
    $\div_{\tilde{\mu}}(f \tilde{\cX}_h)(\zeta_0) >0\,.$
\end{definition} 
 The following result is a small variant of Theorem 6.2 of \cite{Colin_de_Verdi_re_2020}:
\begin{proposition}\label{MS.WH}
Let $\mathscr{F}^h$ be the foliation induced by $\tilde{\cX}_h$ on $Z_{\mathbcal{e}_0}^h$. 
    If $\mathscr{F}^h$ is Morse-Smale, then it admits a weakly hyperbolic simple structure (see Definitions \ref{def_simple_structure} and \ref{weak_hyp}). 
    In particular 
    $K^+$ is the finite union of the stable critical points and closed orbits, and of the unstable manifolds  of the weakly stable saddle points, whereas
         $K^-$ is the finite union of the unstable critical points and closed orbits, and of the stable manifolds  of the weakly unstable saddle points.
\end{proposition}

\begin{proof}
    The proof follows  the same construction of Theorem 6.2 in \cite{Colin_de_Verdi_re_2020}. The only difference is that, due to the difference in the homogeneity classes of the involved functions, we need to show in a different way that all critical points of $\tilde \cX_h$ are either weakly stable or weakly unstable, according to Definition \ref{def.ws}. In order to prove it, we proceed as follows: suppose by contradiction that there exists $\zeta_0 \in Z_{\mathbcal{e}_0}^h$ such that, for any smooth function $f$ defined in a neighborhood of $\zeta_0$ and such that $f(\zeta_0) > 0,$ 
$\div_{\tilde{\mu}}(f \tilde{\cX}_h)(\zeta_0) =0\,.$
Then, choosing $f$ such that $f(\zeta_0) = 1$ and using \eqref{ident}, we get 
$\{h, h_0\}(\zeta_0)=0\,.$
But then, since $\zeta_0$ is a critical point of $\tilde \cX_h$, from \eqref{Xh} we deduce $\cX_h(\zeta_0) = 0$, contradicting that $\mathbcal{e}_0$ is a regular value for $h$.
\end{proof}

\subsubsection{Density of symbols whose characteristic foliation is Morse-Smale}

{In this section we put together the results obtained so far, and we show that, given any symbol $v \in C^0(\T, \bS^0_{cl})$, it is possible to slightly perturb it into a new symbol $v_\e$ whose projected vector field directs a foliation with  weakly hyperbolic simple structure.}
Precisely we have  (recall the notation  \eqref{ext.sym} and \eqref{Zwo}):

\begin{proposition}\label{prop.density}
Let $v \in C^0(\T, \bS^0_{cl})$. Given any $\e >0$, there exists $v_{\epsilon} \in  C^0(\T, \bS^0_{cl})$ of the form $v_\e=v_{\e0} + v_{-\mu}$, with $v_{\e0}$ positively homogeneous and 
$v_{-\mu} \in C^0(\T, \bS^{-\mu})$,  $\mu>0$,  and $\mathbcal{e}_0$ regular value of $[v_{\e0}],$
such that
$\td^{0}_{\T}(v , v_\e) < \e$, 
 and the foliation induced by $\tilde{\cX}_{[v_{\epsilon0}]}$ on $Z_{\mathbcal{e}_0}^{[v_{\epsilon0}]}$  has weakly hyperbolic simple structure.
\end{proposition}

\begin{proof}
First we observe that, without loss of generality, we can assume that the 
{definitely homogeneous symbol}
$v$ admits the decomposition $v = v_0 + v_{-\mu}$, with
\begin{equation}\label{def.14}
v_0(x, \xi) = v_0(\lambda x, \lambda \xi) \quad \forall \lambda \geq 1\,, \ |(x, \xi)|\geq \frac 1 4\,, \quad v_{0} \in C^0(\T; \bS^0)\,, \quad v_{-\mu} \in C^0(\T; \bS^{-\mu})\, \ \mu >0.
\end{equation}
{Note that the homogeneity is here required to hold in $|(x, \xi)|\geq \frac 1 4$ and not only in 
$|(x, \xi)|\geq 1$ as in Definition \ref{def:pos_hom} of definitely homogeneous symbol.}
{However, we can always reduce to \eqref{def.14} by a different definition of the smoothing component.}
Precisely,  since $v \in  C^0(\T; \bS_{cl}^0)$, we have $v= \tilde v_0 + \tilde v_{-\mu}$, with $\tilde v_j \in C^0(\T; \bS^{j}_{cl})$, $j=0, -\mu$, and $\tilde v_0$ definitely homogeneous according to Definition \ref{def:pos_hom}. Then to obtain \eqref{def.14} it is sufficient to choose
\begin{equation}\label{new.dec}
v_{0} := (1-\eta) (\tilde v_0 \circ \boldsymbol{\pi})\,, \quad v_{-\mu} := \tilde{v_0} - (1-\eta)(\tilde v_0 \circ \boldsymbol{\pi})  + \tilde v_{-\mu}\,,
\end{equation}
with $\eta \in C^\infty(\R^4; \R)$ a non-negative radial cut-off function satisfying
$$
\eta(x,\xi) =
\begin{cases}
   1 & |(x, \xi)| \leq \frac 1 8\,,\\
   0 & |(x, \xi)| \geq \frac 1 4\,.
\end{cases}
$$
Note that $v_0(t, \cdot)$ is a smooth function, since the function $\tilde v_0 \circ \boldsymbol{\pi}$ has a singularity only at $(x, \xi) = 0$, where $1-\eta$ identically vanishes; furthermore $v_0$ coincides with $\tilde v_0  \in C^0(\T; \bS^0_{cl})$ outside the ball $B_1(0):=\{(x, \xi) \in \R^4 \ |\ |(x, \xi)| \leq 1 \}$, thus one also has $v_0 \in C^0(\T; \bS^0_{cl})$, and $v_0$ satisfies the homogeneity property in \eqref{def.14} by construction. Similarly, we observe that the function $\tilde{v_0} - (1-\eta) (\tilde v_0 \circ \boldsymbol{\pi}) $ has no singularities and identically vanishes outside $B_1(0)$, thus we deduce that $\tilde{v_0} - (1-\eta) (\tilde v_0 \circ \boldsymbol{\pi}) \in  C^0(\T; \bS^{-\infty}_{cl})$ and $v_{-\mu} \in C^0(\T; \bS^{-\mu}_{cl})$. This proves \eqref{def.14}.

\noindent
Assume now that $v_0$ is as in \eqref{def.14}. We fix a smooth, non-negative, radial cut-off function $\chi$ fulfilling
\begin{equation}\label{cut.chi}
\chi \in C_c^{\infty}(\R^4, [0, 1]), \qquad \chi(x, \xi)= 
\begin{cases}
 1 & |(x, \xi)| \leq \frac14\\
 0 & |(x, \xi)| \geq 1
\end{cases} \ ,
\end{equation}
and we are going to construct $v_\e$ as follows:
\begin{equation}
    v_\e := v_{\e0} + v_{-\mu} \,,
    \end{equation}
    where
\begin{equation}
    v_{\e0}(t, \cdot):= v_0 + \big( (1-\chi) ( h_\e  - [ v_0]) \big)\circ \phi^{-t}_{h_0}\,,
\end{equation}
and $h_\e \in C^\infty(\R^4 \backslash \{0\}; \R)$ is a suitable smooth, positively homogeneous function of degree $0$, that we now are going to exhibit. We shall require that $h_\e$ also satisfies the following properties:
\begin{itemize}
    \item[$(i)$]  $\td_{C^\infty(\S^3)}([v_0], h_\e) < \epsilon$;
    \item[$(ii)$] there exists $\mathbcal{e}_0 \in \Ran(h_\e)$, regular value for $h_\e$, such that the foliation induced by  $\tilde \cX_{h_\e}$ on $Z^{h_\e}_{\mathbcal{e}_0}$ has weakly hyperbolic simple structure.
\end{itemize}
If $h_\e$ is as above, then observing that, on the support of $1-\chi$, both $ h_\e$ and $[v_0]$ are smooth and positively homogeneous functions of degree 0, and using Lemma \ref{composition_flow}, it is clear that $v_{\e0}$ is a smooth, definitely homogeneous function of degree 0. Thus we have $v_{\e} = v_{\e 0} + v_{-\mu} \in C^0(\T; \bS^0_{cl})$ and, by $(i)$,
$$
\td^{0}_{\T}(v, v_\e) \leq  \td^{0}([v_0] (1-\chi), h_\e (1-\chi) ) \leq \e\,.
$$
Furthermore, by construction we have that $\la v_{\e0} \ra = \la v_0\ra + (1-\chi)(h_\e - [v_0])$ (see Definition \ref{resonant_average}), and using that $\chi \equiv 0 $ on $|(x,\xi)| \geq 1$, we also have 
$[v_{\e0}] = h_\e\,,$
thus by $(ii)$ there exists $\mathbcal{e}_0 \in \Ran{[v_{\e0}]}$ such that the foliation induced by $\tilde \cX_{[v_{\e0}]}$ on $Z^{[v_{\e0}]}_{\mathbcal{e}_0}$ has weakly hyperbolic simple structure. 

Then it remains to exhibit $h_\e$ with the above properties in order to prove the proposition.
To this aim, consider the positively homogeneous function $[v_0]$. Two cases are given: either 
 $[v_0]$ is identically constant, or not.
In the first case, 
we define
\begin{equation}
    h_\e:=  [v_0] + \e h_\star\,, 
\end{equation}
with $h_\star$ as in Lemma \ref{lem.ms.ex} and we take $\mathbcal{e}_0 = [v_0] + \frac{\e}{2}$, and items $(i)$ and $(ii)$ follow from Lemma \ref{lem.ms.ex} and the arbitrariness of $\e$.\\ 
In the second case, let $\mathbcal{e}_0$ be a  regular value of $[v_0]$; the existence of a smooth, definitely homogeneous, function $h_\e$ such that item $(i)$ holds and the foliation induced by $\tilde \cX_{h_\e}$ on $Z^{h_\e}_{\mathbcal{e}_0}$ is Morse-Smale is guaranteed by Proposition \ref{close.MS}. Finally, item $(ii)$ holds by Proposition \ref{MS.WH}.
\end{proof}

\subsubsection{Proof of Proposition \ref{WH_fol_generic}}
First we show that $\cA \subseteq \cV$. 
Take $v \in \cA$;
by the very  definition \eqref{setA}, 
  the vector field $\tilde \cX_{[ v_0]} $ induces a foliation on $Z_{\mathbcal{e}_0}^{[v_0]}$ with a weakly hyperbolic simple structure, with $\mathbcal{e}_0$ a regular value of $[v_0]$. 
  By Proposition \ref{simple_structure_thm} and Corollary \ref{interval.escape}, there exist a positively homogeneous   function $\tilde a$ of degree two and an interval $I \ni \mathbcal{e}_0$ and contained in $\Ran([v_0])$, such that $\tilde a$ is an escape function for $[v_0]$ at any $\omega \in I $ and for any  $(x, \xi) \in [v_0]^{-1}(I) \cap \{ |(x, \xi)| \geq 1\}$.

Since $\tilde a$ is only positively homogeneous of degree $2$, we regularize it at 0 with the cut-off in \eqref{cut.chi}, i.e. we put 
$a:= (1-\chi)\tilde a$.
Then 
$$
\{ \la v_0 \ra, a \}  = \{ [v_0], \tilde a \} \geq \delta, \qquad  \forall \ (x, \xi) \in \la v_0 \ra^{-1}(I) \cap \{|(x, \xi)| \geq 1\}.
$$
Since, by Remark \ref{rmk.out.ball}, $\la v_0 \ra^{-1}(I) \cap \{|(x, \xi)| \geq 1\} = [v_0]^{-1}(I) \cap \{|(x, \xi)| \geq 1\}$, we have thus showed that $v \in \cV$.

The density of  $\cal{A}$ in  $C^0(\T, \bS^0_{cl})$ with respect to the distance $\td^{0}_{\T}$ is Proposition \ref{prop.density}.
\qed

\section{Growth of Sobolev norms for potentials in $\cV$}
In this last section we prove Theorem \ref{main_theorem} point $(ii)$, by applying the abstract  Theorem \ref{thm:ab0} with $K_0:= H_0$ and $V(t) =\Opw(v(t))$, $v \in \cV$. 
We verify that  assumptions \bI, \b{II}, \b{III} are met. 

\noindent{\bf Verification of Assumption I:} it holds trivially by the symbolic calculus of pseudodifferential operators (see Proposition \ref{symbolic.calculus}). \\
\noindent{\bf Verification of Assumption II:} $(i)$ holds since  
$\sigma(H_0) = \N+1$, whereas $(ii)$ by the exact Egorov theorem (see Proposition \ref{egorov}). \\
\noindent{\bf Verification of Assumption III:}
First note that, applying again Theorem  \ref{egorov}, if $K_0 = H_0$, then the two operators $\lala V \rara$ and $\la V \ra$ respectively defined in \eqref{Vaverage} and \eqref{average_operator} coincide.\\
$(i)$ Decompose as usual $v = v_0 + v_{-\mu}$, $\mu >0$, then by \eqref{Vaverage}, \eqref{average} and the exact Egorov theorem \ref{egorov}, one has
$$\la V \ra =\la V_0 \ra + \la V_{-\mu} \ra, \quad 
\la V_0 \ra :=\Opw( \la v_0 \ra) \in \cS^0\,,
$$
and $\la V_{-\mu} \ra=\Opw(\la v_{-\mu}\ra) \in \cS^{-\mu}$ a compact operator. In particular 
$
\sigma_{ess}(\la V \ra) = \sigma_{ess}(\la V_0 \ra ) 
$
and the last  spectrum is known to be
\begin{equation}
     \sigma_{ess}(\Opw(\la v_0\ra ))= [\min \la v_0 \ra\vert_{ \S^3}, \max \la v_0\ra \vert_{\S^3}]\ ,
\end{equation}
see e.g \cite[Theorem 2.7]{Mas23} or \cite[Theorem 2.1]{Colin_de_Verdi_re_2020}.
Since for  symbols in $\cV$ the function $\la v_0 \ra$  is not constant {on $\S^3$} (otherwise $\{\la v_0 \ra, a\} \equiv 0$ {in $|(x, \xi)|\geq 1$}),  Assumption \b{III} $(i)$  holds taking 
an interval
$I_0  \subset  \Ran(\la v_0 \ra\vert_{\S^3})$. \\
\noindent
$(ii)$ Put $A:= \Opw(a)$, where $a \in \bS^1_{cl}$ and $I \subset I_0$ are the escape function and the interval in \eqref{setV}.
Let $g \in C_c^{\infty}(\R, \R_{\geq 0})$ be a smooth and compactly supported function such that $\textrm{supp}(g) \subset I\,.$

We verify \eqref{mourre.est} for some $\theta >0 $ and $K$ compact and selfadjoint operator in $L^2(\R^2)$, and eventually shrinking $I$.
It is sufficient to show  \eqref{mourre.est}  for $V_0$, i.e. 
\begin{equation}\label{Mourre.V0}
    g_I(\la V_0 \ra)\, \im\,  [\la V_0 \ra, A]g_I(\la V_0 \ra) \geq \theta g_I^2(\la V_0 \ra)+K \ , 
\end{equation}
since then 
by  Lemma \ref{lemma_compactness_cutoff} and the fact that  $[\la V - V_0 \ra, A]$ is a compact operator
$$
 g_I(\la V \ra)\, \im [\la V \ra, A] \, g_I(\la V \ra) =  g_I(\la V_0 \ra)\im [\la V_0 \ra, A]g_I(\la V_0 \ra)+K,
 $$
with $K$ compact, proving \eqref{mourre.est}.

To show \eqref{Mourre.V0} we  use first  symbolic   and functional calculus (see Section \ref{sec:pseudo})   to write
\begin{equation}\label{mourre.symbols}
\begin{split}
    g_I(\la V_0 \ra)\, \im [\la V_0 \ra, A]\, g_I(\la V_0 \ra)
    &=  \Opw(g_I^2(\la v_0 \ra)\{\la v_0 \ra, a\})+K\\
\end{split}
\end{equation}
with $K$ a pseudodifferential operator of strictly negative order, hence compact. 

 Next we remark that, from the definition of the set $\cV$ in  \eqref{setV} and the support properties of $g_I$: 
$$
g_I^2(\la v_0 \ra)( \{ \la v_0 \ra, a \} -\delta ) \geq 0,  \qquad \forall \ |(x, \xi)| \geq 1,
$$
so the  strong Gårding inequality \ref{strong.garding} yields
\begin{equation}\label{v0.garding}
    \la \Opw(g_I^2(\la v_0 \ra) ( \{ \la v_0 \ra, a \} -\delta ))u, u \ra \geq -C\norm{u}_{-\frac12}^2.
\end{equation}
 Hence
\begin{equation}\label{v0.garding.2}
    \begin{split}
        \la \Opw(g_I^2(\la v_0 \ra) \{ \la v_0 \ra, a \}) u, u \ra 
        & \stackrel{\eqref{v0.garding}}{\geq} \delta \la  \Opw(g_I^2(\la v_0 \ra) u, u \ra -C\norm{u}_{-\frac12}^2 \\
        & \stackrel{\eqref{Sob.Spaces}}{=} \delta \la [ \Opw(g_I^2(\la v_0 \ra)) - C H_0^{-1}]u, u \ra,
    \end{split}
\end{equation}
 and since from \eqref{comp.symbols}
 \begin{equation}\label{v0.garding.3}
      \Opw(g_I^2(\la v_0 \ra)) = \Opw(g_I(\la v_0 \ra ))^2 + K_1, \qquad \mbox{ with } K_1 \in \cS^{-1},
 \end{equation}
we obtain 
\begin{equation}\label{v0.garding.4}
    \Opw(g_I^2(\la v_0 \ra) \{ \la v_0 \ra, a \}) \stackrel{\eqref{v0.garding.2}, \eqref{v0.garding.3}}{\geq} \delta \Opw(g_I(\la v_0 \ra ))^2 + K = \delta g_I^2(\la V_0 \ra) + K,
\end{equation}
where again the compact contribution $K$ changes at each equality. 

Plugging \eqref{v0.garding.4} in \eqref{mourre.symbols} we obtain 
\eqref{Mourre.V0} with $\theta=\delta$, concluding the proof.

\appendix

\section{Proof of Lemma \ref{lem.ms.ex}}\label{app:example}
\begin{figure}[ht]
    \centering
\includegraphics[ width=0.6\linewidth]{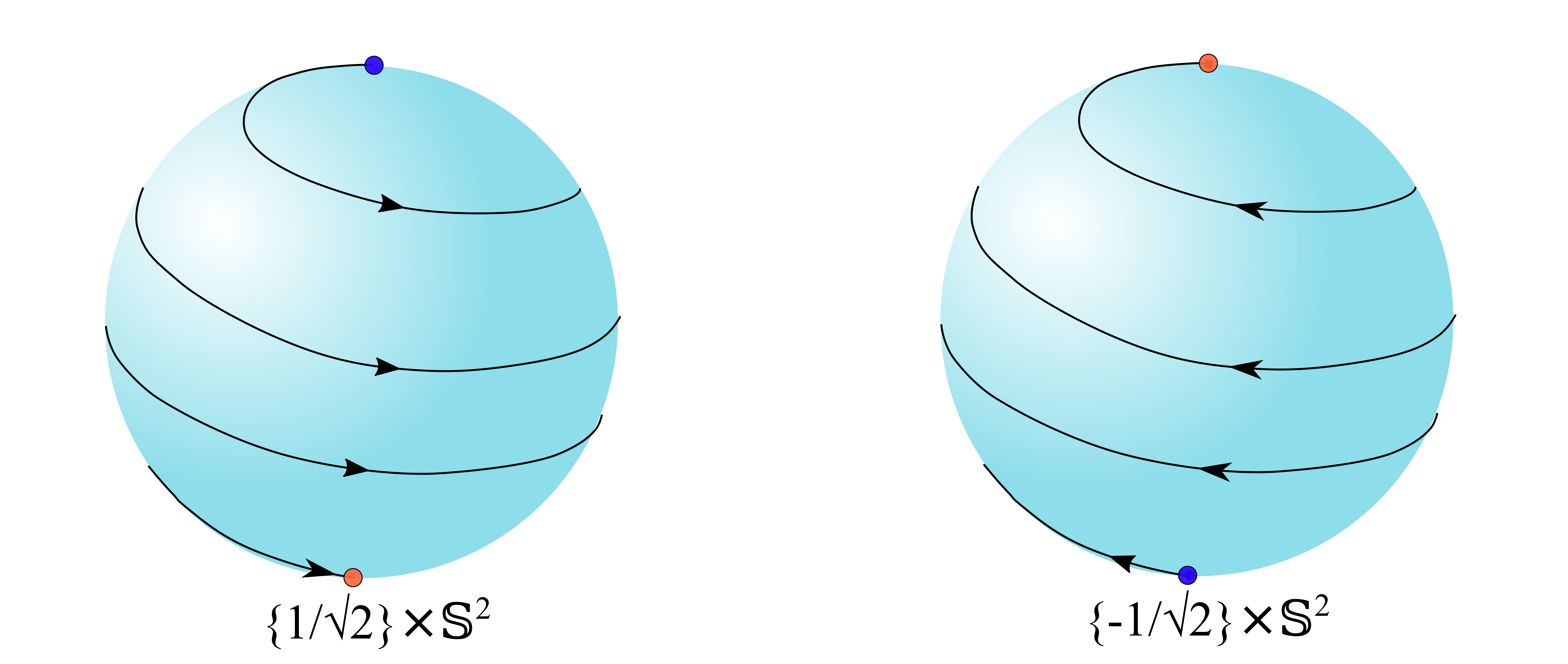}
    \caption{The flow of the vector field $\tilde X_{h_*}$. The repellor $K^-$ consists in the two blue dots, whereas the attractor  $K^+$ is formed by the orange dots.}
    \label{fig:enter-label}
\end{figure}

First note that the energy level $\Sigma_{\frac12}^{h_\star}=\{(x, \xi) \in \R^4 \colon x_1^2 = \xi_1^2 + x_2^2 + \xi_2^2  \}$. Then, by \eqref{Zwo}, 
$$
Z^{h_\star}_{\frac12} = \{ x_1 = \pm \frac{1}{\sqrt 2} \} \cup \{ \xi_1^2 + x_2^2 + \xi_2^2= \frac12 \}
$$
which is topologically $\S^2 \cup \S^2$. 
We now compute the induced vector field $\tilde \cX_{h_\star}$ on one of the spheres, say the one corresponding to $x_1 = \frac{1}{\sqrt{2}}$.
Such sphere is described, in the coordinates $\zeta$ of \eqref{projection}, by
\begin{equation}\label{app.res}
    \zeta_2^2 + \zeta_3^2  + \zeta_4^2 =  \frac12 \ ,
\end{equation}
where we have identified $\left(\frac{x_1}{\rho}, \frac{x_2}{\rho}, \frac{\xi_1}{\rho}, \frac{\xi_2}{\rho}\right) \equiv (\zeta_1, \zeta_2, \zeta_3, \zeta_4)$. We now compute $\tilde \cX_{h_\star}(\zeta)$ via its  definition \eqref{tilde_Xh}, getting 
\begin{equation}
    \tilde \cX_{h_\star}(\zeta)=
     2 \frac{x_1}{\rho^3} \begin{pmatrix}
        0 \\
       \xi_1 x_2 - x_1 \xi_2 \\
       \xi_1^2  + x_1^2 - 2h_0(x, \xi)  \\
        x_1 x_2 + \xi_1 \xi_2 \\
    \end{pmatrix} 
    =
    \sqrt{2} \begin{pmatrix}
        0 \\
       \zeta_3 \zeta_2  - \frac{1}{\sqrt{2}} \zeta_4 \\
       \zeta_3^2 - \frac12 \\
        \zeta_3 \zeta_4  + \frac{1}{\sqrt{2}} \zeta_2 
        \end{pmatrix}  \ .
    \end{equation}
By \eqref{app.res}, each variable $\zeta_j$ is defined in the interval $[-\frac{1}{\sqrt{2}}, \frac{1}{\sqrt{2}}]$. Consider first the equation for $\zeta_3$, given by $\dot \zeta_3  = \sqrt{2}(\zeta_3^2 -\frac12)$. 
The dynamics is clear: the fixed point $-\frac{1}{\sqrt{2}}$ is a global attractor and the fixed point $\frac{1}{\sqrt{2}}$ a global repellor. 
Then, in view of \eqref{app.res},  $\zeta^\pm:= (\frac{1}{\sqrt{2}}, 0, \mp \frac{1}{\sqrt{2}}, 0)$ are fixed points for the dynamics of $\tilde \cX_{h_\star}$. 
Moreover, since  $\zeta_3(t) \to \mp  \frac{1}{\sqrt{2}} $ for $t \to \pm \infty$ whenever $\zeta_3(0) \neq  \pm \frac{1}{\sqrt{2}}$, still from \eqref{app.res}
we deduce that
$(\frac{1}{\sqrt{2}}, \zeta_2(t), \zeta_3(t), \zeta_4(t))  \to \zeta^\pm$ when $t \to \pm \infty$.
This,  putting   $K^\pm := \{ \zeta^\pm \}$, shows that the foliation has a simple structure (Definition \ref{def_simple_structure}).
Note that the flow of $\tilde \cX_{h_\star}$ is the simplest Morse-Smale flow on $\S^2$! See Figure \ref{fig:enter-label}.

To show that the foliation is also weakly hyperbolic, it is easier to use formula
\eqref{ident}, that reads
$$
\div_{\tilde \mu}(\tilde \cX_{{h_\star}})  = -2 \{{h_\star}, h_0\}  =
4 x_1 \xi_1 {\vert_{Z_{\frac12}^{h_\star}} } \ ,
$$
with $\tilde \mu$ as in \eqref{tildemu}. Then, at the point $\zeta^+$, corresponding to $K^+$, the right hand side is strictly negative, hence by continuity the l.h.s is strictly negative also in a neighbourhood of $K^+$. 
Analogously, in a neighbourhood of $K^-$, the l.h.s. is strictly positive, proving the weak hyperbolicity (Definition \ref{def:wh}).

As a final comment, an escape function for $h_\star$ can be directly computed:  put $a(x, \xi):= -x_1 \xi_1$ to get 
$$
\{ {h_\star}, a\}  = 
\frac{\xi_1^2 x_1^2}{2 h_0^2} + \frac{x_1^2}{h_0} \big( 1- \frac{x_1^2}{2h_0}\big) \geq 2{h_\star} (1-{h_\star})
$$
which is strictly positive on any energy level $\Sigma_{\mathbcal{e}_0}^{h_\star}$ with $\mathbcal{e}_0 \in (0,1)$.

\section{Proof of Proposition \ref{global_flow}}\label{app:global_flow}

Before starting the proof, we recall the following topological properties of attractors \cite[Lemma 1.6]{Hurley}.

\begin{lemma}\label{property_attractor}
Let $K^+$ be an attractor and $V$ a  neighbourhood of $K^+$. Then there exists an open neighbourhood $U$ of $K^+$ such that
$K^+ \subset U \subset V$, 
$\bigcap_{t \geq 0}\phi^t_X(\overline{U})=K^+$ and 
$\phi^t_X(U) \subset U$ for all positive times.
\end{lemma}

We divide the proof of Proposition \ref{global_flow} in several steps. A key point will be to observe that, by Lemma \ref{ham_field_in_rho.zeta}, if $z(t)$ solves the Cauchy problem $\dot z = X_{h}(z)$, $z(t_0)=z_0$ for some $t_0 \in \R$, then its  components $\rho(t),\ \zeta(t)$ solve the Cauchy problem
\begin{equation}\label{dot.zeta}
\begin{cases}
  \dot \zeta = \rho^{-2}\tilde X_h(\zeta)\\
  \dot \rho = \rho^{-1}\{h, h_0\} 
\end{cases}\,, \quad
\begin{cases}
\zeta(t_0) = \zeta_0 := \boldsymbol{\pi}(z_0) \\
\rho(t_0) = \rho_0 := \|z_0\|
\end{cases} \,.
\end{equation}
For this reason, we will often reparametrize time as
\begin{equation}\label{rep.time}
    \tau_{t_0}(t):= \int_{t_0}^t \rho^{-2}(s)\ \di s
\end{equation}
and exploit the fact that, for any time interval $I \ni t_0$ such that $\rho(t) >0$ $\forall t \in I,$ one has
\begin{equation}\label{w.tau.t}
    \zeta(t) = w(\tau_{t_0}(t)) \quad\forall\ t \in I\,,
\end{equation}
where $w$ solves the Cauchy problem
\begin{equation}\label{dot.w}
   \frac{\di w}{\di \tau}(\tau) = \tilde{\cX}_h(w(\tau))\,, \quad w(0):= \zeta(t_0)\,.
\end{equation}
Note that the flow of \eqref{dot.w} is globally defined in time, since $Z_{\mathbcal{e}_0}$ is a compact manifold. Moreover, let $U^+$ be the open set given by the definition of weak hyperbolicity \ref{weak_hyp}. Since $K^+$ is an attractor, 
  Lemma \ref{property_attractor} guarantees the existence of an open neighbourhood of $K^+$, $\tilde U^+ \subset U^+$ such that, if $w(0) \in \tilde U^+$,
\begin{equation}\label{lim_wU}
    w(t) \in \tilde U^+ 
    \ \ \forall \ t\geq 0 \qquad \mbox{ and } \qquad w(t) \to K^+ \mbox{ as } t \to +\infty \ .
\end{equation}
Up to shrinking $U^+$ we are going to suppose $\tilde U^+ = U^+$.


\noindent{\bf Step 1:} \emph{Item $(i)$ holds, namely there exist open neighborhoods $U^\pm \subset Z_{\mathbcal{e}_0}$ of $K^{\pm}$ such that for any $z_0 \in \boldsymbol{\pi}^{-1}(U^\pm)$, the flow $\Phi^t_{\cX_h}(z_0)$ is globally defined for $\pm t \geq 0$ and 
 $\Phi^t_{\cX_h}(\boldsymbol{\pi}^{-1}(U^\pm)) \subseteq \boldsymbol{\pi}^{-1}(U^\pm)$ for $\pm t \geq 0$.}
\begin{proof}[Proof of Step 1]
 We prove the statement with ``$+$''.  The idea is to exploit that the projected flow $\zeta(t)$ stays in $U^+$ for all times; the difficulty is that the radial component $\rho(t)$ can either blow up or go to $0$ in finite time.
 While the first possibility is ruled out by homogeneity of the vector field, the second one would contradict the existence of the local escape function $k^+$.

Let  $z_0 = (\rho_0, \zeta_0) \in \boldsymbol{\pi}^{-1}(U^+)$ (recall  \eqref{projection}) with $\rho_0 \neq 0$, and 
denote
 $\Phi^t_{\cX_h}(\rho_0, \zeta_0)=(\rho(t), \zeta(t))$. 
We need to show that for all positive times $\rho(t)\in (0,\ +\infty)$  and $\zeta(t) \in U^+$. \\
Denote by $T_{\max} \in (0; +\infty]$ the maximal time such that $\rho(t) \in (0,\ +\infty)$. 
First remark that $\rho(t)$ can't blow up in finite time, since
$\rho=\sqrt{2h_0}$ and one has
\begin{equation}\label{hzerodot}
    |\dot{h}_0|=|\{h, h_0 \}| \leq \max_{Z_{\mathbcal{e}_0}} |\{h, h_0\}| =:  \nu
\end{equation}
that leads to the upper bound
\begin{equation}\label{eq.nu}
|h_0(t) - h_0(0)| \leq \nu  t   \ , \quad \forall t \geq 0 \ . 
\end{equation}
Hence the only possibility to have $T_{\max}<+\infty$ is that $\rho(t) \to 0$ as $t \to T_{\max}^-$.
This is not possible as long as 
$\zeta(t)$ belongs to $U^+$, as we now show.
Denote
by $t_*$ the escape time from $\boldsymbol{\pi}^{-1}(U^+)$, namely
$$
t_* := \inf\{ t \geq 0 \colon \zeta(t) \notin U^+ \}  \ . 
$$
As $\zeta(0) \in U^+$, $t_* >0$.
Now, by \eqref{hh0}, the escape function $k^+(t)$ is strictly increasing as long as $\zeta(t)$ belongs to $U^+$ and, by \eqref{def:kpm}, it is equivalent to $h_0$:
\begin{equation}\label{k+,h_0}
	c_2 h_0(t) \geq k^+(t) \geq k^+(0) \geq c_1 h_0(0) >0  \ , \quad \forall \ t \in [0, t_*] \ .
\end{equation}
This gives also that $T_{\max} > t_*$.

We are now going to prove that $t_* = + \infty$, implying both that $\Phi_{\cX_h}^t(\boldsymbol{\pi}^{-1}(U^+))$ is defined for all $t \geq 0$, and that $\boldsymbol{\pi}^{-1}(U^+)$ is invariant along $\Phi_{\cX_h}^t$ for $t \geq 0$. 

Assume by contradiction that $t_* < + \infty$. Then, recalling 
$t_* < T_{\rm max}$, take  $\e>0$ so small that  
$ t_* + \epsilon < T_{\rm max}$, so that we have 
\begin{equation}\label{309:1800}
\rho(t) \geq \bar \rho > 0  \ , \quad \forall t \in 
[t_*-\epsilon ,\ t_* + \epsilon ] \ .
\end{equation}
Then the  time reparametrization $\tau_{t_* - \epsilon}(t)$ in \eqref{rep.time} is well defined in the time interval of \eqref{309:1800}.
On this time interval, the identity  \eqref{w.tau.t} holds with $t_0 := t_*-\epsilon$ and $w(0) = \zeta(t_*-\epsilon) \in U^+$.
Since $w(\tau)$ is globally defined in time and $w(\tau) \in U^+$ for any $\tau \geq 0$, we get
$\zeta(t_* + \epsilon) = w\big(\tau_{t_*-\epsilon}(t_*+\epsilon) \big) \in U^+$ contradicting the maximality of $t_*$.
 \end{proof}
\noindent{\bf Step 2:} \emph{The set $\Gamma^\pm$ is invariant for all $\pm t \geq 0$. 
}
\begin{proof}[Proof of Step 2]
We prove the statement with ``$+$''.     
Denote again by 
 $\Phi^t_{\cX_h}(\rho_0, \zeta_0)=(\rho(t), \zeta(t))$ the flow of $\cX_h$ with initial datum $z_0=(\rho_0, \zeta_0)$, written in the coordinates  \eqref{projection}. We point out that $\rho(t)$ and $\zeta(t)$ are well defined for all $t\geq 0$ from Step 1 and 
\begin{equation}\label{309:1805}
\rho(s) \geq \bar \rho_t > 0  \ , \quad \forall s \in [0,\ t] \ .
\end{equation}
In particular the function $\tau_{0}(t)$ in \eqref{rep.time} is well defined for all $t\geq 0$.

We must show that whenever $\rho_0 >0$ and $\zeta_0 \in K^+$, 
$   \zeta(t) \in K^+$ for any $t >0$.
We have that
$\zeta(t) = w(\tau_0(t))$ and since   $  w(0) = \zeta(0) \in K^+\,$
and  
 $K^+$ is invariant for all $t>0$ under the flow of $\tilde{\cX_h}$, $\zeta(t) \in K^+$ for all $t\geq 0$, which  proves Step 2.
     \end{proof}
\noindent{\bf Step 3:}
{\em If $z_0 \in \Gamma^\pm$, there exists a finite time $\pm T^\pm(z_0)<0$ such that
	\begin{gather}
	\label{G.inv.neg}
	\Phi^t_{\cX_h}(z_0) \in \Gamma^{\pm} \quad \forall t \quad \text{s.t.} \quad  \pm T^{\pm}(z_0) < \pm t\,,\\
	\label{G.finite.death}
	\Phi^t_{\cX_h}(z_0) \to 0 \quad  \text{as} \quad  t \to T^\pm (z_0)^\pm\,.
	\end{gather}
	In particular, the maximal time of existence of $\Phi^t_{\cX_h}(z_0)$ for $z_0 \in \Gamma^+$ is $(T^+(z_0), +\infty)$, and it is $(-\infty, T^-(z_0))$ for $z_0 \in \Gamma^-$.}
\begin{proof}[Proof of Step 3] We prove the statement with ``$+$''. Let $T^+(z_0)< 0$ be the smallest time such that $\Phi^t_{\cX_h}(z_0)$ exists.
The fact that the invariance property \eqref{G.inv.neg} holds for all $t> T^+(z_0)$ can be shown exactly as in Step 2. 

We now prove that \eqref{G.finite.death} holds for a finite time $T^+(z_0)$. Pick $z_0 \in \Gamma^+$, and consider the escape function $k^+$ defined in \eqref{def:kpm}. The evolution of $k^+$ along the flow is given by the equation 
$
\dot{k}^+(t)=\{h, k^+\}$ for all
$t \in (T^+(z_0),\ \infty)$. 
From \eqref{hh0}, in $\Gamma^+$ we have a strictly positive bound from below on $\{h, k^+\}$, and since $\Gamma^+$ is invariant for the flow $\Phi^t_{\cX_h}$ as proved in  \eqref{G.inv.neg}, we have: 
$$
\dot{k}^+(t) \geq \delta, \qquad \forall t \in (T^+(z_0),\ \infty). 
$$
Integrating both sides between $t<0$ and $0$ we have:
\begin{equation}\label{cresce.poco}
\int_t^0 \dot k^+(s) \di s \geq -\delta t \implies k^+(t) \leq \delta t +k^+(0)\,.
\end{equation}
Then, combining \eqref{k+,h_0}
with \eqref{cresce.poco}, we deduce the existence of positive constants $c_1,\ c_2>0$ such that
\begin{equation}\label{bound.h0}
h_0(t) \leq c_1 \delta t + c_2 h_0(0) \,.
\end{equation}
The right hand side above tends to $0$ as $ t \to -\frac{c_2 h_0(0)}{c_1 \delta}<0$ (recall $h_0(0) >0$).
Thus  there exists a negative time $T^+(z_0) \in \big[-\frac{c_2 h_0(0)}{c_1 \delta}, 0\big)$ such that 
$\rho(t) = \sqrt{2 h_0(t)} \to 0$  as $ t \to T^+(z_0)\,,$
proving \eqref{G.finite.death}.
\end{proof}
\noindent Steps 2 and 3 prove Items $(ii)$ and $(iii)$ of Proposition \ref{global_flow}.

\vspace{15pt}

\noindent{\bf Step 4:} {\em  If  $z_0 \in \Sigma_{\mathbcal{e}_0}$ and  $\Phi^t_{\cX_h}(z_0)$ is defined for all $t \geq 0$, resp. $t \leq 0$, then $z_0 \in B(\Gamma^+)$, resp. $B(\Gamma^-)$.}
\begin{proof}[Proof of Step 4]
 We prove the statement with ``$+$''.     
It  suffices to show that
   \begin{equation}\label{ci.va}
       \operatorname{dist}(\zeta(t),\ K^+) \to 0 \quad \text{as} \quad t \to \infty\,.
   \end{equation}
By Step 2, this  holds for $z_0 \in \Gamma^+$ and we  restrict to the case $z_0 \notin \Gamma^+$.
To show \eqref{ci.va}, we
note that, since $\rho^2 = 2 h_0$ and $\Phi^t_{X_h}(z_0)$ is defined for all times,
$$
\tau_0(t) = \int_0^t \rho^{-2}(s) \di s \stackrel{\eqref{eq.nu}}{\geq } 
\int_0^t \frac{1}{2(h_0(0) + \nu s) } \di s \to \infty \mbox{ as } t \to \infty \ .
$$
Hence, since $w(\tau) \to K^+$ as  $\tau \to +\infty$, so does $\zeta(t)=w(\tau_0(t))$.

\end{proof}
\noindent{\bf Step 5:}   {\em If  $z_0 \in \Sigma_{\mathbcal{e}_0} \backslash(\Gamma^+ \cup \Gamma^-)$, then $\Phi^t_{\cX_h}(z_0)$ exists for all times $t \in \R$.}
 \begin{proof}[Proof of Step 5]
We show that if $z_0 \in \Sigma_{\mathbcal{e}_0} \setminus \Gamma^-$ then $\Phi^t_{\cX_h}(z_0)$ exists for any $t \geq 0$. 
Analogously one can prove that if $z_0 \in \Sigma_{\mathbcal{e}_0} \setminus \Gamma^+$ then $\Phi^t_{\cX_h}(z_0)$ exists for  any $t \leq 0$ concluding the proof. \\
To this aim, we prove that there exists a sequence of sets  $\{V_k\}_{k \in \N_0}\subset Z_{\mathbcal{e}_0}$ satisfying 
\begin{equation}\label{preso.tutto}
\Sigma_{\mathbcal{e}_0}\setminus \Gamma^- =  \bigcup_{k \in \N_0} \boldsymbol{\pi}^{-1} (V_k)\,,
\end{equation}
and we prove, inductively on $k$, that
\begin{equation}\label{well.def.k}
	\forall k \in \N_0\,, \quad \Phi_{\cX_h}^t(\boldsymbol{\pi}^{-1}(V_k)) \quad \text{is well defined }\  \forall t\geq 0\,.
\end{equation}
\underline{The set $V_0$.} 
To ensure \eqref{well.def.k} for $k=0$, it is sufficient to choose
$V_0 := U^+$,
with $U^+$ the set of Step 1.

\noindent\underline{The set $V_k$.} 
For $k \geq 1$, we recursively
 define 
\begin{equation}\label{V.k}
V_k:=\Phi^{-\underline{t}}_{\tilde{\cX}_h}(V_{k-1}),
\quad \mbox{ where }\quad  \underline{t} := \frac{1}{4\tilde{\nu}}\ln\left(\frac{3}{2}\right)\,, \quad \tilde{\nu}:= \max\left\lbrace 1,\  \nu\right\rbrace \ ,
\end{equation}
with $\nu$ defined in \eqref{hzerodot}. Suppose now that \eqref{well.def.k} holds for $k-1$. We are going to prove that \eqref{well.def.k} holds for $k$. Let $z_0 \in \boldsymbol{\pi}^{-1}(V_k)$ and suppose by contradiction that there exists $0<T_{\max}:= T_{\max}(z_0)<+\infty$ such that 
\begin{equation}\label{t.max}
\rho(t) \to 0, \ \mbox{ as } t \to T_{\max}^-\,.
\end{equation}
Then for all $t \in [0, T_{\max}/2]$, $\Phi_{X_h}^t(z_0)$ is well defined
and $\rho(t) \geq \underline{\rho}>0$ for some $\underline{\rho}>0$, thus on such time interval
one also has $\zeta(t) = w(\tau_0(t))$. We now claim that
\begin{equation}\label{lwb}
\tau_0\left(\frac{T_{\max}}{2}\right) \geq 2 \underline{t}\,.
\end{equation}
Indeed, passing to the limit $t\to T_{\max}^-$ in \eqref{eq.nu}, one has
\begin{equation}\label{Trex}
T_{\max} \geq \frac{h_0(z_0)}{\tilde{\nu}}\,.
\end{equation}
This, recalling that, again by \eqref{eq.nu},  $\rho^2(t) = 2 h_0(z(t)) \leq 2(h_0(z_0) + \tilde{\nu} t)$, gives
\begin{equation}
\begin{aligned}
\tau_0\left(\frac{T_{\max}}{2}\right)& = \int_0^{\frac{T_{\max}}{2}} \frac{1}{\rho^2(s)}\ \di s \stackrel{\eqref{Trex}}{\geq} \int_0^{\frac{h_0(z_0)}{2\tilde{\nu}}} \frac{1}{\rho^2(s)}\ \di s\\
& \geq \int_0^{\frac{h_0(z_0)}{2\tilde{\nu}}} \frac{1}{2(h_0(z_0) + \tilde{\nu} s)} \ \di s\\
& = \frac{1}{2 \tilde{\nu}} \int_{0}^{\frac{1}{2}} \frac{1}{1+s}\ \di s = \frac{1}{2\tilde{\nu}} \ln\left(\frac{3}{2}\right) = 2 \underline{t}\,,
\end{aligned}
\end{equation}
which proves the claim in \eqref{lwb}. Now, one observes that $\tau_0$ is an increasing continuous function on $[0, T_{\max}/2]$ satisfying
$\tau_0(0)=0$ and $\tau_0\left(\frac{T_{\max}}{2}\right)> \underline{t}$, thus $\exists!\ t_\star \in  \left(0, T_{\max}/2\right)$ such that $\tau_0(t_\star) = \underline{t}$.
But then, recalling $\zeta(0) = w(0) \in V_{k}:= \Phi^{-\underline{t}}_{\tilde X_h}(V_{k-1})$, one has
$$
\zeta(t_\star) = w(\tau_0(t_\star)) = w(\underline{t}) \in V_{k-1}\,,
$$
namely $\Phi^{t_\star}_{X_h}(z_0) \in \pi^{-1}(V_{k-1})$.
Then, for any $t \geq t_\star$, one has
$$
\Phi^{t}_{\cX_{h}}(z_0) = \Phi^{t- t_\star}_{\cX_h}(\Phi^{t_\star}_{\cX_h}(z_0))\in \Phi^{t-t_\star}_{\cX_{h}}\left(\boldsymbol{\pi}^{-1}\left(V_{k-1}\right)\right)\,,
$$
which is well defined by the inductive hypothesis. This then implies that $T_{\max}(z_0)= +\infty$, which gives the contradiction.\\
It remains to show that \eqref{preso.tutto} holds; but this is a consequence of the fact that $K^+$ is a global attractor, thus for any $z_0 \in \Sigma_{\mathbcal{e}_0} \setminus \Gamma^-$, there exists ${k}_0 \in \N_0$ such that $\Phi^{k_0 \underline{t}}_{\tilde X_h}(\zeta_0) \in U^+$, namely $z_0 \in \boldsymbol{\pi}^{-1} (V_{k_0})\,.$
\end{proof}
\noindent Steps 4 and 5 prove Item $(iv)$ of Proposition \ref{global_flow}.

\section{Proof of Lemma  \ref{prop:V-}}\label{app:V-}
 The idea is to construct $V^-$ as the sublevel of a smooth Lyapunov function defined on $Z_{\mathbcal{e}_0}$ for the flow of $\tilde \cX_h$. 
 The first step is to construct  a continuous  Lyapunov function, for which we invoke the following classical result, see e.g.  \cite[Proposition 1.4.21]{fisher2019hyperbolic}.

\begin{proposition}
There exists a continuous Lyapunov function $\cL: Z_{\mathbcal{e}_0} \to [0,1]$ associated to the flow of $\tilde{X}_h$, such that $\cL^{-1}(\{0\}) = K^+$, $\cL^{-1}(\{1\}) =K^-$, and $\cL$ is strictly decreasing along orbits lying outside $K^+ \cup K^-$.
\label{prop:LyapC0}
\end{proposition}
The next step is to  regularize the Lyapunov function. 
We shall exploit a result by Fathi-Pageault  \cite{Fathi2018SmoothingLF} 
that  adapts to our case. 
First we need the following definitions:

\begin{definition}
Given a Lyapunov function $\cL: M \to \R$ for a flow $\phi^t$ on  a compact manifold $M$, the 
 {\em neutral set} $N(\cL)$ is defined as
\[
N(\cL):=\{x \in M: \exists \ t>0 : \cL(\phi^t(x))=\cL(x) \}. 
\]
\end{definition}

\begin{theorem}{\cite[Theorem 2.1]{Fathi2018SmoothingLF}}
Let $M$ be a smooth  manifold, and let $\phi^t$ be a continuous flow on $M$. 
Assume that 
\begin{itemize}
    \item[(i)] the vector field $X(x):=\frac{\di}{\di t}\phi^t(x)|_{t=0}$ exists everywhere, is continuous and is a Lipschitz vector field; 
    \item[(ii)] the flow $\phi^t$ admits a 
continuous Lyapunov function
 $\cL: M \to \R$ such that $\cL(N(\cL))$ is contained in a closed subset of $\R$ with empty interior.
\end{itemize}
 Then we can approximate $\cL$ in the compact open topology of $C^0(M, \R)$ by a $C^{\infty}$ Lyapunov function $\tilde{\cL}: M \to \R$, arbitrarily close to $\cL$ and such that $N(\tilde{\cL})=Crit(\tilde{\cL})$, the set of critical points of $\tilde{\cL}$, and $X[\tilde{\cL}] <0$ on $M\backslash N(\tilde{\cL})$. 
\label{smoothlyap}
\end{theorem}
We apply this result to the continuous Lyapunov function $\cL\colon Z_{\mathbcal{e}_0} \to \R$ constructed in Proposition \ref{prop:LyapC0}. Clearly item $(i)$ of Theorem \ref{smoothlyap} is fulfilled since the vector field $\tilde \cX_h$ is smooth. 
To check $(ii)$, note that 
$N(\cL) = K^+ \cup K^-$, and again by Proposition \ref{prop:LyapC0} we have $\cL(N(\cL))=\{0, 1\}$ which is  a closed subset of $\R$ with empty interior. 
Hence Theorem \ref{smoothlyap} yields a smooth Lyapunov function
$  \tilde{\cL}: Z_{\mathbcal{e}_0} \to \R$
 such that 
\begin{equation}
    ||\tilde{\cL}-\cL||_{C^0}<\frac{1}{10} \ . 
    \label{distC_0}
\end{equation}
Note that  $\tilde{\cL}(Z_{\mathbcal{e}_0}) \subseteq  [-\frac{1}{10}, \frac{11}{10}]$.

\begin{proof}[Proof of Lemma \ref{prop:V-}]
Since $\cL^{-1}(1)=K^-$, we deduce from \eqref{distC_0} that  $K^- \subset \tilde{\cL}^{-1}([\frac{9}{10}, \frac{11}{10}])$, and analogously $K^+ \subset \tilde{\cL}^{-1}([-\frac{1}{10}, \frac{1}{10}])$. 
Then choose  $\varepsilon_1 \in (\frac{1}{10}, \frac{9}{10})$ which is not a critical value of $\tilde{\cL}$; such number exists by Sard's theorem being $\tilde{\cL}$  non-constant. 
Put  $V^-:=\tilde{\cL}^{-1}((\varepsilon_1, \frac{11}{10}])$; 
then  $V^-$ contains $K^-$ and $(V^-)^c$ contains $K^+$, $V^-$ is  invariant for negative times being the superlevel of a Lyapunov function, and its boundary $\partial V^- = \tilde{\cL}^{-1}(\epsilon_1)$ is transverse to $\tilde{X}_h$ since, by Theorem \ref{smoothlyap}, $\tilde{X}_h [\tilde{\cL}] <0$.
\end{proof}


\small

\bibliographystyle{plain}

\begin{thebibliography}{10}


\bibitem{BGMRgrowth}
D. Bambusi, B. Gr{\'e}bert, A.~Maspero, and D.~Robert.
\newblock Growth of {S}obolev norms for abstract linear {S}chr\"{o}dinger equations.
\newblock {\em JEMS,} 23(2):557--583, 2021.

\bibitem{BGMR}
D.~Bambusi, B.~Gr{\'e}bert, A.~Maspero, and D.~Robert.
\newblock Reducibility of the quantum harmonic oscillator in d-dimensions with
  polynomial time-dependent perturbation.
\newblock {\em Analysis {\&} {PDE}}, 11(3):775--799, 2018.

\bibitem{BGMRV}
D. Bambusi, B. Grébert, A. Maspero, D. Robert, and C. Villegas-Blas.
\newblock Longtime dynamics for the Landau Hamiltonian with a time dependent magnetic field.
\newblock  {\em arXiv}, arXiv:2402.00428, 2024.

\bibitem{QN}
D. Bambusi and B. Langella.
\newblock Growth of Sobolev norms in quasi
integrable quantum systems.
\newblock Accepted for publication on {\em Ann. Sci. Éc. Norm. Supér. arXiv,} arXiv:2202.04505, 2022.

\bibitem{QNbari}
D. Bambusi and B. Langella.
\newblock Globally integrable quantum systems and their perturbations.
\newblock {\em accepted for publication on “Singularities, Asymptotics and Limiting Models"
INdAM Springer Volume, pp. 64-103(2024).} arXiv:2403.18670, 2024.

\bibitem{BLMgrowth}
D. Bambusi, B. Langella and R. Montalto.
\newblock Growth of {S}obolev norms for unbounded perturbations of the {S}chr\"{o}dinger equation on flat tori.
\newblock {\em J. Differential Equations}, 318:344--358, 2022.

\bibitem{BLR2022}
D. Bambusi, B. Langella and M. Rouveyrol.
\newblock On the stable eigenvalues of perturbed anharmonic oscillators in dimension two.
\newblock {\em Comm. Math. Phys.}, 390(1): 309--348, 2022.

\bibitem{BanVega2022}
	V. Banica and L. Vega.
	\newblock Unbounded growth of the energy density associated to the Schr\"{o}dinger map and the binormal flow.
	\newblock {\em Annales de l’Institut Henri Poincaré C,  Analyse non linéaire}, 39:(4) 927–-946, 2022.

\bibitem{BertiMaspero}
M. Berti and A. Maspero.
\newblock Long time dynamics of {S}chr\"{o}dinger and wave equations on flat tori.
\newblock {\em J. Differential Equations}, 267(2):1167--1200, 2019.


\bibitem{BourgainGSN1999}
J.~Bourgain.
\newblock On growth of {S}obolev norms in linear {S}chr\"odinger equations with
  smooth time dependent potential.
\newblock {\em Journal d'Analyse  Mathématique}, 77(1):315--348, 1999.

\bibitem{Chabert}
A. Chabert.
\newblock Weakly turbulent solution to Schrödinger equation on the two-dimensional torus with real potential decaying at infinity.
\newblock{\em arXiv}, 	arXiv:2305.15939, 2024.


\bibitem{Chabert2}
A. Chabert.
\newblock A Weakly Turbulent solution to the cubic Nonlinear Harmonic Oscillator on $\R^2$ perturbed by a real smooth potential decaying to zero at infinity.
\newblock{\em Comm PDE}, 49(3): 185-216, 2024.

\bibitem{CKSTT}
J.~Colliander, M.~Keel, G.~Staffilani, H.~Takaoka, and T.~Tao.
\newblock Transfer of energy to high frequencies in the cubic defocusing
  nonlinear {S}chr\"odinger equation.
\newblock {\em Inventiones Mathematicae}, 181(1):39--113, 2010.



\bibitem{Davies.functional.calculus}
E.~B. Davies.
\newblock {The Functional Calculus}.
\newblock {\em J.  London Math.  Soc.}, 52(1):166--176, 1995.

\bibitem{Colin_de_Verdi_re_2020}
Y.~C. de~Verdi{\`{e}}re.
\newblock Spectral theory of pseudodifferential operators of degree 0 and an application to forced linear waves.
\newblock {\em Analysis {\&} {PDE}}, 13(5):1521--1537, 2020.

\bibitem{CdVSR}
Y. C. de Verdière and L. Saint-Raymond.
\newblock Attractors for two‐dimensional waves with homogeneous hamiltonians of degree 0.
\newblock {\em CPAM}, 73(2):421--462, 2020.

\bibitem{Delort2010} J.-M. Delort.
\newblock Growth of {S}obolev norms of solutions of linear {S}chr\"{o}dinger equations on some compact manifolds.
\newblock {\em Int. Math. Res. Not. IMRN}, 2010(2):2305--2328, 2010.

\bibitem{del}
J.-M. Delort.
\newblock Growth of {S}obolev norms for solutions of time dependent
  {S}chr\"odinger operators with harmonic oscillator potential.
\newblock {\em Comm PDE}, 39(1):1--33, 2014.

\bibitem{DuclosLevStovicek} P. Duclos, O. Lev and P.Stovicek.
\newblock On the energy growth of some periodically driven quantum systems with shrinking gaps in the spectrum.
\newblock {\em J. Stat. Phys.}, 130(1):169--193, 2008.

\bibitem{FaouRaphael}
E. Faou and P. Raphaël.
\newblock On weakly turbulent solutions to the perturbed linear harmonic oscillator.
\newblock {\em American Journal of Mathematics}, 145(5): 1465--1507, 2023.

\bibitem{Fathi2018SmoothingLF}
A. Fathi and P. Pageault.
\newblock Smoothing Lyapunov functions.
\newblock {\em Trans AMS}, 2018.

\bibitem{fisher2019hyperbolic}
T.~Fisher and B.~Hasselblatt.
\newblock {\em Hyperbolic Flows.}
\newblock Zurich lectures in advanced mathematics. European Mathematical Society, 2019.


\bibitem{geiges2008introduction}
H.~Geiges.
\newblock {\em An Introduction to Contact Topology}.
\newblock Cambridge Studies in Advanced Mathematics. Cambridge University Press, 2008.

\bibitem{gerard_grellier}
P.~G\'erard and S.~Grellier.
\newblock The cubic {S}zeg{\H o} equation and {H}ankel operators.
\newblock {\em Ast\'erisque}, (389):vi+112, 2017.


\bibitem{GL}
P.~G\'erard and  E. Lenzmann.
\newblock{The Calogero--Moser Derivative Nonlinear Schr\"odinger Equation}.
\newblock{\em CPAM}, 77(10):4008-4062, 2024.

\bibitem{Giroux1991}
E.~Giroux.
\newblock Convexité en topologie de contact.
\newblock {\em Commentarii mathematici Helvetici}, 66(4):637--677, 1991.

\bibitem{GG}
F. Giuliani and M. Guardia.
\newblock{Sobolev norms explosion for the cubic NLS on irrational tori}.
\newblock{\em Nonlinear Analysis}, 220,  2022. DOI:10.1016/j.na.2022.112865.

\bibitem{GGG}
F. Giuliani, and M. Guardia.
\newblock {Arnold diffusion in Hamiltonian systems on infinite lattices}.
\newblock{\em Commun. Pure Appl. Math.}, 2023. DOI 10.1002/cpa.22191.

\bibitem{GrePat}
B. Gr\'{e}bert, E. Paturel.
\newblock On reducibility of quantum harmonic oscillator on {$\Bbb R^d$} with quasiperiodic in time potential.
\newblock {\em Ann. Fac. Sci. Toulouse Math. (6)}, 28(5):977--1014, 2019.

\bibitem{guardia_haus_procesi16}
M.~Guardia, E.~Haus, and M.~Procesi.
\newblock Growth of {S}obolev norms for the analytic {NLS} on {$\Bbb{T}^2$}.
\newblock {\em Advances in Mathematics}, 301:615--692, 2016.

\bibitem{GHMMZ}
M. Guardia, E. Haus, Z. Hani, A Maspero, and M. Procesi. 
\newblock Strong nonlinear instability and growth of
Sobolev norms near quasiperiodic finite-gap tori for the 2D cubic NLS equation.
\newblock {\em JEMS}, 25(4):1497–1551, 2022. 

\bibitem{guardia_kaloshin}
M. Guardia and V. Kaloshin. 
\newblock Growth of Sobolev norms in the cubic defocusing nonlinear Schr\"odinger 
equation. 
\newblock{\em JEMS}, 17(1):71--149, 2015.



\bibitem{hani14}
Z.~Hani.
\newblock Long-time instability and unbounded {S}obolev orbits for some
  periodic nonlinear {S}chr\"odinger equations.
\newblock {\em ARMA}, 211(3):929--964, 2014.

\bibitem{hani15}
Z.~Hani, B.~Pausader, N.~Tzvetkov, and N.~Visciglia.
\newblock Modified scattering for the cubic {S}chr\"odinger equation on product
  spaces and applications.
\newblock {\em Forum of Mathematics, Pi}, 3:e4, 63, 2015.

\bibitem{HausMaspero}
E. Haus and  A. Maspero.
\newblock{ Growth of Sobolev norms in time dependent semiclassical anharmonic oscillators}. 
\newblock{\em J.  Func.  Anal.} , 278(2), 108316, 2020.

\bibitem{haus_procesi15}
E.~Haus and M.~Procesi.
\newblock Growth of {S}obolev norms for the quintic {NLS} on {$\T^2$}.
\newblock {\em Anal. {PDE}}, 8(4):883--922, 2015.

\bibitem{Hurley}
M. Hurley.
\newblock Attractors: Persistence, and density of their basins.
\newblock {\em Trans AMS}, 269(1):247--271, 1982.

\bibitem{Kuk}
S. Kuksin. 
\newblock{Growth and oscillations of solutions of nonlinear Schr\"odinger equation}.
\newblock{\em CMP,}178(2):265--280, 1996.

\bibitem{Kuk2}
S. Kuksin.  
\newblock{On turbulence in nonlinear Schr\"odinger equations.}
\newblock{\em GAFA,} 7(4):783--822, 1997.

\bibitem{lee2003introduction}
J.M. Lee.
\newblock {\em Introduction to Smooth Manifolds}.
\newblock Graduate Texts in Mathematics. Springer, 2003.

\bibitem{LLZ2}
Z. Liang, J. Luo, and Z. Zhao. 
\newblock{Symplectic Normal Form and Growth of Sobolev Norm. }
\newblock{\em arXiv,}
arXiv:2312.16492, 2023.

\bibitem{LLZ1}
Z. Liang, Z. Zhao and Q. Zhou.
\newblock{1-d quantum harmonic oscillator with time quasi-periodic quadratic perturbation: reducibility and growth of Sobolev norms. }
\newblock{\em J.   Math.  Pures  Appl.} 146(1), 158--182 (2021). 



\bibitem{LLZ}
Z. Liang, Z. Zhao and Q. Zhou. 
\newblock{Almost reducibility and oscillatory growth of Sobolev norms.}
\newblock{\em Advances in Mathematics} 436, 109417 (2024). 


\bibitem{LLZ3}
 J. Luo, Z. Liang and Z. Zhao.
 \newblock{Growth of Sobolev Norms in 1-d Quantum Harmonic Oscillator with Polynomial Time Quasi-periodic Perturbation.}
 \newblock{\em Commun. Math. Phys.} 392, 1--23, 2022.

\bibitem{Maspero2018LowerBO}
A. Maspero.
\newblock Lower bounds on the growth of Sobolev norms in some linear time dependent Schr{\"o}dinger equations.
\newblock {\em MRL}, 26(4):1197--1215, 2019.

\bibitem{Mas22}
A.~Maspero.
\newblock Growth of {S}obolev norms in linear {S}chr\"odinger equations as a dispersive phenomenon.
\newblock {\em Advances in Mathematics}, 411:Paper No. 108800, 48, 2022.

\bibitem{Mas23}
A.~Maspero.
\newblock Generic transporters for the linear time-dependent quantum harmonic oscillator on {$\Bbb{R}$}.
\newblock {\em IMRN}, (14):12088--12118, 2023.

\bibitem{MasperoMurgante}
A. Maspero, and F. Murgante.
\newblock One dimensional energy cascades in a fractional quasilinear NLS.
\newblock {\em arXiv}, 	arXiv:2408.01097,  2024.

\bibitem{Maspero_2017}
A.~Maspero and D.~Robert.
\newblock On time dependent Schrödinger equations: Global well-posedness and growth of Sobolev norms.
\newblock {\em J.  Func.  Anal. }, 273(2):721--781, 2017.

\bibitem{Nenciu} G. Nenciu.
\newblock Adiabatic theory: stability of systems with increasing gaps.
\newblock {\em Ann. Inst. Poincar\'{e}}, 67(4):411--424, 1997.

\bibitem{palis1998geometric}
J.~Palis and W.~de~Melo.
\newblock {\em Geometric Theory of Dynamical Systems: An Introduction}.
\newblock Selected monographies. Collaege Press, University of Beijing, 1998.

\bibitem{Peixoto1962StructuralSO}
M.~M. Peixoto.
\newblock Structural stability on two-dimensional manifolds.
\newblock {\em Topology}, 1:101--120, 1962.

\bibitem{shubin2001pseudodifferential}
M.A. Shubin.
\newblock {\em Pseudodifferential Operators and Spectral Theory}.
\newblock Springer Berlin Heidelberg, 2001.


\bibitem{Thomann2020GrowthOS}
L. Thomann.
\newblock Growth of Sobolev norms for linear Schr{\"o}dinger operators.
\newblock {\em Ann. H. Lebesgue}, 4:1595-1618, 2021.


\end{thebibliography}

\end{document}